\documentclass[12pt]{amsart}

\usepackage[dvipsnames]{xcolor}
\usepackage{algorithmic}
\usepackage{amsbsy}
\usepackage{amsfonts}
\usepackage{amsmath}
\usepackage{amssymb}
\usepackage{amsthm}
\usepackage{anyfontsize}
\usepackage[title]{appendix}
\usepackage[font = footnotesize]{caption}
\usepackage{enumitem}
\usepackage{float}
\usepackage[T1]{fontenc}
\usepackage{graphicx}
\usepackage{hyperref}
\usepackage[noabbrev]{cleveref}
\usepackage[numbers,sort&compress]{natbib}
\usepackage{nicematrix}
\usepackage{pgfplots}
\usepackage{rotating}
\usepackage[font=footnotesize]{subcaption}
\usepackage{tikz}
\usepackage{verbatim}
\usepackage{wavelet}
\usepackage{xparse}
\usepackage[top=0.54in,bottom=0.6in,left=0.8in,right=0.8in]{geometry}

\Crefname{thm}{Theorem}{Theorems}
\Crefname{cor}{Corollary}{Corollaries}
\Crefname{prop}{Proposition}{Propositions}
\Crefname{example}{Example}{Examples}
\Crefname{lemma}{Lemma}{Lemmas}

\newcommand {\ii} {\mathbf{i}} %for imaginary unit

\newcommand{\cut}[1]{ \textcolor{red}{} } %cut or not
\renewcommand{\SS}{\mathcal{S}} %stencil
\newcommand{\spt}{\mathbf{c}^*} %stencil center
\newcommand{\bpt}{\mathbf{b}^*} %base points for Taylor expansion
\newcommand{\tbpt}{\tilde{\mathbf{b}}^*} %another base points for Taylor expansion
\newcommand{\gx}{\beta}
\newcommand{\gy}{\gamma}
\newcommand{\sv}{\mathbf{s}} %shifting vector
\newcommand{\psv}{p^{\sv}}
\newcommand{\RN}{\square} %rectangle indx set in N_0^2
\newcommand{\LN}{\mathsf{\Gamma}}
\newcommand{\nv}{\mathbf{n}} %unit normal vector
\newcommand{\hp}[1]{H_{#1}} % left half plane of a line

\DeclareDocumentCommand {\dd} {m o}
    { \IfValueTF{#2}
    { \frac{\mathrm{d}^{#2}} {\mathrm{d} #1^{#2}} }
    { \frac{\mathrm{d}} {\mathrm{d} #1} }
}

\newtheorem{thm}{Theorem}[section]
\newtheorem{remark}[thm]{Remark}
\newtheorem{prop}[thm]{Proposition}
\newtheorem{lemma}[thm]{Lemma}
\newtheorem{cor}[thm]{Corollary}
\newtheorem{example}[thm]{Example}
\newtheorem{definition}[thm]{Definition}
\numberwithin{equation}{section}

\allowdisplaybreaks

\hypersetup {
    colorlinks = true,
    linkcolor = black,
    filecolor = magenta,
    urlcolor = NavyBlue,
    citecolor = NavyBlue
}

\begin{document}

\title[Convergent Sixth-order FDMs for Elliptic PDEs]{Convergent Sixth-order Compact Finite Difference Method for Variable-Coefficient Elliptic PDEs in Curved Domains}

\author{Bin Han and Jiwoon Sim}

\address{Department of Mathematical and Statistical Sciences, University of Alberta, Edmonton, Alberta, Canada T6G 2G1.}
\email{bhan@ualberta.ca, jiwoon2@ualberta.ca}

\thanks{Research supported in part by
Natural Sciences and Engineering Research Council (NSERC) of Canada under grant RGPIN-2024-04991}

\makeatletter \@addtoreset{equation}{section} \makeatother

\begin{abstract}
    Finite difference methods (FDMs) are widely used for solving partial differential equations (PDEs) due to their relatively simple implementation. However, they face significant challenges when applied to non-rectangular domains and in establishing theoretical convergence, particularly for high-order schemes. In this paper, we focus on solving the elliptic equation $-\nabla \cdot (a\nabla u)=f$ in a two-dimensional curved domain $\Omega$, where the diffusion coefficient $a$ is variable and smooth. We propose a sixth-order $9$-point compact FDM on uniform Cartesian grids within the domain, not relying on ghost points or information outside $\overline{\Omega}$. All the boundary stencils near $\partial \Omega$ have at most $6$ different configurations and use at most $8$ grid points inside $\Omega$. We rigorously establish the sixth-order convergence of the numerically approximated solution $u_h$ in the $\infty$-norm. Additionally, we derive a gradient approximation $\nabla u$ directly from $u_h$ without solving auxiliary equations. This gradient approximation achieves proven accuracy of order $5+\frac{1}{q}$ in the $q$-norm for all $1\le q\le \infty$ (with a logarithmic factor $\log h$ for $1\le q<2$). To validate our proposed sixth-order compact finite different method, we provide several numerical examples that illustrate the sixth-order accuracy and computational efficiency of both the numerical solution and the gradient approximation for solving elliptic PDEs in curved domains.
\end{abstract}


    \maketitle
    \pagenumbering{arabic}

    \section{Introduction}
    \label{sec: intro}

    The finite difference method (FDM) is a widely used tool for numerically solving partial differential equations, largely due to its simplicity and straightforward implementation on Cartesian grids. However, it faces significant challenges when applied to irregular domains with curved boundaries, particularly at grid points near the boundary (e.g., see \cite{jensen1972finite}). On the other hand, high-order FDM schemes are highly desired for their efficiency and high accuracy. However, high-order FDM schemes are considerably more difficult to construct with small stencils and proven theoretical convergence. This paper addresses these challenging issues by developing an efficient and reliable finite difference scheme tailored for variable-coefficient elliptic PDE in curved domains.

    In this paper, we consider the following boundary value problem:
    \begin{equation}\label{eq:PDE}
        \begin{cases}
            -\nabla \cdot (a \nabla u) = f &\mbox{in } \Omega, \\
            u= g &\mbox{on } \partial \Omega,
        \end{cases}
    \end{equation}
    where $\Omega \subset  \R^2$ is a bounded open domain with smooth boundary $\partial \Omega$, and the diffusion coefficient $a>0$ is a smooth function in $\Omega$. In this paper, we are particularly interested in high-order compact FDMs with small stencils and proven theoretical convergence for the above elliptic PDEs in curved domains with variable diffusion coefficient $a$. The precise assumptions on $a, f, \Omega$ for our developed schemes and proven theoretical convergence rates will be stated in \Cref{sec: convergence}.

    It is well known that higher-order FDMs necessarily require larger stencils. But FDMs with small stencils are of fundamental importance and interest in computational mathematics, because small stencils facilitate implementation, lead to small bandwidth and improved sparsity of the stiffness matrices, and more importantly, significantly reduce the number of exceptional boundary stencils with required modified stencil coefficients near the curved boundaries. As a consequence, compact FDMs (i.e., schemes having $1$-ring stencils) are highly sought in the literature of numerical PDEs. In this paper, we are only interested in $9$-point compact FDMs with the highest possible accuracy order for the elliptic PDE with variable coefficients in curved domains.

    Cartesian grids are particularly desired for the convenience of setting up FDMs and are well suited for rectangular/regular domain $\Omega$. Our developed sixth-order FDM shall use the grids generated from Cartesian grids, more precisely, for any given mesh size $h>0$ and any point $p\in \R^2$, we shall only use the grid $\Omega_h:=\Omega \cap (p+h\Z^2)$, where $p+h \Z^2:=\{p+(ih, jh) \setsp i,j\in \Z\}$. Without loss of generality, we shall always take $p=(0,0)$ for the purpose of simple presentation. That is, for any given mesh size $h>0$, we define the computational grids
    \be \label{eq: Omegah}
    \Omega_h:=\Omega \cap (h \Z^2)\quad \mbox{with}\quad
    \Omega_h^{\circ}:=\{ p\in \Omega_h \setsp p+[-h,h]^2 \subseteq \Omega\},\quad
    \partial \Omega_h:=\Omega_h \bs \Omega_h^{\circ},
    \ee
    where $\Omega_h^{\circ}$ is for interior stencils and $\partial \Omega_h$ is for boundary stencils near $\partial \Omega$. It is important to notice that grid points in $\partial \Omega_h$ are not lying on the boundary $\partial \Omega$ of the problem domain $\Omega$ but within at most $\sqrt{2}h$ distance to the boundary $\partial \Omega$. For each point $p\in \Omega_h^\circ$, we shall use a $9$-point compact stencil whose center is $p$. For each boundary point in $\partial \Omega_h$, we shall use no more than $8$-point stencils and we have no more than six special types of boundary stencils. Our proposed method achieves sixth-order consistency and never uses ghost points or information outside the closure of $\Omega$. In addition, we rigorously establish the sixth-order convergence of our proposed scheme by ensuring the discrete maximum principle. Furthermore, we derive a fifth-order accurate approximation of the gradient $\nabla u$ from the numerically approximated solution without solving additional equations.

    Because there is a huge literature on various finite difference methods, here we only review the literature related to the particular elliptic PDEs \eqref{eq:PDE} for a domain $\Omega$ to be either rectangular or curved. Because the proof of theoretical convergence of FDMs is often challenging, while we are reviewing the literature on FDMs for the elliptic PDE \eqref{eq:PDE}, we shall also discuss when their convergence has been established or not in the literature.
    Let $\Omega$ be a rectangular domain (or a cube in three-dimensional space). For the constant diffusion coefficient $a=1$,
    compact FDMs up to sixth order have been extensively studied and developed in \cite{li2023high, settle2013derivation, wang2009sixth, zhai2013family, zhai2014new, feng2021sixth} and many references therein. Higher consistency order is achieved with non-compact stencils, e.g., \cite{feng2020fft}. Now we review the literature for the diffusion coefficient $a$ to be a smooth function. Ma and Ge \cite{ma2020high} proposed blended compact difference schemes that have up to sixth-order consistency for 3D elliptic equations. Wang et al. \cite{wang2014fourth} constructed a fourth-order scheme for semilinear elliptic problems. The FDM proposed by Shi et al. \cite{shi2021high} reaches fourth-order accuracy for both the function $u$ and its gradient. For elliptic interface problems, Feng et al. \cite{FHM22} obtained a compact FDM with fourth-order accuracy of the solution and third-order accuracy of its gradient. The convergence is proven in \cite{FHM22,shi2021high,wang2014fourth}. Feng et al. \cite{feng2024sixth} provided sixth-order methods for \cref{eq:PDE} with interfaces. When no interface exists, the proposed method is proven to achieve sixth-order convergence in \cite{feng2024sixth}. According to the existing literature (e.g., \cite{feng2021sixth,feng2024sixth}), for a rectangular domain $\Omega$, six is the highest possible accuracy order for compact stencils.

    We now review the literature when $\Omega$ is a smooth curved domain. For $a = 1$, the classical approach is the Shortley-Weller method \cite{sw38}, where one directly modifies stencil coefficients for stencils near $\partial \Omega$. This method achieves convergent second-order accuracy. Bramble and Hubbard \cite{bramble1964new} and Price \cite{price1968monotone} proposed fourth-order FDMs for the Poisson equation and the convection-diffusion equation, respectively. These methods have a relatively small stencil and the convergence is proven. Esmaeilzadeh and Barron \cite{eb22} transformed each stencil near the boundary to the standard 5-point stencil and derived a fourth-order FDM. Pan et al. \cite{pan2021high} enlarged the computational domain and used the techniques of immersed interface method to derive third-order schemes. Using fictitious values formulation and ray-casting matched interface and boundary (MIB) method, \cite{RFZ22,Li2025spatially} proposed fourth-order FFT accelerated schemes which successfully handle sharply curved boundaries. The convergence of the last three methods is not established yet.

    There are much fewer papers in the literature addressing the case that the diffusion coefficient $a$ is smooth and $\Omega$ is a smooth curved domain. Samarskii and Fryazinov \cite{samarskii1971finite} proposed a second order convergent scheme using non-uniform mesh. Ito et al. \cite{ilk05} proposed FDMs with up to fourth-order consistency by approximating the solution near the boundary via polynomial interpolation. In \cite{gf2005,gfck2002}, the authors extrapolated the solution onto ghost cells to the other side of the boundary, which results in second order convergent and fourth-order consistent FDMs, respectively. A similar strategy is considered by Clain et al. \cite{clp21}, which is able to achieve arbitrary consistency order with large stencils. However, the convergence of the numerical solution is not proven in the above FDMs with consistency order higher than 2, and these methods employ large stencils to obtain a desired approximation to the solution near the boundary. As a consequence of using large stencils, one often has to consider many specially designed stencil configurations with modified coefficients near the boundary curves. Besides, the resulting linear system becomes much less sparse, leading to increased computational complexity and implementation difficulties of a FDM scheme with large stencils.

    The major contribution of this article is to provide a reliable scheme that is proven to have sixth-order convergence. The convergence of FDM is typically proved via the discrete maximum principle, which requires that the discretization of the differential operator is a monotone matrix \cite{varga1966discrete}. In practice, such a matrix is provided with a nonsingular M-matrix, or a weakly chained diagonally dominant matrix with nonpositive off-diagonal entries (see \cite{shivakumar1974sufficient,plemmons1977m} for the definition and equivalence of these matrices). However, as indicated in \cite{li2020monotonicity}, except for certain 9-point finite difference methods, almost all high-order schemes produced by finite difference or finite element methods do not result in an M-matrix due to positive off-diagonal entries. In the present paper, we ensure the monotone property by carefully constructing the stencil near the boundary. Based on the sixth-order convergence of the numerical solution, we derive a fifth-order approximation of the gradient $\nabla u$ in the $\infty$-norm without solving auxiliary equations. Furthermore, we observe that the numerical solution exhibits certain regularity, which enables us to prove a superconvergence of order $5+\frac{1}{q}$ in the $q$-norm for all $1\le q \le \infty$ (with a logarithmic factor $\log h$ for $1\le q<2$).

    The paper is organized as follows. In \Cref{sec: notations}, we introduce complex partial derivatives and discuss their property and advantages for solving \eqref{eq:PDE} in a smooth curved domain. The sixth-order $9$-point compact FDM at interior grid points is developed in \Cref{sec: interior}. In \Cref{sec: boundary} we construct the fourth-order FDM at boundary grid points with emphasis on small boundary stencils using at most $8$ grid points near $\partial \Omega$ and having at most $6$ different boundary stencil configurations. \Cref{sec: convergence} deals with the theoretical convergence analysis of our method for both the numerical solution and gradient approximation. The stencil coefficients of the proposed method consist of high-order derivatives of the functions in \cref{eq:PDE}. In \Cref{sec: derivatives}, we will provide an efficient way to evaluate these derivatives using only function values. For the rest of \Cref{sec: numerical} we provide some useful details to implement the proposed method and test it in diverse scenarios with oscillating functions and domain boundaries. Concluding remarks are given in \Cref{sec: conclusion}.

%%%%%%%%%%%%%%%%%%%%%%%%%%%%%%%%%%%%%%%%%%%%%%%%%%%%%%%%%%%%%%%%%%%%%%%
%%%%%%%%%%%%%%%%%%%%%%%%%%%%%%%%%%%%%%%%%%%%%%%%%%%%%%%%%%%%%%%%%%%%%%%
    \section{Auxiliary Results Using Complex Partial Derivatives for Constructing FDMs}
    \label{sec: notations}

    To present our construction of compact FDMs in later sections, it is very helpful for us to introduce some notations, necessary definitions, and auxiliary results here.

    To avoid complexity of presentation, in \Cref{sec: notations,sec: interior,sec: boundary} we assume that all involved functions are smooth enough; the formal assumptions are given in \Cref{sec: convergence}. Define $\N_0 := \N\cup\{0\}$. For a smooth function $v$ and $(k, \ell) \in \NN^2$, the ordinary partial derivative $\partial^{(k,\ell)} v$ and the so-called ``complex'' partial derivative $\partial^{(k,\ell)}_{\C} v$ are defined by
    \begin{equation}\label{deriv:C}
        \partial^{(k, \ell)} v := \frac{\partial^{k + \ell} v}{\partial^k x \partial^\ell y}
        \quad \mbox{and} \quad
        \partial_{\C}^{(k, \ell)} v := \frac{1}{2^{k+\ell}} \left( \frac{\partial}{\partial x} - \ii  \frac{\partial}{\partial y} \right)^k \left( \frac{\partial}{\partial x} + \ii  \frac{\partial}{\partial y} \right)^\ell v,
    \end{equation}
    where $\ii $ is the imaginary unit.
    To understand the definition \eqref{deriv:C},
    we shall see how the standard Taylor expansion can be equivalently expressed by using the complex partial derivatives $\partial^{(k,\ell)}_\C$. Throughout the paper, the notation $\bo(h^M)$ with various subscripts refers to a function that is bounded by $C h^M$ as $h \to 0^+$, where the constant $C$ only depends on the expressions and their derivatives in the subscript, and $C$ remains positive and bounded if its dependencies are bounded. For example, the remainder term in the standard Taylor expansion for a smooth function $v$ can be denoted as $\bo_{v} (h^n)$.

    \begin{prop}\label{prop:taylor:C}
        Let $v: \R^2 \to \R$ be a smooth function in a neighborhood of a base point $\bpt\in \R^2$. For any $n \in \mathbb{N}$, $p \in \R^2$ and  sufficiently small $h\in \R$, we have
        \begin{equation}\label{taylor:C}
            v(\bpt + ph) = \sum_{0 \leq k + \ell<n} \frac{1}{k! \ell!} (p_r + \ii  p_i)^k (p_r - \ii  p_i)^\ell h^{k+\ell} \partial_{\C}^{(k, \ell)} v(\bpt) + \bo_v(h^{n}),
        \end{equation}
        where $(p_r, p_i):=p$, i.e., we identify the point $p\in \R^2$ with the complex number $p_r+\ii p_i\in \C$.
    \end{prop}

    \begin{proof}
Consider the transform $z:=x+\ii y$ and $\bar{z}:=x-\ii y$. Then $x=\frac{1}{2}(z+\bar{z})$ and $y=\frac{1}{2\ii}(z-\bar{z})$. Using the transform, we can define a bivariate function $V(z,\bar{z}):=v(x,y)$. Noting that
\[
\frac{\partial}{\partial z}=
\frac{\partial}{\partial x} \frac{\partial x}{\partial z}+\frac{\partial}{\partial y} \frac{\partial y}{\partial z}=\frac{1}{2}\left(\frac{\partial}{\partial x}-\ii\frac{\partial}{\partial y}\right),
\qquad
\frac{\partial}{\partial \bar{z}}=
\frac{\partial}{\partial x} \frac{\partial x}{\partial \bar{z}}+\frac{\partial}{\partial y} \frac{\partial y}{\partial \bar{z}}=\frac{1}{2}\left(\frac{\partial}{\partial x}+\ii\frac{\partial}{\partial y}\right),
\]
we observe from the definition \eqref{deriv:C} that $\partial^{(k,\ell)}_\C v(x,y)=(\frac{\partial}{\partial z})^k (\frac{\partial }{\partial \bar{z}})^\ell v(x,y)=
(\frac{\partial}{\partial z})^k (\frac{\partial }{\partial \bar{z}})^\ell V(z,\bar{z})=
\partial^{(k,\ell)} V(z,\bar{z})$, which is just the standard $(k,\ell)$-th partial derivative of $V$.

Note that the standard Taylor expansion of $v(\bpt+ph)$ at the base point $\bpt$ is just the Taylor expansion of the one-dimensional function $v(\bpt+ph)$ of variable $h$ at the base point $h=0$. Similarly write $\bpt=(\bpt_r,\bpt_i)$ as in $p=(p_r,p_i)$. Note that $v(\bpt+ph)=V((\bpt_r+\ii \bpt_i)+(p_r+\ii p_i)h, (\bpt_r-\ii \bpt_i)+(p_r-\ii p_i)h)$, which can be regarded as a function of $h$ and whose Taylor expansion at the base point $h=0$ is just the right-hand side of \eqref{taylor:C}.
\end{proof}

    In sharp contrast to all papers in the literature on FDMs, in this paper we shall use complex partial derivatives $\partial^{(k,\ell)}_{\C}$ in \eqref{deriv:C}, which offer us a different perspective and a key advantage of symmetry over our previous approach in \cite{feng2021sixth,FHM22,feng2024sixth}.
    To develop FDMs for the elliptic equation \eqref{eq:PDE} in curved domains, this approach using complex partial derivatives is necessary and critical for us to avoid complicated expressions arising from geometries of curved boundaries for building finite difference schemes at boundary stencils.
    It is also very important to keep in mind that even though complex numbers will appear in our construction, the coefficients in all our constructed FDM schemes through complex partial derivatives are real numbers (see
    \Cref{sec: interior,sec: boundary}).

    Throughout the paper, for simplicity of presentation, we often drop the base point $\bpt$ in $\partial^{(k,\ell}_\C v(\bpt)$ of \eqref{taylor:C} in \Cref{prop:taylor:C}
    if the base point is clear in the context.
    For any $M\in \NN$ and a smooth function $u:\R^2\rightarrow \C$, now applying \Cref{prop:taylor:C} with $v=u$ and $n=M+2$, we have the following Taylor expansion at a base point $\bpt$:
    \begin{equation}
        \label{u:taylor0}
        u(\bpt+ p h)= \sum_{0 \leq k + \ell \leq M + 1}
        N^{k, \ell}(p)
        h^{k + \ell} \partial_{\C}^{(k, \ell)} u + \bo_u(h^{M + 2}),
    \end{equation}
    where we omitted the base point $\bpt$ after the function $u$ for simplicity, and we define
    \begin{equation} \label{def:N}
        N^{k,\ell}(p):=N^{k,\ell}(p_r, p_i):= \frac{(p_r + \ii  p_i)^k (p_r - \ii  p_i)^\ell}{k! \ell!}\quad \mbox{with}\quad (p_r, p_i):=p\in \R^2.
    \end{equation}
    We now study the Taylor expansion of a smooth exact solution $u$ of the model problem \eqref{eq:PDE} by using complex partial derivatives $\partial^{(k,\ell)}_\C u$ for $(k,\ell)\in \NN^2$. To make our presentation simpler, we note that the model problem \eqref{eq:PDE} can be simply rewritten as follows:
    \begin{equation}\label{pde:mod}
        \Delta u = \nabla \tilde{a} \cdot \nabla u + \tilde{f} \quad \mbox{ with }\quad \tilde{a} := -\ln a \quad \mbox{ and }\quad \tilde f := -\frac{f}{a}.
    \end{equation}
    Using complex partial derivatives, the above equation \eqref{pde:mod} can be equivalently transformed into
    \begin{equation}\label{pde:11}
        \partial_{\C}^{(1, 1)} u = \frac{1}{2} \partial_{\C}^{(0, 1)} \tilde{a} \, \partial_{\C}^{(1, 0)} u + \frac{1}{2} \partial_{\C}^{(1, 0)} \tilde{a} \, \partial_{\C}^{(0, 1)} u + \frac{1}{4}\tilde{f}.
    \end{equation}
    Taking complex partial derivatives to both sides of \eqref{pde:11}, for $k, \ell \geq 1$, we deduce that
    {\footnotesize
    \begin{align*}
        \partial_{\C}^{(k, \ell)} u
        &= \frac{1}{2} \sum_{\substack{0 \leq m \leq k - 1 \\ 0 \leq n \leq \ell - 1}} \binom{k - 1}{m} \binom{\ell - 1}{n} \left( \partial_{\C}^{(k - 1 - m, \ell - n)} \tilde{a} \, \partial_{\C}^{(m + 1, n)} u + \partial_{\C}^{(k - m, \ell - 1 - n)} \tilde{a} \, \partial_{\C}^{(m, n + 1)} u \right) + \frac{1}{4} \partial_{\C}^{(k - 1, \ell - 1)} \tilde{f}.
    \end{align*}
    }%
    Taking into account of the identity \eqref{pde:11},
    we shall define two index subsets of $\NN^2$ as follows:
    \begin{equation}\label{rectangle:ind}
        \RN^{k,\ell}_{m, n}:=
        \{(i, j) \in \NN^2 \setsp m \leq i \leq k, n \leq j\leq \ell, (i, j) \neq (k, \ell) \},
    \end{equation}
    i.e., the index set $\RN^{k,\ell}_{m, n}$ is the rectangle $[m,k]\times [n,\ell]$ in $\N_0^2$ but without the corner $(k,\ell)$, and we define an index subset $\LN^k_\ell$ of $\N_0^2$  (with points only sitting on the nonnegative $x$-axis or $y$-axis) by
    \begin{equation*}
        \LN^{k}_{\ell}
        := \{ (m, 0) \in \NN^2 \setsp m=0,\ldots, k \} \cup \{ (0, n) \in \NN^2 \setsp n=1,\ldots, \ell \}
    \end{equation*}
    for $k, \ell, m, n\in \NN$.
    One can check that the above expression of $\partial^{(k,\ell)}_\C u$ can be simplified into
    \begin{equation}
        \label{eq: partial_C^{k, l} 1}
        \partial_{\C}^{(k, \ell)} u
        = \sum_{(m, n) \in
        %R^*_{0, k, 0, l}
        \RN^{k,\ell}_{0,0}}
        \tilde{a}^{k, \ell}_{m, n}
        \partial_{\C}^{(m, n)} u + \frac{1}{4} \partial_{\C}^{(k - 1, \ell - 1)} \tilde{f},
    \end{equation}
    where
    \begin{equation}\label{ta}
        \tilde{a}^{k, \ell}_{m, n} := \left( \frac{m}{2k} + \frac{n}{2\ell} - \frac{mn}{k\ell} \right) \binom{k}{m} \binom{\ell}{n} \partial_{\C}^{(k - m, \ell - n)} \tilde{a}.
    \end{equation}

    Let $\td$ be the sequence such that
    \begin{equation}\label{delta}
    \td(0):=1 \quad  \mbox{and}\quad \td(k):=0\quad \mbox{ for } k\ne 0.
    \end{equation}
    The identity \eqref{eq: partial_C^{k, l} 1} implies that $\partial_{\C}^{(k, \ell)} u$ can be eventually represented in terms of $\partial_{\C}^{(m, n)} u$ for $(m, n) \in \LN^{k}_{\ell}$. More precisely,
    \begin{equation}
        \label{eq: partial_C^{k, l} 2}
        \partial_{\C}^{(k, \ell)} u
        = \sum_{(m, n) \in \LN^{k}_{\ell}} \tilde{A}^{k,\ell}_{m,n}
        \partial_{\C}^{(m, n)} u + \tilde{F}_{k, \ell}
    \end{equation}
    for uniquely determined coefficients $\tilde{A}^{k, \ell}_{m, n}$ and  $\tilde{F}_{k, \ell}$ defined through the following recursive formulas:
    \begin{equation}\label{tA0}
        \tilde{A}^{k, \ell}_{m, n} := \td(k - m) \td(\ell - n),
        \qquad
        \tilde{F}_{k, \ell} := 0\qquad \mbox{ if }\; k \ell=0,
    \end{equation}
    where $k,\ell,m,n\in \mathbb{N}_0$,
    and the other values for $k\ell\ne 0$ are recursively defined through
    \begin{equation}
        \label{tA}
        \tilde{A}^{k,\ell}_{m, n} := \sum_{(i, j) \in \RN^{k,\ell}_{m, n}} \tilde{a}^{k, \ell}_{i, j} \tilde{A}^{i, j}_{m, n},
        \qquad
        \tilde{F}_{k, \ell} := \sum_{(i, j) \in \RN^{k,\ell}_{0, 0}} \tilde{a}^{k, \ell}_{i, j} \tilde{F}_{i, j} + \frac{1}{4}\partial_{\C}^{(k - 1, \ell - 1)} \tilde{f}.
    \end{equation}

    Therefore, using the identity \eqref{eq: partial_C^{k, l} 2}, we can reformulate the Taylor expansion in \eqref{u:taylor0} of the solution $u$ to the model problem \eqref{pde:mod}
    at a base point $\bpt\in \Omega$ as follows:
    \begin{equation}
        \label{u:taylor}
        u(\bpt+ ph)
        = \sum_{(m, n) \in \LN^{M+1}_{M + 1}} \sum_{k = m + n}^{M + 1}
        A^k_{m,n}(p) \partial_{\C}^{(m, n)} u(\bpt)  h^k + F(p) + \bo_u(h^{M + 2}),
    \end{equation}
    for $p\in \R^2$ with the line segment $[\bpt,\bpt+p h]$ inside $\Omega$, where $A^k_{m,n}(p)$ and $F(p)$ are defined below:
    \begin{equation}\label{Akmnp}
    A^k_{m,n}(p):= \sum_{j = m}^{k - n} N^{j,k-j}(p) \tilde{A}^{j, k - j}_{m,n} \quad \mbox{ for  } (m,n)\in \LN^{M+1}_{M+1}\quad \mbox{and}\quad k=0,\ldots,M+1,
    \end{equation}
    with the convention $A^k_{m,n}(p):=0$ for $k=0,\ldots, m+n-1$ because $\sum_{j=m}^{k-n}$ is empty, and
    \begin{equation}\label{Fp}
    F(p):= \sum_{0 \leq k + \ell \leq M + 1}
    N^{k,\ell}(p) \tilde{F}_{k, \ell} h^{k + \ell}.
    \end{equation}

    We finish this section by making some remarks. By equations \eqref{def:N}, \eqref{tA}, \eqref{Akmnp} and $\tilde{a}^{k, l}_{0, 0} = 0$ for $k, l \geq 1$ in \eqref{ta}, we get $A^k_{0, 0} = \td(k)$ for $k \in \NN$. By \cref{tA0,def:N,Akmnp}, we have
    \begin{equation}
        \label{Amm0}
        A^{m}_{m,0}(p) =
        N^{m,0}(p) \tilde{A}^{m, 0}_{m, 0} =
        N^{m,0}(p)= \frac{(p_r+\ii p_i)^m}{m!}, \quad
        A^m_{0,m}(p) = N^{0,m}(p)= \frac{(p_r-\ii p_i)^m }{m!},
    \end{equation}
    for $(p_r,p_i):=p\in \R^2$. For real-valued functions $u, \tilde{a}$ and $\tilde{f}$, from definitions and $\partial_{\C}^{(k, \ell)} = \overline{\partial_{\C}^{(\ell,k)}}$, one can directly check that
    {\small
    \begin{equation}\label{C:pair}
        N^{k,\ell}(p)=\overline{N^{\ell,k}(p)}, \
        \tilde{a}^{k,\ell}_{m,n} = \overline{\tilde{a}^{\ell,k}_{n,m}}, \
        \tilde{A}^{k,\ell}_{m,n} = \overline{\tilde{A}^{\ell,k}_{n,m}}, \
        \tilde{F}_{k,\ell} = \overline{\tilde{F}_{\ell,k}}, \
        A^k_{m,n}(p) = \overline{A^k_{n,m}(p)}, \
        F(p) = \overline{F(p)}.
    \end{equation}
    }
    From the definition of $F(p)$ in \eqref{Fp}, one concludes from \eqref{C:pair} that $F(p)$ is real-valued.
    Hence, \eqref{u:taylor} can be rewritten as the following Taylor expansion using real-valued coefficients:
    \begin{equation}
    \label{u:taylor2}
        u(\bpt + ph)
        = u(\bpt) + \sum_{m = 1}^{M + 1} \sum_{k = m}^{M + 1} 2 \re \left( A^{k}_{m,0}(p) \partial_{\C}^{(m, 0)} u(\bpt) \right) h^k + F(p) + \bo_u(h^{M + 2}).
    \end{equation}
    %

%%%%%%%%%%%%%%%%%%%%%%%%%%%%%%%%%%%%%%%%%%%%%%%%%%%%%%%%%%%%%%%%%%
%%%%%%%%%%%%%%%%%%%%%%%%%%%%%%%%%%%%%%%%%%%%%%%%%%%%%%%%%%%%%%%%%%
    \section{Construction of Compact $9$-point FDM Schemes at Interior Grid Points}
    \label{sec: interior}

    We shall develop our FDM schemes separately according to whether the stencil center is an interior or boundary grid point. In this section, we deal with sixth-order $9$-point compact interior stencils, while the boundary stencils will be handled in the next section.

    Let $\SS$ be the reference stencil $[-1,1]^2 \cap \Z^2$ centered at $(0,0)$. By definition of $\Omega_h^\circ$ in \eqref{eq: Omegah}, each grid point $\spt\in \Omega_h^\circ$ will serve as the stencil center and all its $1$-ring neighboring grid points $\spt+ ph$, $p \in \SS$ lie inside $\Omega$. Now we expand the solution $u$ in \eqref{u:taylor2} at each point $\spt+ ph$ for $p \in \SS$ at the base point $\bpt:=\spt$. In view of this, for each stencil point $\spt+ ph$ we aim to find the stencil coefficient $C_p(h) \in \R$, a real polynomial of variable $h$, such that for a given positive integer $M\in \mathbb{N}$,
    \begin{equation}
        \label{eq: interior scheme 0}
        \sum_{p \in \SS} C_p(h) u(\spt + ph) = \sum_{p \in \SS} C_p(h) F(p) + \bo(h^{M + 2}),
    \end{equation}
    where $F(p)$ is defined in \eqref{Fp} and is real-valued. Here and afterwards, any summation $\sum_{k = m}^n$ with $m > n$ is treated as $0$.
    The conditions on $C_p(h)$ in \eqref{eq: interior scheme 0} are given by the following lemma.

    \begin{lemma}
        \label{lem: interior constraint of c}
        Let $M\in \mathbb{N}$ and define
        $C_p(h) := \sum_{k = 0}^{M + 1} c_{p, k} h^k$ with $c_{p, k} = \bo_{\tilde{a}} (1)$ for $p \in \SS$. Then the linear system in \eqref{eq: interior scheme 0} with the remainder term $\bo_{\tilde{a}, u} (h^{M + 2})$ holds if and only if
           {\footnotesize{\begin{equation}\label{eq: interior constraint of c}
            \begin{aligned}
                & \sum_{p \in \SS} \re\left(A^m_{m, 0}(p)\right)  c_{p, j}
                = -\sum_{k = 0}^{j - 1} \sum_{p \in \SS} \re\left(A^{m + j - k}_{m, 0}(p) \right) c_{p, k}, \quad \forall \, j = 0,\ldots, M + 1, \ m = 0, \ldots, M + 1 - j, \\
                & \sum_{p \in \SS} \im\left(A^m_{m, 0}(p) \right) c_{p, j}
                = -\sum_{k = 0}^{j - 1} \sum_{p \in \SS} \im\left(A^{m + j - k}_{m, 0}(p) \right) c_{p, k}, \quad \forall \, j = 0,\ldots, M + 1, \ m = 1, \ldots, M + 1 - j,
            \end{aligned}
        \end{equation}
        }}%
        where the quantities $A^k_{m,n}$ are defined in \eqref{Akmnp}. Note that $A^k_{0,0}=\td(k)$ for all $k\in \NN$.
    \end{lemma}

    \begin{proof}
        By expanding $u(\bpt + ph)$ at the base point $\bpt$ via \eqref{u:taylor2}, we obtain
        \begin{align*}
        \sum_{p \in \SS} &C_p(h) u(\bpt + ph)
        = u(\bpt) \sum_{p\in \SS} C_p(h)+
        \sum_{p \in \SS} C_p(h) F(p) + \bo_{\tilde{a}, u} (h^{M + 2})\\
        &+\sum_{m=1}^{M+1} \sum_{k= m}^{M + 1} \sum_{p \in \SS} 2\left[\re(A^k_{m, 0}(p)) \re(\partial_{\C}^{(m, 0)} u)-
        \im(A^k_{m, 0}(p)) \im(\partial_{\C}^{(m, 0)} u)\right]C_p(h) h^k.
        \end{align*}
        Treating all $u$, $\re(\partial^{(m,0)}_\C u)$ and $\im(\partial^{(m,0)}_\C u)$ for $m=1,\ldots, M+1$ as independent variables, we deduce from the above identity that \eqref{eq: interior scheme 0} becomes
        \begin{equation}
            \label{Cp:interior}
            \sum_{k = m}^{M + 1} \sum_{p \in \SS} A^k_{m, 0}(p)C_p(h)h^k = \bo_{\tilde{a}, u} (h^{M+2}),\quad m=0,\ldots, M+1,
        \end{equation}
        where we used the fact $A^k_{0, 0} = \td(k)$. Now plugging $C_p(h)=\sum_{j=0}^{M+1} c_{p,j} h^j$ into \eqref{Cp:interior}, we have
        \[
        \bo_{\tilde{a}, u} (h^{M+2})=
        \sum_{k=m}^{M+1} \sum_{j=0}^{M+1} \sum_{p\in \SS} A^k_{m,0}(p) c_{p,j} h^{j+k}
        =\sum_{j=m}^{M+1} \sum_{k=0}^{j-m} \sum_{p\in \SS} A^{j-k}_{m,0}(p) c_{p,k} h^j+\bo_{\tilde{a}, u} (h^{M+2}).
        \]
        Because $h$ is independent, we conclude that the above identity is just $\sum_{k=0}^{j-m} \sum_{p\in \SS} A^{j-k}_{m,0} c_{p,k}=0$ for $0\le m\le M+1$ and $m\le j\le M+1$, which is equivalent to \eqref{eq: interior constraint of c} by replacing $j-m$ with the new index $j$.
    \end{proof}

    The constraint in \Cref{lem: interior constraint of c} is further investigated in the following proposition, which also provides a constructive way of generating stencil coefficients. Note that we can arrange the elements in the reference stencil $\SS:=[-1,1]^2\cap \Z^2$ with $\#\SS=9$ in the following order:
    \begin{equation}\label{S:order}
    (-1, -1), \quad (-1, 0),\quad (-1, 1),\quad (0, -1),\quad (0, 0),\quad (0, 1),\quad (1, -1),\quad (1, 0),\quad (1, 1).
    \end{equation}
    Throughout the paper, we shall always use this ordering of $\SS$ to translate the set $\{c_{p,j} \setsp p\in \SS\}$ into a column vector $\vec{c}_j\in \R^9$. Recall that we identify a point $(p_r, p_i):=p\in \R^2$ with the complex number $p_r+\ii p_i\in \C$ in our calculation.

    \begin{prop} \label{prop: interior c}
    Let $M\in \mathbb{N}$ and $\spt \in \Omega_h^\circ$ be a stencil center. Then
    the linear system \eqref{eq: interior constraint of c} with $j=0$ has a nonzero solution with $\vec{c}_0\ne 0$ if and only if $M \leq 6$. Moreover, for $M = 6$,
    there always exist real-valued coefficients $c_{p, j}\in \R$ for $p \in \SS$ and $j=0,\ldots, 7$ such that
    \begin{enumerate}
        \item[(i)] $\{ c_{p, j}: p \in \SS, j=0,\ldots, 7 \}$ is a real-valued solution to \eqref{eq: interior constraint of c} with $\vec{c}_0 \ne 0$ and $c_{p,j} =\bo_{\tilde{a}} (1)$, and \eqref{eq: interior scheme 0} holds with $C_p(h) := \sum_{j = 0}^7 c_{p, j} h^j\in \mathbb{R}$, $p \in \SS$ and the remainder term $\bo_{\tilde{a}, u} (h^8)$;
        \item[(ii)] These real numbers $\{ c_{p, j}\}_{p \in \SS}$ for $j=0,\ldots,7$ satisfy the following sign condition:
        \begin{equation}\label{sign:cond}
            c_{(0, 0), 0} > 0, \ c_{p, 0} < 0
            \quad \mbox{and}\quad
            c_{(0,0),j}\ge 0, \ c_{p, j} \le 0, \ j=1,\ldots,7
            \quad \mbox{for all } p \in \mathring{\SS} := \SS \setminus \{(0,0)\};
        \end{equation}
        In particular, $C_{(0,0)}(h)>0$ and $C_p(h) < 0$ for all $p \in \mathring{\SS}$.
        \item[(iii)] For all $j=0,\ldots,7$, these real numbers $\{ c_{p, j}\}_{p \in \SS}$ satisfy the sum condition $\sum_{p \in \SS} c_{p, j} = 0$.
    \end{enumerate}
    \end{prop}

    \begin{proof}
        Consider $\bpt:=\spt$ as the base point.
        For each $j=0,\ldots,M + 1$, \eqref{eq: interior constraint of c} consists of $2M + 3 - 2j$ linear equations with $9$ unknowns $\{c_{p,j}\}_{p\in \SS}$. Using the default ordering of the set $\SS$ given above in \eqref{S:order}, the linear equations \eqref{eq: interior constraint of c} can be equivalently expressed in the matrix form $\mathbb{A}_j \vec{c}_j = \vec{b}_j$ for $j=0,\ldots, M+1$, where $\mathbb{A}_j$ is an $(2M + 3 - 2j)\times 9$ matrix and $\vec{c}_j, \vec{b}_j\in \mathbb{R}^9$. By \eqref{Amm0}, for each  $(p_r, p_i):=p \in \SS$,
        the entries of the $(2M+3)\times 9$ matrix $\mathbb{A}_0$ are given by
        {\small{ \begin{equation}\label{A0}
            \mathbb{A}_0(1, p) = 1, \; \mathbb{A}_0(2m, p) = \re \frac{(p_r + \ii p_i)^m}{m!}, \; \mathbb{A}_0(2m + 1, p) = \im \frac{(p_r + \ii p_i)^m}{m!}, \quad m=1,\ldots, M + 1
        \end{equation}
        }}%
        and
        \begin{equation}\label{Aj}
            \mathbb{A}_j(k,p)=\mathbb{A}_0(k,p), \qquad k=1,\ldots, 2M+3-2j, \ j=1,\ldots, M+1.
        \end{equation}
        That is, the $(2M + 3 - 2j)\times 9$ matrix $\mathbb{A}_j$
        is just the submatrix of $\mathbb{A}_0$ by taking its first $2M + 3 - 2j$ rows.
        Moreover, the vector $\vec{b}_0$ is identically zero, and for each $j=1,\ldots, M+1$,
        {\footnotesize
        \begin{equation}
            \label{eq: proof 2 interior c}
            \vec{b}_j (1) =0, \quad
            \vec{b}_j (2m) = -\sum_{k = 0}^{j - 1} \sum_{p \in \SS} \re \left(A^{j+m-k}_{m, 0}(p) \right) c_{p, k}, \quad
            \vec{b}_j (2m + 1) = -\sum_{k = 0}^{j - 1} \sum_{p \in \SS} \im\left(A^{j+m-k}_{m, 0}(p) \right) c_{p, k},
        \end{equation}
        }%
        for $m=1,\ldots, M+1-j$. It is very important to notice that all the entries of $\vec{b}_j$ only depend on previous $\vec{c}_0,\ldots, \vec{c}_{j-1}$. Hence, it is not surprising that we solve the linear systems $\mathbb{A}_j \vec{c}_j=\vec{b}_j$ in the natural ordering $j=0,\ldots,M+1$.
        By symbolic calculation, the ranks of the $(2M+3)\times 9$ matrices $\mathbb{A}_{0}$ of constants for $M=0,\ldots,7$ are $3,5,7,8,8,8,8,9$. Because $\vec{b}_0=0$, as a consequence,
        the homogeneous linear system $\mathbb{A}_0 \vec{c}_0 = 0$ has a nontrivial solution $\vec{c}_0$ if and only if $M \leq 6$. Moreover, for $M=6$, up to a multiplicative constant, all the solutions to $\mathbb{A}_0 \vec{c}_0 = 0$  is given by
        \begin{equation}\label{vec:c0}
        \vec{c}_0 = [-1, -4, -1, -4, 20, -4, -1, -4, -1].
        \end{equation}

        Now we only consider $M=6$ for solving the linear systems \eqref{eq: interior constraint of c}. We solve $\mathbb{A}_j \vec{c}_j=\vec{b}_j$ in the order of $j=0,\ldots, 7$ via symbolic calculation and present in Appendix \ref{app: Interior stencil coefficients} one possible real-valued solution with $c_{p, j} = \bo_{\tilde{a}} (1)$ and $\vec{c}_0$ given in \eqref{vec:c0}.
        By \Cref{lem: interior constraint of c}, we conclude that
        \eqref{eq: interior scheme 0} must hold with $C_p(h) := \sum_{j = 0}^7 c_{p, j} h^j$, $p \in \SS$ and the remainder term $\bo_{\tilde{a}, u} (h^8)$. Hence, item (i) holds.

        Because $\mathbb{A}_j(1,p)=\mathbb{A}_0(1,p)=1$ for all $p\in \SS$ and $\vec{b}_j(1)=0$, every solution to \eqref{eq: interior constraint of c} implies
        \[
        \sum_{p\in \SS} c_{p,j}=\sum_{p\in \SS} \mathbb{A}_j(1,p) c_{p,j}=[\mathbb{A}_j \vec{c}_j]_1=\vec{b}_j(1)=0.
        \]
        This proves that item (i) always guarantees the sum condition in item (iii). Unfortunately, the sign condition in \eqref{sign:cond} is only satisfied for $j=0$ by \eqref{vec:c0}. We now modify it so that all items (i)-(iii) are satisfied.
        For any real numbers $q_0,\ldots,q_7\in \mathbb{R}$,
        we define
        \begin{equation}\label{tcpj}
        \tilde{C}_p(h) := \sum_{j=0}^7 \tilde{c}_{p,j} h^j\quad \mbox{ and }\quad
        \tilde{c}_{p, j}:= \sum_{k = 0}^j q_{j - k} c_{p, k},\quad p\in \SS, \ j=0,\ldots, 7.
        \end{equation}
        Then we trivially have $\tilde{C}_p(h)= C_p(h) Q(h) + \bo_{\tilde{a}} (h^8)$ for all $p\in \SS$ with $Q(h):= \sum_{j = 0}^7 q_j h^j$. Because $Q$ is independent of $p\in \SS$, by \Cref{lem: interior constraint of c}, items (i) and (iii) must be satisfied with the original solution $c$ being replaced by the modified $\tilde{c}$. We now choose $q_0,\ldots, q_7$ so that item (ii) is also satisfied.
        By \eqref{vec:c0}, we see that $c_{p,0}\ne 0$ and we can define $q_j$, $j=0,\ldots,7$ by
        \begin{equation}\label{q:modcoeff}
            q_0 := 1 \quad \mbox{ and }\quad
            q_j := \sum_{k = 1}^j \lambda_k q_{j - k}
            \quad \mbox{with} \quad
            \lambda_j := \max \left\{ \max_{p \in \SS} \left( -\frac{c_{p, j}}{c_{p, 0}} \right), 0 \right\}
            ,\quad j=1,\ldots,7.
        \end{equation}
        Then we can prove by induction and \cref{tcpj} that item (ii) holds for $\tilde{c}_{p, j}$.
    \end{proof}

    We finish this section by discussing the special case $-\Delta u=f$.
    Then $\Delta u = \tilde{f}$ in \eqref{pde:mod} with $\tilde{f}:=-f$ and $\tilde{a}:=0$. Due to $\tilde{a}=0$, for $p\in \R^2$,  we can easily obtain
    \[
    \tilde{A}^{k, \ell}_{m, n} = \td(k-m) \td(\ell-n),
    \quad
    \tilde{F}_{k, \ell} = \td(k) \td(\ell) \partial_{\C}^{(k - 1, \ell - 1)} \tilde{f},
    \quad
    A^k_{m,n}(p) = \td (k-m-n) N^{m,n}(p)
    \]
    for $k, \ell \in \NN$, $(m, n) \in \LN^{k}_{\ell}$, and $F(p) = \sum_{0\le k + \ell \leq M - 1} N^{k+1,\ell+1}(p) h^{k + \ell + 2} \partial_{\C}^{(k, \ell)} \tilde{f}$.
    Hence, for each stencil point $p\in \SS$, the linear equations in \eqref{eq: interior constraint of c} become
    \[
    \sum_{p \in \SS} A^{m+n}_{m,n}(p) c_{p, j} = 0,
    \quad \forall \, j=0,\ldots, M + 1, \ (m, n) \in \LN^{M + 1 - j}_{M + 1 - j}.
    \]
    For $M = 6$, up to a nonzero multiplicative constant to all real numbers $c_{p,j}$, all the real-valued solutions $\{c_{p,j}\; : \; p\in \SS, j=0,\ldots, 7\}$ to the above linear system are given by
    \begin{align*}
        & \vec{c}_0=\kappa_0 \vec{v}_1, \
        \vec{c}_1=\kappa_1 \vec{v}_1, \
        \vec{c}_2=\kappa_2 \vec{v}_1, \
        \vec{c}_3=\kappa_3 \vec{v}_1
        \quad \mbox{with}\quad \vec{v}_1 := [-1, -4, -1, -4, 20, -4, -1, -4, -1], \\
        & \vec{c}_4 = \kappa_{4} \vec{v}_2 + \kappa_{5} (\vec{v}_1 - \vec{v}_2)
        \quad \mbox{with}\quad \vec{v}_2 := [-1, 0, -1, 0, 4, 0, -1, 0, -1], \\
        &\vec{c}_5 = \kappa_{6} \vec{v}_2 + \kappa_{7} [0, -2, 1, 0, 4, -2, -1, 0, 0] \\
        & \qquad + \kappa_{8} [0, 0, -1, -2, 4, 0, 1, -2, 0] + \kappa_{9} [1, -1, 0, -1, 0, 1, 0, 1, -1], \\
        & \vec{c}_6 = [-\kappa_{10},\, -2\kappa_{11},\, -\kappa_{12},\, -2\kappa_{13},\, 2(\kappa_{10} + \kappa_{12} + \kappa_{13} + \kappa_{14}) + 4\kappa_{11},\, -2\kappa_{14},\, \\
        & \hspace{32pt} \kappa_{13} + \kappa_{15} - \kappa_{11} - \kappa_{12} - \kappa_{14},\, -2\kappa_{15},\, \kappa_{14} + \kappa_{15} - \kappa_{10} - \kappa_{11} - \kappa_{13}], \\
        & \vec{c}_7 = [-\kappa_{16}, -\kappa_{17}, -\kappa_{18}, -\kappa_{19}, \kappa_{16}+\cdots+\kappa_{23}, -\kappa_{20}, -\kappa_{21}, -\kappa_{22}, -\kappa_{23}],
    \end{align*}
    where $\kappa_0,\ldots,\kappa_{23} \in \R$ are free parameters. Moreover, all the items (i)--(iii) of \Cref{prop: interior c} are satisfied if $\kappa_0 > 0$, $\min \{ \kappa_{6}, \frac{1}{2} \kappa_{7}, \frac{1}{2} \kappa_{8} \} \geq |\kappa_{9}|$ and all the remaining free parameters $\kappa_j\ge 0$.

    By the definition of $F(p)$ in \eqref{Fp}, the right-hand side of \eqref{eq: interior scheme 0} without $\bo(h^{M+2})$ becomes
    \[
    \sum_{p\in \SS} C_p(h)F(p)=h^2 \sum_{0\le k + l \leq 5} \sum_{p \in \SS} C_p(h) N^{k+1,\ell+1}(p) h^{k + \ell} \partial_{\C}^{(k, \ell)} \tilde{f}.
    \]
    Using $\tilde{f}=-f$ and the definition \eqref{def:N}, we obtain from \eqref{eq: interior scheme 0} the general sixth-order finite difference scheme for the Poisson equation $-\Delta u=f$, where $\partial^{(k,\ell)}f$ are evaluated at the base point $\bpt \in \Omega^\circ_h$:
    \begin{align*}
        & h^{-2} \sum_{p \in \SS} C_p(h) u_h(\bpt+ph) \\
        &= (\kappa_0 + \kappa_1 h) \left( 6f + \frac{1}{2} h^2 (\partial^{(2, 0)} f + \partial^{(0, 2)} f) + \frac{1}{60} h^4 (\partial^{(4, 0)} f + 4\partial^{(2, 2)} f + \partial^{(0, 4)} f) \right) \\
        & \quad \, + h^2 (\kappa_2 + \kappa_3 h) \left( 6f + \frac{1}{2} h^2 (\partial^{(2, 0)} f + \partial^{(0, 2)} f) \right) + h^4 (2 \kappa_{4} + \kappa_{5}) f + h^5 (2 \kappa_{6} + \kappa_{7} + \kappa_{8}) f.
    \end{align*}
    Hence, we constructed all possible sixth-order compact FDMs satisfying items (i)--(iii) of \Cref{prop: interior c} with $M=6$ in the sense that we ignored the terms of $\bo(h^6)$ on the left-hand side (which do not affect the order $6$ of the scheme). Setting all free parameters to $0$ except for $\kappa_0 = 1$ in the above stencil coefficients, we obtain the known sixth-order finite difference scheme (e.g., see \cite{settle2013derivation,wang2009sixth,zhai2013family} in the literature) for the Poisson equation $-\Delta u=f$.

%%%%%%%%%%%%%%%%%%%%%%%%%%%%%%%%%%%%%%%%%%%%%%%%
%%%%%%%%%%%%%%%%%%%%%%%%%%%%%%%%%%%%%%%%%%%%%%%%
    \section{Construction of the FDM Schemes at Boundary Grid Points}
    \label{sec: boundary}

    We now develop our finite difference schemes for a boundary grid point $\spt \in \partial \Omega_h$ using its associated nearby base point $\bpt \in \partial \Omega$. Because the boundary curve $\partial \Omega$ is smooth, we can obtain a parametric equation in a neighborhood of the base point $\bpt$ on $\partial\Omega$:
    \begin{equation}\label{bdcurve}
    x=\gx(t),\; y=\gy(t),\quad t\in (t^*-\gep, t^*+\gep) \quad \mbox{with}\quad \bpt=(\gx(t^*), \gy(t^*)), \; (\gx'(t^*), \gy'(t^*)) \ne (0, 0)
    \end{equation}
    for some $\gep>0$. For example, if $\partial \Omega$ is given by a level set $\Phi(x,y)=0$. Then we may obtain $y=\varphi(x)$ in a neighborhood of $\bpt \in \partial \Omega$ such that $\Phi(x,\varphi(x))=0$. Hence, we may employ the parametric equation $\gx(t)=t^*+t$, $\gy(t)=\varphi(\gx(t))$ for $t\in (t^*-\gep,t^*+\gep)$, where $t^*$ is the $x$-coordinate of the base point $\bpt\in \partial \Omega$.

    Let $\theta$ be the tangent angle at $\bpt\in \partial \Omega$. More precisely,
    \begin{equation}\label{theta}
        \theta := \operatorname{Arg}(z_0) \in (-\pi,\pi]
        \quad \mbox{with}\quad
        z_0 := \gx'(t^*)+\ii\gy'(t^*)\ne 0.
    \end{equation}
    Then one can observe that
    \[
    2 e^{\ii\theta} \partial^{(1,0)}_\C=(\cos \theta +\ii \sin \theta)\left(\frac{\partial}{\partial x}-\ii \frac{\partial}{\partial y}\right)=
    \left(\cos \theta \frac{\partial}{\partial x}+\sin \theta \frac{\partial}{\partial y}\right)+\ii\left( \sin \theta \frac{\partial}{\partial x}-\cos \theta \frac{\partial}{\partial y}\right),
    \]
    where the real and imaginary parts are the directional derivatives along the tangent direction and the normal direction, respectively. Hence, it is very natural to consider $2^n e^{\ii n\theta}\partial^{(n,0)}_\C=[2e^{\ii\theta} \partial^{(1,0)}_\C]^n$.

    \subsection{Constraints on stencil coefficients of boundary stencils}
    \label{sec: Constraints on the stencil coefficients}

    In this section, we aim to derive an analog of equations \eqref{eq: interior constraint of c} for the stencil coefficients at the boundary grid point. We start from the representation \eqref{u:taylor2}, where the functions are expanded at a base point $\bpt\in \partial \Omega$. In this representation, there are altogether $2M + 4$ ``unknowns'':     $\partial^{(k,\ell)}_\C u, (k,\ell)\in \LN^{M+1}_{M+1}$ and $h$. \Cref{lem: B tilde} shows how we can differentiate the boundary condition $u(\gx(t), \gy(t)) = g(\gx(t),\gy(t))$ to get the constraints on the complex partial derivatives $\partial^{(k,\ell)}_\C u$. These constraints help us eliminate roughly half of the unknowns in \eqref{u:taylor2}. As a remark, this elimination process cannot be successfully carried out if we adopt standard partial derivatives instead.

    \begin{lemma}\label{lem: B tilde}
        Using the parametric equation \eqref{bdcurve} of the boundary $\partial \Omega$, we define a one-dimensional function $\tilde{g}(t):=g(\gx(t),\gy(t))$ for $t\in (t^*-\gep,t^*+\gep)$. Then
        for every $m\in \N$,
        \begin{equation}\label{deriv:g:n}
            \frac{\tilde{g}^{(m)}(t^*)}{m!}
            =\Big[\frac{1}{m!} \dd{t}[m] u(\gx(t),\gy(t))\Big]\Big|_{t=t^*}
            = \sum_{n=1}^m 2 \re \left( \tilde{B}_{m, n}(t^*) e^{\ii n\theta} [\partial_{\C}^{(n, 0)} u](\bpt) \right) + \tilde{G}_m (t^*),
        \end{equation}
        where $\bpt:=(\gx(t^*),\gy(t^*))$,
        the quantities $\tilde{B}_{m, n}(t^*)$ and $\tilde{G}_m(t^*)$ are defined by
        \begin{equation}\label{tBnm}
        \tilde{B}_{m, n}(t^*):=
        \sum_{j=n}^m \sum_{\ell=0}^{m-j} N^{j,m-\ell-j}_\ell \tilde{A}^{j,m-\ell-j}_{n,0}e^{-\ii n\theta}
        \quad \mbox{and}\quad
        \tilde{G}_m(t^*):=\sum_{j=0}^m \sum_{\ell=0}^{m-j} N^{j,m-\ell-j}_\ell \tilde{F}_{j,m-\ell-j},
        \end{equation}
        where all the complex numbers $N^{k,\ell}_j$ for $j,k,\ell\in \NN$ are recursively defined by
    \begin{equation}\label{Nklj}
        N^{0,0}_j:=\td(j), \quad
        N^{k,\ell}_j:=\frac{1}{k}\sum_{n=0}^j N^{k-1,\ell}_n z_{j-n}(t^*), \quad
        N^{\ell,k}_j=\overline{N^{k,\ell}_j},
        \quad j, \ell\in \NN, k\in \N.
    \end{equation}
    with $z_j(t^*) := \frac{\gx^{(j+1)}(t^*)+\ii \gy^{(j+1)}(t^*)}{(j+1)!}$ for all $j\in \NN$. Note that all $\tilde{G}_m$ in \eqref{tBnm} are real-valued.
    \end{lemma}

    \begin{proof}
        Define $p(t):=(\frac{\gx(t)-\gx(t^*)}{t-t^*},\frac{\gy(t)-\gy(t^*)}{t-t^*})$. Note that $\tilde{g}(t)=u(\gx(t),\gy(t))=u(\bpt + p(t) \tilde{h})$ with $\tilde{h} := t-t^*$. Using the Taylor expansion in \eqref{u:taylor2} with $M=\infty$, we have $\tilde{g}(t) = u(\gx(t),\gy(t))$ and
        \begin{equation}
            \label{eq: tg0}
            u(\gx(t),\gy(t))=
            u(\bpt) + \sum_{n = 1}^{\infty} 2 \re \left(\left( \sum_{k = n}^{\infty} A^{k}_{n,0}(p(t)) \tilde{h}^k \right) [\partial_{\C}^{(n, 0)} u](\bpt) \right) + F(p(t)),
        \end{equation}
        where $A^k_{n,m}$ is defined in \eqref{Akmnp} and $F(p(t))$ is defined in \eqref{Fp}. Note that
        \[
            p(t)=\left(\sum_{j=0}^\infty \frac{\gx^{(j+1)}(t^*)}{(j+1)!} \tilde{h}^j,
            \sum_{j=0}^\infty \frac{\gy^{(j+1)}(t^*)}{(j+1)!} \tilde{h}^j\right)=
            \left(\sum_{j=0}^\infty \frac{1}{2} (z_j(t^*)+\overline{z_j(t^*)}) \tilde{h}^j,
            \sum_{j=0}^\infty \frac{1}{2\ii}(z_j(t^*)-\overline{z_j(t^*)}) \tilde{h}^j\right)
        \]
        and by \eqref{Akmnp}, $A^{k}_{n,0}(p(t))=\sum_{j=n}^k N^{j,k-j}(p(t)) \tilde{A}^{j,k-j}_{n,0}$.
        Now from \eqref{def:N} and \eqref{Nklj}, we have
        \[
            N^{k,\ell}(p(t))=
            \frac{(\sum_{j=0}^\infty z_j(t^*) \tilde{h}^j)^k (\sum_{j=0}^\infty \overline{z_j(t^*)} \tilde{h}^j)^\ell}{k! \ell!}
            =\sum_{j=0}^{\infty} N^{k,\ell}_j \tilde{h}^j,
        \]
        where we used the fact that $N^{k,\ell}(p)=\overline{N^{\ell,k}(p)}$ in \eqref{C:pair} and hence $N^{k,\ell}_j=\overline{N^{\ell,k}_j}$. Therefore,
        \[
            \sum_{k = n}^{\infty} A^{k}_{n,0}(p(t)) \tilde{h}^k
            =\sum_{k = n}^{\infty} \sum_{j=n}^k \sum_{\ell=0}^\infty N^{j,k-j}_\ell \tilde{A}^{j,k-j}_{n,0} \tilde{h}^{k+\ell}
            =\sum_{m=n}^\infty \left(\sum_{j=n}^m \sum_{\ell=0}^{m-j} N^{j,m-\ell-j}_\ell \tilde{A}^{j,m-\ell-j}_{n,0} \right) \tilde{h}^m,
        \]
        which is just $\sum_{m=n}^\infty \tilde{B}_{m,n}(t^*) e^{\ii n\theta} \tilde{h}^m$ by \eqref{tBnm}.
        On the other hand, we deduce from \eqref{Fp} that
        \[
            F(p(t))=\sum_{j,k=0}^\infty N^{j,k}(p(t))\tilde{F}_{j,k} \tilde{h}^{j+k}
            =\sum_{j,k=0}^\infty \sum_{\ell=0}^\infty N^{j,k}_\ell \tilde{F}_{j,k} \tilde{h}^{j+k+\ell}=
            \sum_{m=0}^\infty \left(\sum_{j=0}^m \sum_{\ell=0}^{m-j} N^{j,m-\ell-j}_\ell \tilde{F}_{j,m-\ell-j}\right) \tilde{h}^m,
        \]
        which is just $\sum_{m=0}^\infty \tilde{G}_m(t^*) \tilde{h}^m$ by \eqref{tBnm}.
        That is, by $\tilde{h}=t-t^*$, we proved
        \[
            \tilde{g}(t) = u(\bpt)+\sum_{n=1}^\infty
            \sum_{m=n}^\infty 2 \re \left( \tilde{B}_{m, n}(t^*) e^{\ii n\theta} [\partial_{\C}^{(n, 0)} u](\bpt) \right) (t-t^*)^m + \sum_{m=0}^\infty \tilde{G}_m (t^*) (t-t^*)^m,
        \]
        from which we have \eqref{deriv:g:n}.
        All $\tilde{G}_m$ in \eqref{tBnm} are real-valued due to \eqref{C:pair} and $N^{k,\ell}_j=\overline{N^{\ell,k}_j}$.
    \end{proof}

    By the definition of $\tilde{B}_{n,n}$ in \eqref{tBnm}, noting $\tilde{A}^{m,0}_{m,0}=1$ by \eqref{tA0} and $N^{m,0}_0=\frac{z_0^m}{m!}$ by \eqref{Nklj}, we have
    \[
    \tilde{B}_{m,m}(t^*) = e^{-\ii m\theta} N^{m,0}_0 \tilde{A}^{m,0}_{m,0}= \frac{|z_0(t^*)|^m}{m!}
    \quad \mbox{with}\quad z_0(t^*) := \gx'(t^*)+\ii \gy'(t^*)\ne 0.
    \]
    Consequently, dropping $t^*$ and $\bpt$ for simplicity, we can rewrite \eqref{deriv:g:n} as
    \begin{equation*}
        \re\left(e^{\ii m\theta}\partial_{\C}^{(m, 0)} u\right) = \frac{m!}{2|z_0|^m}\left[\frac{\tilde{g}^{(m)}}{m!}
        -\tilde{G}_m-
        \sum_{n = 1}^{m-1} 2\re\left(\tilde{B}_{m, n}e^{\ii n\theta}\partial_{\C}^{(n, 0)} u\right)\right].
    \end{equation*}
    Now we can recursively deduce that
    \begin{equation}\label{Re:Im:u}
    \re\left(e^{\ii m\theta}\partial_{\C}^{(m, 0)} u\right)=\sum_{n=1}^{m-1} B_{m,n}\im\left(e^{\ii n\theta}\partial_{\C}^{(n, 0)} u\right)+G_m,\qquad m=1,\ldots, M+1,
    \end{equation}
    where the real-valued quantities $B_{m,n}$ and $G_m$ are defined to be
    \begin{equation}
        \label{Bmn}
        B_{m,n} := \frac{m!}{|z_0|^m} \left[\im(\tilde{B}_{m,n})-\sum_{\ell=n+1}^{m-1} \re(\tilde{B}_{m,\ell}) B_{\ell,n}\right],
        \quad m \geq 2, \ n=1,\ldots,m-1
    \end{equation}
    and
    \begin{equation}\label{Gm}
    G_m:=\frac{m!}{2|z_0|^m}\left[ \frac{\tilde{g}^{(m)}}{m!}-\tilde{G}_m -\sum_{n=1}^{m-1} 2\re(\tilde{B}_{m,n})G_n\right], \quad m\in \N.
    \end{equation}
    Using \eqref{Re:Im:u} and the boundary condition $u(\bpt) = g(\bpt)$, for real-valued data $u,\tilde{a}, \tilde{f}$ and $g$, we obtain from the Taylor expansion in \eqref{u:taylor2} that
    \begin{equation}\label{u:taylor:bdr}
        u(\bpt+ p h)
        = \sum_{m = 1}^{M + 1}\sum_{k=m}^{M+1} A^k_m(p) h^k \im(e^{\ii m\theta}\partial^{(m,0)}_\C u)
         + G(p) + \bo_u (h^{M + 2}),
    \end{equation}
    where $\bpt+p h\in \Omega$ for $p\in \mathbb{R}^2$, and the real-valued quantities $A^k_m(p)$ and $G(p)$ are defined by
    \begin{equation}\label{Akm:Gm}
    \begin{split}
        &A^k_m (p) :=
        -2\im(A^k_{m,0} (p) e^{-\ii m\theta})+
        \sum_{n=m+1}^{k} 2 \re(A^k_{n,0}(p) e^{-\ii n\theta}) B_{n,m},\\
        &G(p) := g(\bpt) + F(p)+\sum_{m = 1}^{M+1} \sum_{k=m}^{M+1} 2\re(A^k_{m,0} (p) e^{-\ii m\theta}) G_m h^k.
        \end{split}
    \end{equation}

    Consider a boundary stencil center $\spt\in \partial \Omega_h$ and a reference stencil
    $\SS_{\spt} \subseteq \mathbb{Z}^2$ with $(0,0)$ referring to the stencil center $\spt$ such that $\SS_{\spt}$ has at most $8$ points of $\Z^2$.
    We shall consider a base point $\bpt\in \partial \Omega$ near the stencil center $\spt$ and then we define a shifting vector $\sv$ and its shift operator by
    \be
        \label{vec:shift}
        \sv:= (\spt-\bpt)/h
        \quad \mbox{with}\quad \|\sv\|\le \sqrt{2} \quad \mbox{and}\quad
        \psv:=p+\sv,\qquad p\in \mathbb{R}^2.
    \ee
    Note that $\spt+ph=\bpt+\psv h$ .
    In view of the identity \eqref{u:taylor:bdr}, we aim to find
    stencil coefficients $C_p(h):=\sum_{j=0}^M c_{p,j} h^j\in \mathbb{R}$ with $c_{p,j}\in \mathbb{R}$ for $p \in \SS_{\spt}$ such that for a given positive integer $M\in \mathbb{N}$,
    \begin{equation}
        \label{eq: boundary scheme 0}
        \sum_{p \in \SS_{\spt}} C_p(h) u(\spt + ph)
        = \sum_{p \in \SS_{\spt}} C_p(h) G(\psv) + \bo(h^{M + 2}),
    \end{equation}
    Note that $C_p(h)$ has one degree order lower than the interior stencil coefficients due to the use of Dirichlet boundary condition.
    The conditions on $C_p(h)$ in \eqref{eq: boundary scheme 0} are given by the following lemma.

    \begin{lemma}
        \label{lem: boundary constraint of c}
        Let $M\in \mathbb{N}$, $\spt\in \partial \Omega_h$, $\bpt\in \partial \Omega$ and $\sv$ as in \eqref{vec:shift}. Define
        $C_p(h) := \sum_{j = 0}^M c_{p, j} h^j$ with real-valued numbers $c_{p, j} = \bo_{\tilde{a}, \gx, \gy} (1)$ for $p \in \SS_{\spt}$. Then \cref{eq: boundary scheme 0} with the remainder term $\bo_{\tilde{a}, \gx, \gy, u} (h^{M + 2})$ holds if and only if
        \begin{equation}
            \label{eq: boundary constraint of c}
            -\sum_{p\in \SS_{\spt}}A^m_m(\psv) c_{p,j}=\sum_{k=0}^{j-1} \sum_{p\in \SS_{\spt}} A^{m+j-k}_m(\psv) c_{p,k},\qquad j=0,\ldots, M, \ m=1,\ldots, M+1-j.
        \end{equation}
        Note that $A^m_m(\psv) = -\frac{2}{m!} \im((\psv_r+\ii \psv_i)e^{-i\theta})^m$, where $(\psv_r, \psv_i)= \psv:=p+\sv=p+(\spt-\bpt)/h$.
    \end{lemma}

    \begin{proof}
        From \cref{u:taylor:bdr} and the fact $\spt + ph = \bpt + \psv h$, we have
        \[
        \sum_{p\in \SS_{\spt}} C_p(h) u(\spt+p h)=\sum_{m=1}^{M+1}\sum_{k=m}^{M+1} \sum_{p\in \SS_{\spt}} C_p(h) A^k_m(\psv) h^k \im(e^{\ii m\theta} \partial^{(m,0)}_{\C} u)+\sum_{p\in \SS_{\spt}} C_p(h) G(\psv).
        \]
        Treating $ \im(e^{\ii m\theta} \partial^{(m,0)}_{\C} u)$ for $m=1,\ldots, M+1$ as independent variables and using the definition of $C_p(h)$, we observe that
        \eqref{eq: boundary scheme 0} becomes
        \[
        \sum_{k=m}^{M+1} \sum_{j=0}^M \sum_{p\in \SS_{\spt}}  A^k_m(\psv) c_{p,j} h^{j+k}=\bo_{\tilde{a}, \gx, \gy, u}(h^{M+2}),\qquad m=1,\ldots, M+1.
        \]
        Changing the summation index as in \Cref{lem: interior constraint of c}, we obtain \eqref{eq: boundary constraint of c}. By \eqref{Amm0} and \eqref{Akm:Gm}, we have
        \[
            A^m_m(\psv)=-2\im(A^m_{m,0} e^{-\ii m\theta})
            =-\frac{2}{m!} \im ((\psv_r+\ii \psv_i)e^{-i\theta})^m.
        \]
        This completes the proof.
    \end{proof}

    \subsection{Construction of boundary stencils and their coefficients}
    \label{sec: choice of boundary stencil}

    From now on, we fix $M = 4$ and
    take a boundary grid point $\spt\in \partial \Omega_h$ as the stencil center. By the definition in \eqref{eq: Omegah}, there must exist $\bpt\in (\spt+h [-1,1]^2)\cap \partial \Omega\ne \emptyset$ and $\|\bpt-\spt\|\le \sqrt{2}h$. In practical implementation, we further require that the vector from $\spt$ to $\bpt$ be horizontal, vertical, or $\pm 45^\circ$. If not unique, then we take the one with the smallest $\|\bpt-\spt\|$. Note that $\bpt=(\gx(t^*),\gy(t^*))$ in \eqref{bdcurve}.

    Then the directed tangent line $L_{\bpt}$ at $\bpt$ to the boundary curve $\partial \Omega$ is given by
    \begin{equation}\label{L:line}
    L_{\bpt}:=\{(x,y)\in \mathbb{R}^2 \; : \; (x,y)=\bpt+t (\gx'(t^*),\gy'(t^*)),\; t\in \R\},
    \end{equation}
    and we define $\hp{L_{\bpt}}$ to be the open half plane on the left-hand side of the directed line $L_{\bpt}$. Without loss of generality, we can assume that $\spt\in \hp{L_{\bpt}}$; otherwise, we just change the variable $t$ into $-t$. We deduce from the Taylor expansion of the parametric equation \eqref{bdcurve} of $\partial \Omega$ at $\bpt$ that the distance between $L_{\bpt}$ and $\partial \Omega$ is bounded by $C h^2$ for all $t\in (t^*-2h, t^*+2h)$ with $C$ depending on the curvature of $\partial \Omega$ at $\bpt$.

\begin{figure}[htb!]
        \subfloat[Case I]{\includegraphics[width = 0.2 \linewidth, page = 1]{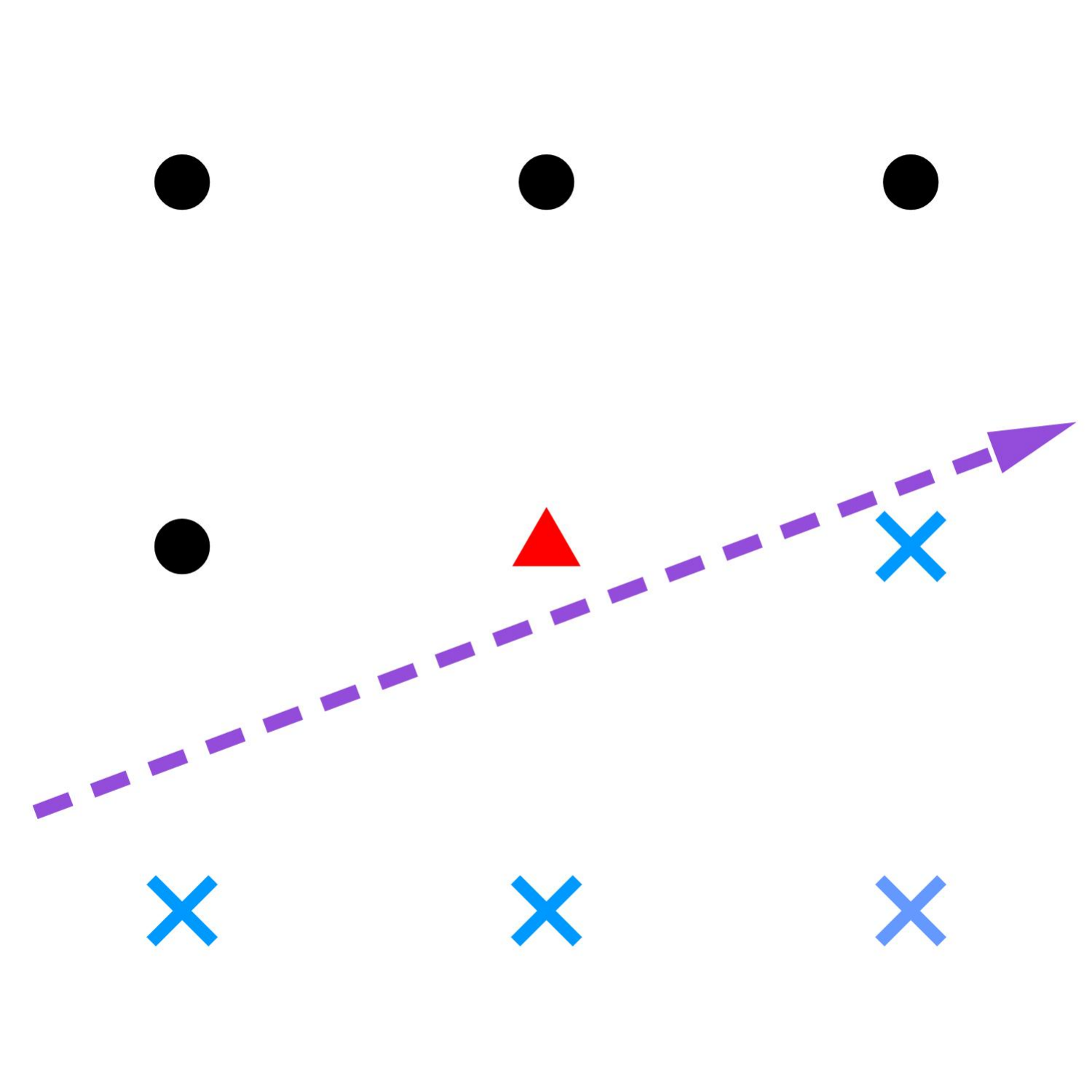}}
        \subfloat[Case II]{\includegraphics[width = 0.2 \linewidth, page = 2]{./figures/stencil.pdf}}
        \subfloat[Case III]{\includegraphics[width = 0.2 \linewidth, page = 3]{./figures/stencil.pdf}}
        \subfloat[Case IV]{\includegraphics[width = 0.2 \linewidth, page = 4]{./figures/stencil.pdf}}
        \subfloat[Case V]{\includegraphics[width = 0.2 \linewidth, page = 5]{./figures/stencil.pdf}}

        \caption{Five cases under condition $(\spt+h \SS)\cap \hp{L_{\bpt}}\subset \Omega$, where $\spt$ is the red triangle. The tangent line $L_{\bpt}$ is the purple dashed line with the arrow indicating the direction. $\hp{L_{\bpt}}$ is the left-hand region of $L_{\bpt}$. For simplicity, the base point $\bpt$ and $\partial \Omega$ are not explicitly shown above. All black dots belong to $\hp{L_{\bpt}}\cap \Omega$, while
        some blue crosses in $[h\mathbb{Z}^2]\setminus \hp{L_{\bpt}}$ may belong to $\overline{\Omega}$.}
        \label{fig: inside grid points cases}
    \end{figure}

\begin{figure}[htb!]
        \subfloat[{\tiny{Type 1 for Case I}}]{\includegraphics[width = 0.225 \linewidth, page = 6]{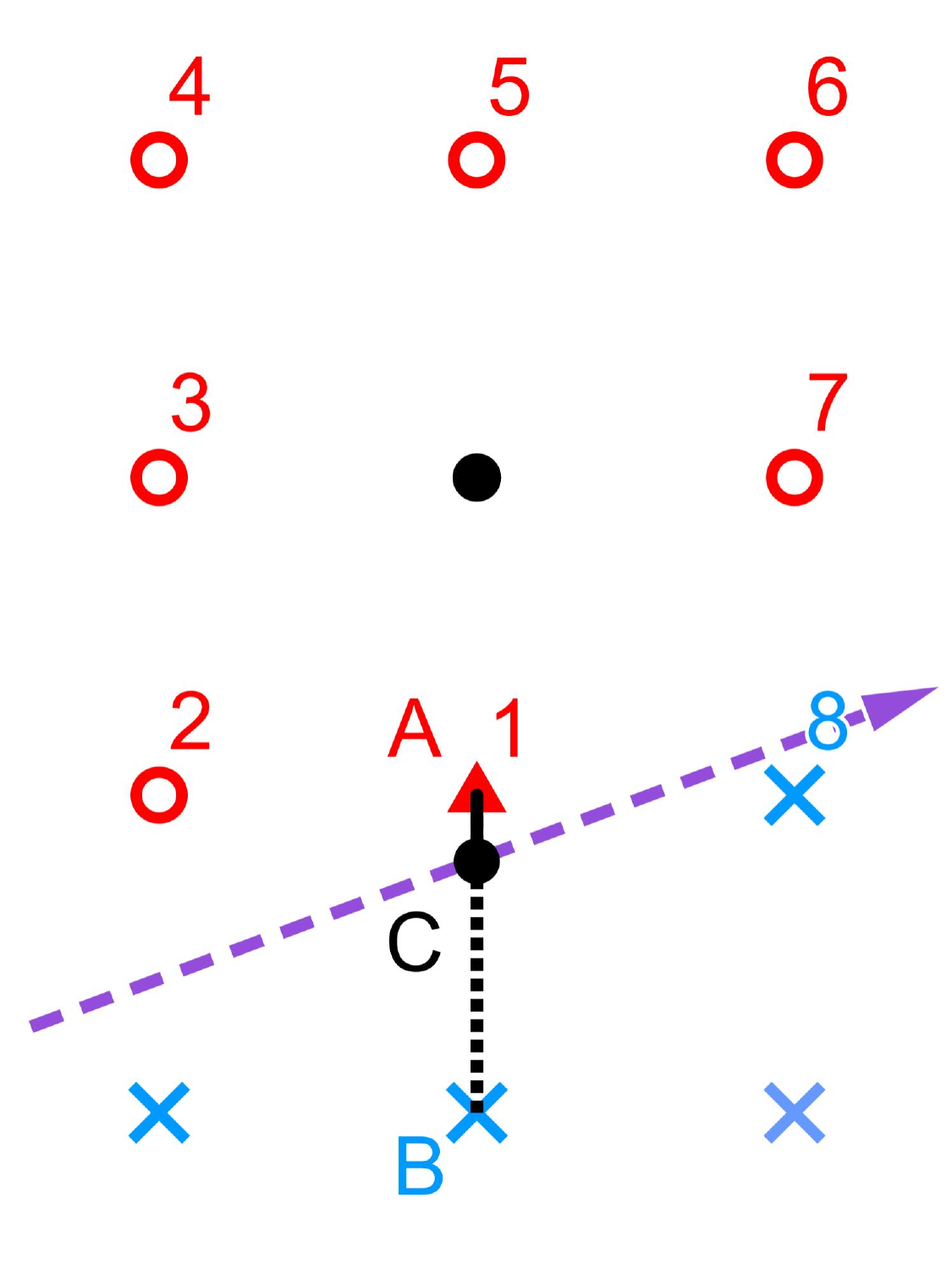}}
        \subfloat[{\tiny{Type 2 for Cases II \& IV}}]{\includegraphics[width = 0.225 \linewidth, page = 7]{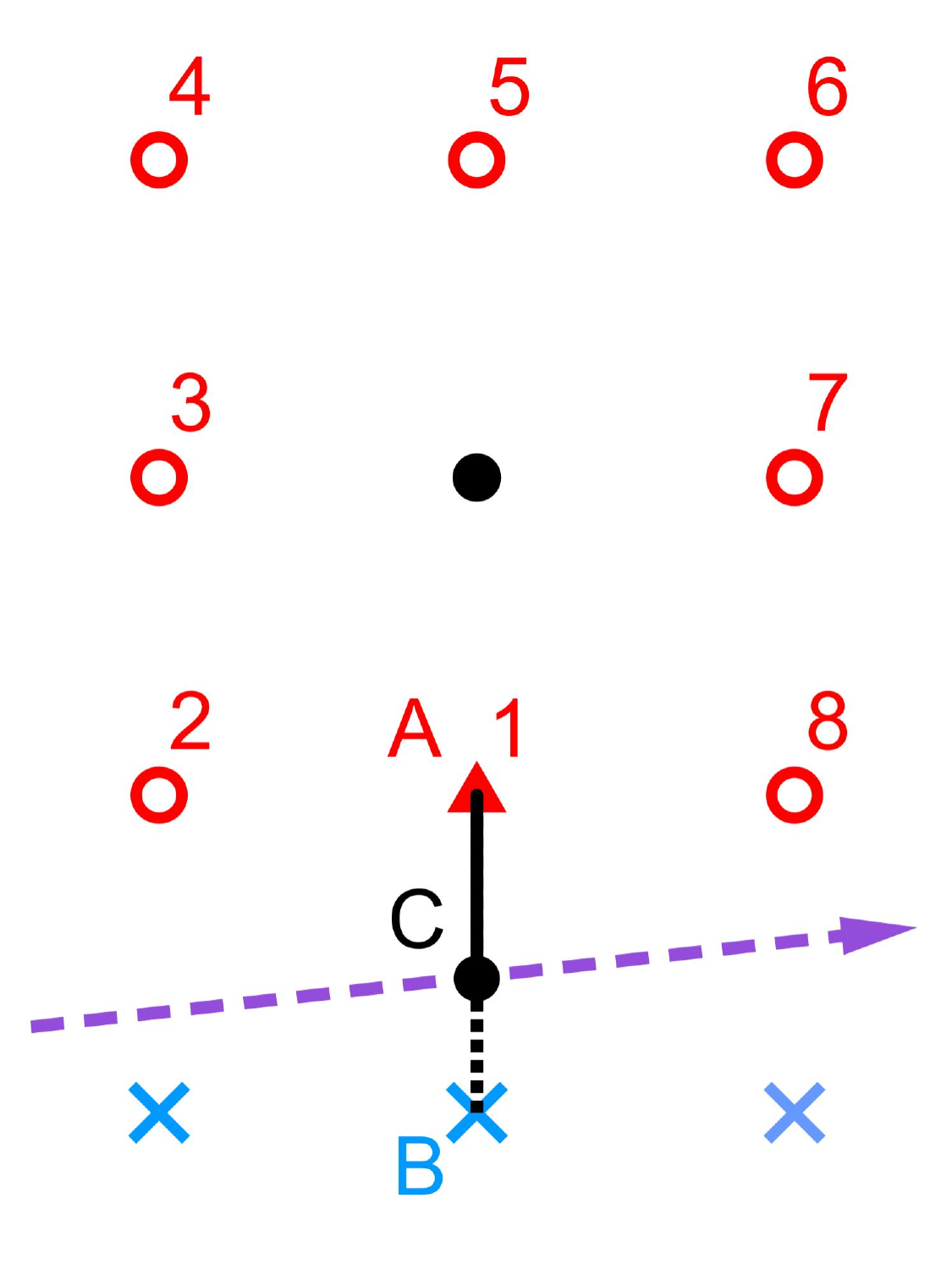}}
        \hspace{0.025\linewidth}
        \subfloat[{\tiny{Type 3 for Case III}}]
        {\includegraphics[width = 0.285 \linewidth, page = 8]{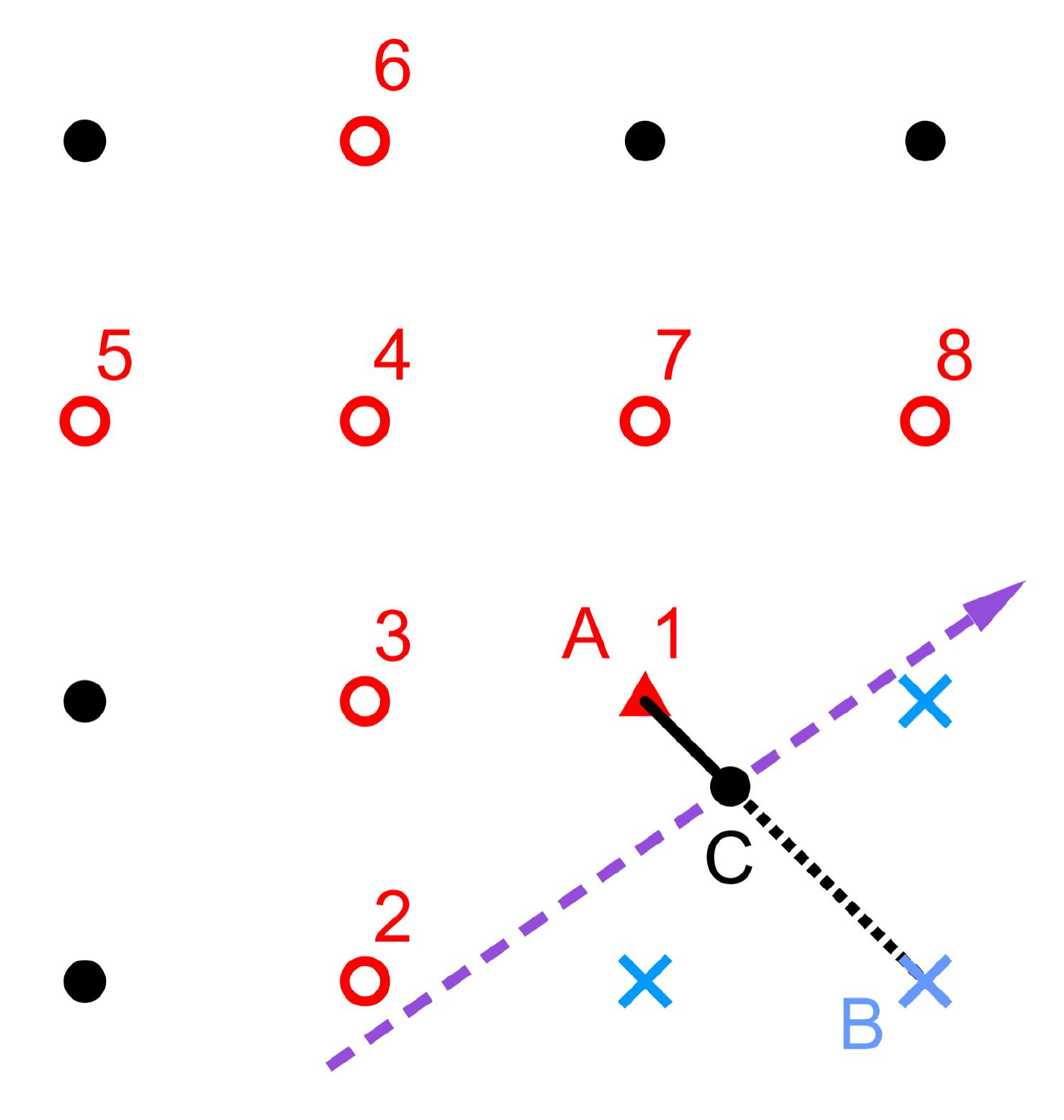}}
        \subfloat[{\tiny{Type 4 for Case V}}]{\includegraphics[width = 0.225 \linewidth, page = 9]{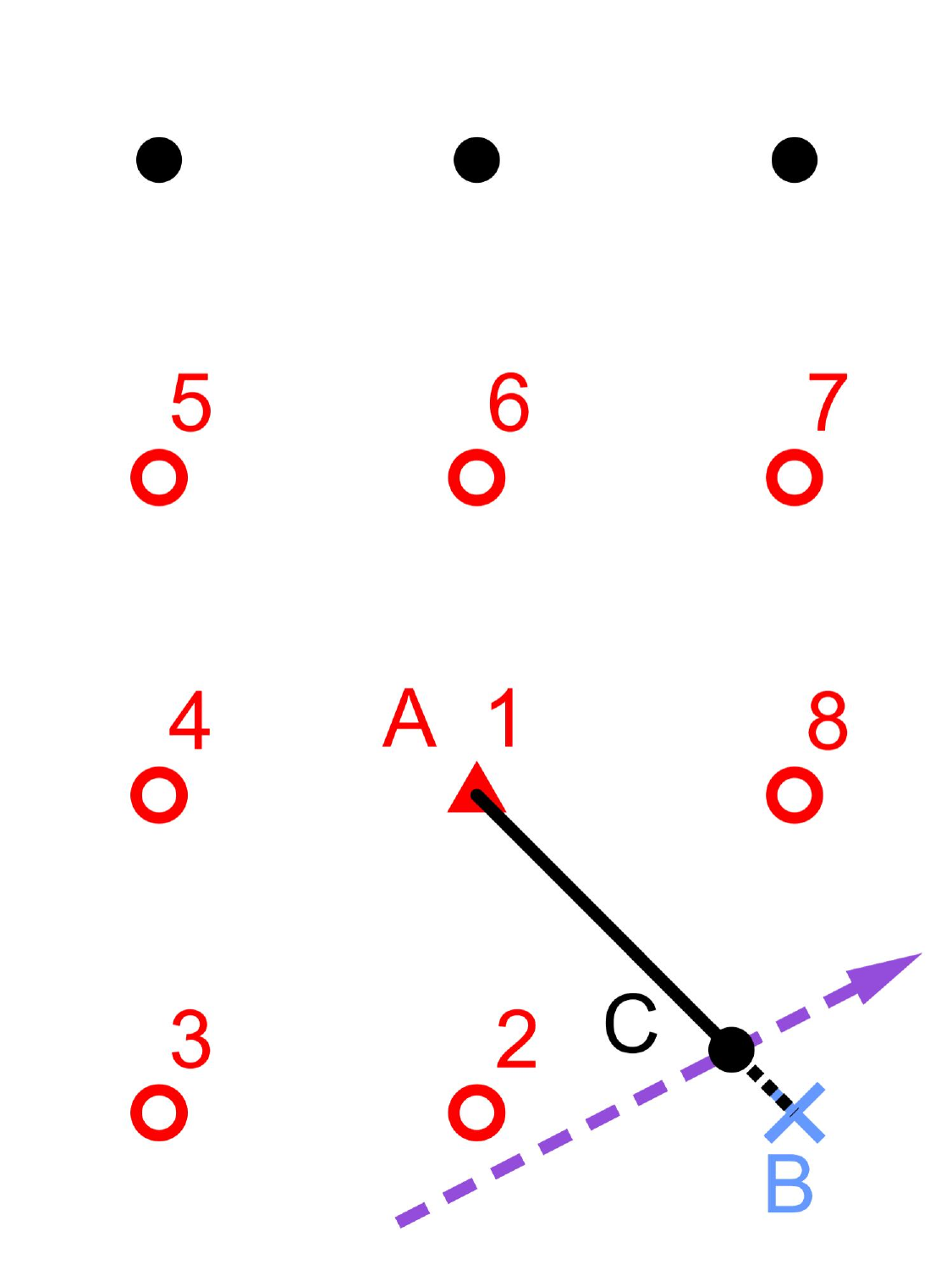}}
        \caption{
        Under the condition $(\spt+h\SS)\cap \hp{L_{\bpt}}\subset \Omega$, four boundary stencil types $\spt+h\SS_{\spt} \subseteq \Omega_h$, consisting of all red grid points ordered by labels
        $1,\ldots, \#\SS_{\spt}$.
        The red triangle is the stencil center $\spt$. Other symbols have the same meaning as in \Cref{fig: inside grid points cases}.
        Cases II and IV in \Cref{fig: inside grid points cases} share the same stencil Type 2 in (B). The points $C\in L_{\bpt}$ will be used in \eqref{extra:eqs} for extra equations.}
        \label{fig: stencil types}
    \end{figure}

    We now consider the $9$ points in $\spt+h\SS$ with $\SS:=[-1,1]^2 \cap \mathbb{Z}^2$ and two cases whether all the points $(\spt+h \SS)\cap \hp{L_{\bpt}}$  belong to $\Omega$ or not.
    If $(\spt+h \SS)\cap \hp{L_{\bpt}}\subset \Omega$, up to flipping and rotation, we have a total of five configurations of $\spt+ h\SS$ with respect to $\Omega$, as illustrated in \Cref{fig: inside grid points cases}. In this case, it is not necessary for us to explicitly indicate the base point $\bpt$ in \Cref{fig: inside grid points cases,fig: stencil types}.

\begin{figure}[h!]
    \centering
        \subfloat[\tiny{Type 5 for Case VI}]{\includegraphics[width = 0.225 \linewidth, page = 10]{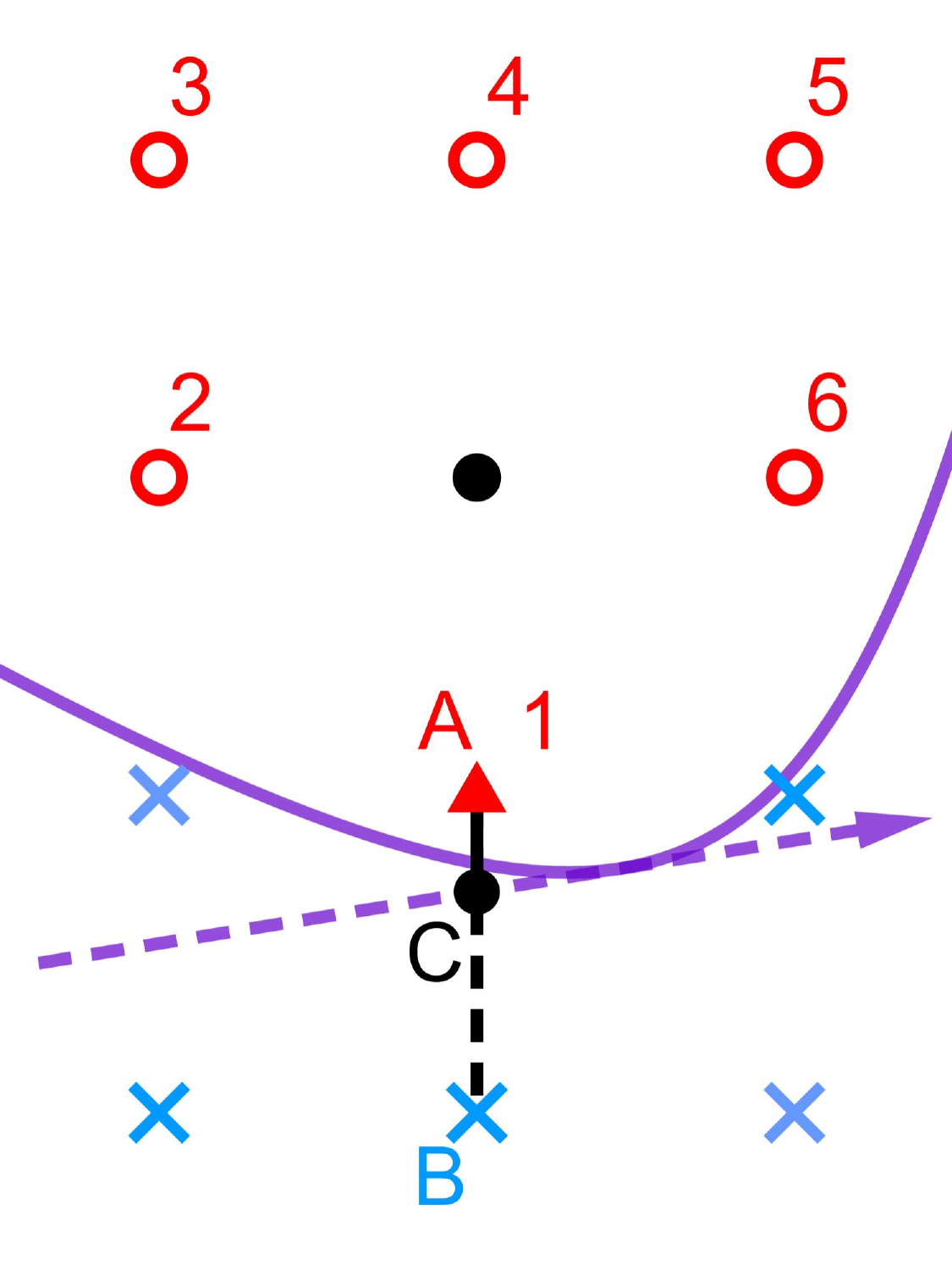}}
        \hspace{0.05\linewidth}
        \subfloat[\tiny{Type 6 for Case VII}]{\includegraphics[width = 0.285 \linewidth, page = 11]{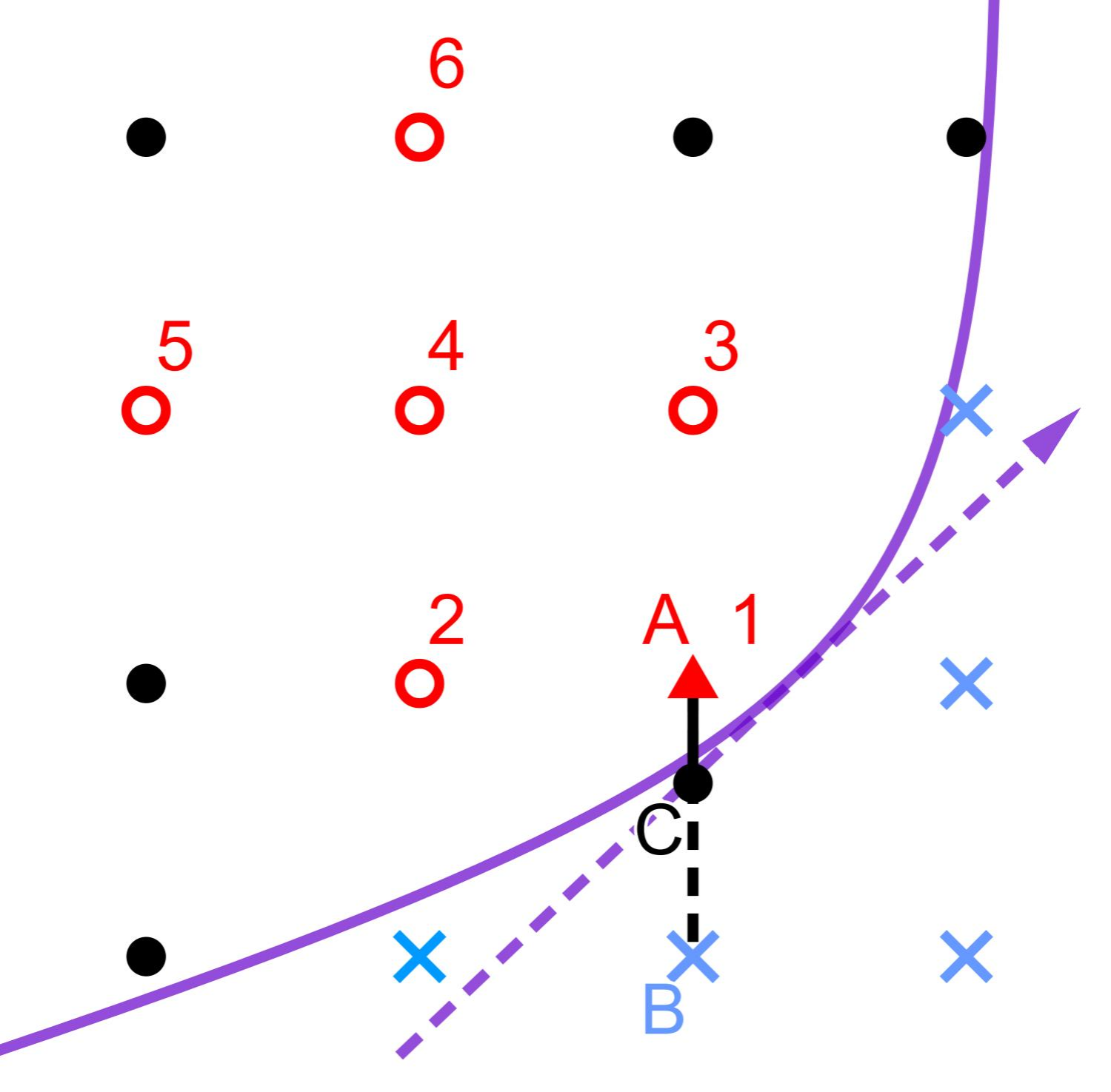}}
        \hspace{0.05\linewidth}
        \subfloat[\tiny{Type 2 for Case VIII}]{\includegraphics[width = 0.225 \linewidth, page = 12]{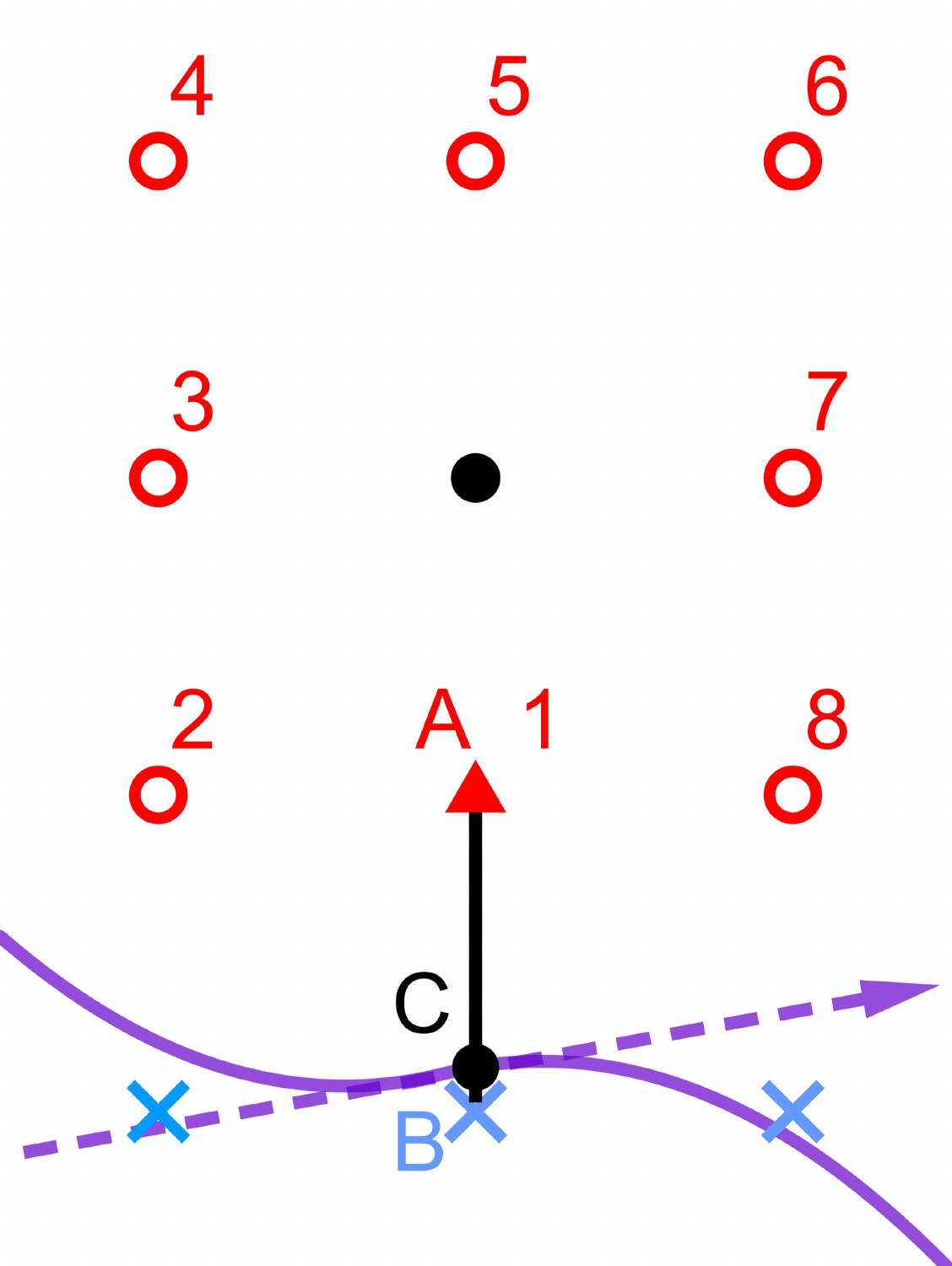}}
        \caption{
        The additional three boundary stencil configurations when the condition $(\spt+h \SS)\cap \hp{L_{\bpt}} \subset \Omega$ fails.
        The directed tangent line $L_{\bpt}$ is the purple dashed line with the arrow indicating the direction. $\hp{L_{\bpt}}$ is the left-hand region of the directed line $L_{\bpt}$.
        The purple solid curve is the actual boundary $\partial \Omega$. The base point $\bpt\in \partial \Omega$ is not shown but on $L_{\bpt}\cap \partial \Omega$.
        The stencil $\SS_{\spt}$ consists of all the red dot grid points in $\Omega_h$.
        All the blue crosses are outside $\Omega$ but may lie on $\partial \Omega$.}
        \label{fig: other inside point cases}
    \end{figure}

    We now discuss how to build a suitable boundary stencil $\SS_{\spt}$ with stencil coefficients having desired properties for consistency order $M + 2 = 6$. When $j = 0$, \cref{eq: boundary constraint of c} is a homogeneous linear system $\mathbb{A}_0 \vec{c}_0 = 0$ of size $5 \times (\#\SS_{\spt})$, whose solution space generally has dimension $(\#\SS_{\spt})-5$. Regardless of the geometry of $\partial \Omega$,
    we could just use the smallest possible $\#\SS_{\spt} = 6$ for all stencil centers $\spt\in \partial \Omega_h$. However, due to the curvature of $\partial \Omega$ near $\spt$, this often leads to many cases of special stencil shapes/configurations $\spt+h \SS_{\spt}\subset \Omega$; consequently, the constructed scheme becomes very complicated to be practically implemented for treating many special cases. As an effort to keep both $\#\SS_{\spt}$ and the number of the special cases of boundary stencil shapes $\SS_{\spt}$ as small as possible, it turns out that we take $\#\SS_{\spt}\in \{6,7,8\}$ depending on the geometry of $\partial \Omega$ near $\spt$ and the tangent line $L_{\bpt}$, and we consider in total only $6$ special cases of stencil shapes showing in \Cref{fig: stencil types,fig: other inside point cases}.

    For the five cases in \Cref{fig: inside grid points cases}, we have a total of four stencil configurations as illustrated in \Cref{fig: stencil types}. To reduce the number of stencil types,  we combine cases II and IV as one configuration by treating the bottom-left dark dot in (D) of \Cref{fig: inside grid points cases}, though inside the domain $\Omega$, as a blue cross in (B) of \Cref{fig: stencil types}. We also select a point $C\in L_{\bpt}$ in \Cref{fig: stencil types} for computing stencil coefficients.

    Once a boundary stencil $\SS_{\spt}$ is selected, we now discuss how to obtain stencil coefficients $c_{p,j}$ satisfying the linear system $\mathbb{A}_j \vec{c}_j = \vec{b}_j$ in \eqref{eq: boundary constraint of c}. This linear system is solved in the order of $j = 0, \ldots, 4$, and inspired by the proof of \Cref{prop: interior c}, we look for \emph{admissible} zeroth-order coefficients $\vec{c}_0$ for proving theoretical convergence later. The admissibility conditions are defined as follows.

    \begin{definition}\label{def:admissible}
    A column vector $\vec{c}_0:=\{c_{p, 0}\}_{p\in \SS_{\spt}}$ is said to be an admissible solution if
    \begin{itemize}[noitemsep]
        \item[(i)] $\vec{c}_0$ is a real-valued solution to
        $\mathbb{A}_0 \vec{c}_0 = \vec{b}_0$ (hence, $\vec{c}_0$ satisfies \eqref{eq: boundary constraint of c} for $j = 0$ and $m=1,\ldots,5$) such that all the coefficients $c_{p, 0}$ for $p\in \SS_{\spt}$ are bounded by a universal constant;
        \item[(ii)] $c_{(0, 0), 0} = 1$ and $c_{p, 0} \leq 0$ for all $p \in \mathring{\SS}_{\spt} := \SS_{\spt} \setminus \{(0, 0)\}$;
        \item[(iii)] $\sum_{p \in \SS_{\spt}} c_{p, 0} \geq \mu_c$ for some positive constant $\mu_c > 0$ independent of $\spt \in \partial \Omega_h$.
    \end{itemize}
    \end{definition}

    To obtain admissible zeroth-order coefficients $\vec{c}_0$, we consider an augmented linear system $\mathbb{A}^*_0 \vec{c}_0 = \vec{b}^*_0$ with an $(\#\SS_{\spt})\times (\#\SS_{\spt})$ matrix $\mathbb{A}^*_0$ by prepending $(\#\SS_{\spt}-5)$ extra linear equations to $\mathbb{A}_0\vec{c}_0=\vec{b}_0$:
    \begin{equation}\label{extra:eqs}
    \mathbb{A}^*_{0;k} \vec{c}_0=\vec{b}^*_0(k),\quad k=1,\ldots, \#\SS_{\spt}-5
    \quad \mbox{with}\quad
    \mathbb{A}^*_{0;1}:=[1,0\ldots,0],\quad \vec{b}^*_0(1):=1.
    \end{equation}
    For each stencil type in \Cref{fig: stencil types,fig: other inside point cases}, we will provide the extra equations in \eqref{extra:eqs} so that the augmented linear system has a unique solution that is admissible and numerically stable. The previous statement will be discussed in detail and verified rigorously in \Cref{app: verification of stencil coefs}.

    For the four stencil types in \Cref{fig: stencil types} with the selected point $C\in L_{\bpt}$, we list the extra $(\#\SS_{\spt}-5)$ linear equations in \eqref{extra:eqs} explicitly in \Cref{table: stencil information}. These extra equations only involve the distance between points $A$ (i.e., the stencil center $\spt$) and $C\in L_{\bpt}$.

    \begin{table}[htb!]
        \begin{center}
        \begin{NiceTabular}{|c|c|c|c|}[cell-space-limits=4pt]
            \hline
            Stencil type & Cases & $\#\SS_{\spt}$ & Extra equations in \eqref{extra:eqs} \\ \hline\hline
            I & 1   & 7 & $\vec{c}_0(2)-\vec{c}_0(7) = -\frac{|\overrightarrow{AC}|}{4h}$ \\ \hline
            \Block{2-1}{II} & 2 & \Block{2-1}{8} & $\vec{c}_0(2) - \vec{c}_0(7) = -\frac{|\overrightarrow{AC}|}{5h}$ \\ \cline{2-2}
                & 4 &  & $\vec{c}_0(8) - \vec{c}_0(3) = -\frac{|\overrightarrow{AC}|}{5h}$ \\ \hline
            III & 3 & 8 & \Block{1-1}{$\vec{c}_0(2) - \vec{c}_0(4) = -\frac{ (1 + |\tau|) |\overrightarrow{AC}|}{15 \sqrt{2} h}$ \\ $\vec{c}_0(8)- \vec{c}_0(4) = -\frac{ (1 + |\tau|)|\overrightarrow{AC}|}{15\sqrt{2} h}$ \\ $\tau = \tan (\theta - \frac{\pi}{4})$} \\ \hline
            IV & 5   & 8 & \Block{1-1}{$\vec{c}_0(2) - \vec{c}_0(7) = 0$ \\ $\vec{c}_0(8) - \vec{c}_0(3) = 0$} \\ \hline
        \end{NiceTabular}
        \end{center}
        \vspace{-9pt}
        \caption{\footnotesize Information on different stencil types (Part I): The corresponding cases (see \Cref{fig: stencil types}), the size of the stencil, and the required extra equation in \eqref{extra:eqs}. The angle $\theta$ is defined in \cref{theta}. The number $|\overrightarrow{AC}|$ above is the distance between points $A$ and $C$ in \Cref{fig: stencil types}.}
        \label{table: stencil information}
    \end{table}

    The above constructions of boundary stencils $\spt+h\SS_{\spt}\subset \Omega_h$ in \Cref{fig: stencil types} require the condition $(\spt+h \SS)\cap \hp{L_{\bpt}} \subset \Omega$. We now consider the case that this condition fails, i.e., we always have some grid points $q\in (\spt+h \SS)\cap \hp{L_{\bpt}}$ but $q\not \in \Omega$. For small enough $h>0$, there are at most two ``trouble'' points $q\in \spt+h \SS$, often near $L_{\bpt}$, such that $q\not\in \Omega$ but $q\in \hp{L_{\bpt}}$. Hence, for small enough $h$, according to the curvature of $\partial \Omega$ at $\bpt$, we have in total three additional cases, as illustrated in \Cref{fig: other inside point cases}. The corresponding selected boundary stencils are also illustrated in \Cref{fig: other inside point cases}, where the base point $\bpt\in \partial \Omega$ is not explicitly given but satisfies $\bpt\in \partial \Omega \cap L_{\bpt}$.

    Using the point $C\in L_{\bpt}$ shown in \Cref{fig: other inside point cases}, the extra $(\#\SS_{\bpt}-5)$ linear equations in \eqref{extra:eqs} are presented in \Cref{table: other stencil information} below. Such extra equations allow us to obtain admissible coefficients $\vec{c}_0$ for sufficiently small $h$. This issue of admissible $\vec{c}_0$ will be fully addressed and proved in \Cref{app: verification of stencil coefs}.

    \begin{table}[h!]
        \begin{center}
        \begin{NiceTabular}{|c|c|c|c|}[cell-space-limits=4pt]
            \hline
            Stencil type & Case & $\#\SS_{\spt}$ & Extra equations in \eqref{extra:eqs} \\ \hline\hline
            5 & VI   & 6 & N/A \\ \hline
            6 & VII  & 6 & N/A \\ \hline
            2 & VIII & 8 & \Block{1-1}{$\vec{c}_0(2)- \vec{c}_0(7) = -\frac{|\overrightarrow{AC}|}{5h}$ \\ $\vec{c}_0(8) - \vec{c}_0(3) = -\frac{|\overrightarrow{AC}|}{5h}$} \\ \hline
        \end{NiceTabular}
        \end{center}
        \vspace{-9pt}
        \caption{\footnotesize Information on different stencil types (Part II): The corresponding cases (see \Cref{fig: other inside point cases}), the size of the stencil, and the required extra equation in \eqref{extra:eqs}. The angle $\theta$ is defined in \cref{theta}. The number $|\overrightarrow{AC}|$ above is the distance between points $A$ and $C$ in \Cref{fig: other inside point cases}.}
        \label{table: other stencil information}
    \end{table}

    When $h$ is not sufficiently small (in particular, when the curvature of $\partial \Omega$ at $\bpt$ is large), it is possible that more points in $(\spt+h \SS)\setminus \Omega$ may belong to $(\spt+h\SS)\cap \hp{L_{\bpt}}$, and hence the above constructed six stencil types will be invalid. In this case, because $h$ is not sufficiently small, we can simply pick $5$ points near $\spt$ from $\Omega_h \cup \partial \Omega$ to solve $\mathbb{A}_0 \vec{c}_0=\vec{b}_0$ without adding any extra equations in \eqref{extra:eqs}. This is because the proof of convergence only deals with small $h\to 0^+$.

    Now we have fixed the boundary stencil, and the following result follows in parallel with \cref{prop: interior c}.

    \begin{prop}\label{prop: boundary c}
        There exists a positive $h_0 = \bo_{\gx, \gy} (1)$ such that for all $0 < h < h_0$, the solution $\vec{c}_0$ to $\mathbb{A}^*_0 \vec{c}_0=\vec{b}^*_0$, which is augmented from $\mathbb{A}_0\vec{c}_0=0$ with extra equations in \eqref{extra:eqs} being stated in \Cref{sec: choice of boundary stencil}, must be real-valued and admissible. Let $\mu_c$ be as in \Cref{def:admissible}.
        Furthermore, for any $0<h<h_0$, there exist real-valued $\vec{c}_j:=\{c_{p,j}\; : \; p\in \SS_{\spt}\}$ for $j=1,\ldots, 4$ such that
        \begin{enumerate}[noitemsep]
            \item[(i)] The real-valued coefficients $\vec{c}_j, j=1,\ldots, 4$ satisfy the linear system \eqref{eq: boundary constraint of c} with $M=4$ and all $c_{p, j} =\bo_{\tilde{a}, \gx, \gy} (1)$. In addition, the equations \eqref{eq: boundary scheme 0} hold for $C_p(h):= \sum_{j = 0}^4 c_{p, j} h^j$, $p\in \SS_{\spt}$ with the remainder term $\bo_{\tilde{a}, \gx, \gy, g}(h^4)$.

            \item[(ii)] For all $j=1,\ldots, 4$, $c_{(0,0), j} \geq 0$ and $c_{p, j} \leq 0$ for all $p \in \mathring{\SS_{\spt}}$;

            \item[(iii)] For all $j=1,\ldots, 4$, $\sum_{p \in \SS_{\spt}} c_{p, j} > 0$ and consequently, $\sum_{p \in \SS_{\spt}} C_p(h) \geq \mu_c>0$.
        \end{enumerate}
    \end{prop}
    
    \begin{proof}
        For each stencil shape, we shall prove in \Cref{app: verification of stencil coefs} the existence and construction of an admissible unique solution $\vec{c}_0$ satisfying $\mathbb{A}^*_0 \vec{c}_0=\vec{b}^*_0$. Then we can further solve \eqref{eq: boundary constraint of c} for the higher-order coefficients $\vec{c}_j$ and use least squares minimization techniques to make the solution unique (mainly to keep the magnitude of stencil coefficients under control). It is easy to see that if the unique solution $\vec{c}_{0}$ to $\mathbb{A}^*_0\vec{c}_0=\vec{b}^*_0$ is admissible by satisfying all conditions in \Cref{def:admissible}, then all the obtained coefficients in $\vec{c}_{j}, j=1,\ldots, 4$ are of order $\bo_{\tilde{a}, \gx, \gy}(1)$ and item (i) holds. After that, we perform a procedure analogous to item (ii) of \Cref{prop: interior c} to modify the higher-order stencil coefficients and achieve properties items (ii) and (iii). For this purpose, one only needs to repeat the proof of \Cref{prop: interior c} and replace $\lambda_j$ in \cref{q:modcoeff} by
        \begin{equation}
            \label{eq: lambda_m boundary}
            \tilde{\lambda}_j := \max \left\{ \lambda_j, -\frac{\sum_{p \in \SS_{\spt}} c_{p, j}} {\sum_{p \in \SS_{\spt}} c_{p, 0}} \right\},
            \quad j=1,\ldots,4.
        \end{equation}
        In summary, the fact $\tilde{\lambda}_j \geq \lambda_j$ will guarantee that item (ii) is true, and the second term in the definition of $\tilde{\lambda}_j$ guarantees item (iii). See the proof of \Cref{prop: interior c} for the detailed argument.
    \end{proof}

%%%%%%%%%%%%%%%%%%%%%%%%%%%%%%%%%%%%%%%%%%%%%%
%%%%%%%%%%%%%%%%%%%%%%%%%%%%%%%%%%%%%%%%%%%%%%
    \section{The Sixth-order Convergence of the Numerical Solution and Gradient $\nabla u$}
    \label{sec: convergence}

    For our proposed FDM scheme, in this section we rigorously prove the sixth-order convergence of the numerically approximated solution $u_h$ in the $\infty$-norm. Then we shall derive a gradient approximation $\nabla u$ directly from $u_h$ without solving auxiliary equations. Finally, we prove that the gradient approximation achieves a superconvergence of
    order $5+\frac{1}{q}$ in the $q$-norm for all $1 \le q \le \infty$ (with a logarithmic factor $\log h$ for $1 \le q < 2$).

    \subsection{Sixth-order convergence of the numerically approximated solution $u_h$}
    \label{sec: convergence of u}

    In \Cref{sec: interior,sec: boundary} we have described in detail the construction of the FDM scheme at interior and boundary grid points. We have spent much effort on proving the admissibility of the solution $\vec{c}_0$ to $\mathbb{A}^*_0 \vec{c}_0=\vec{b}^*_0$ in \Cref{app: verification of stencil coefs}. This then leads to \Cref{prop: interior c,prop: boundary c} on extra properties of the stencil coefficients. Then we shall use these properties to prove that our proposed scheme achieves sixth-order convergence. 

    We begin by explicitly stating the assumptions on the bounded domain $\Omega$ and various functions in the model problem \eqref{eq:PDE}.
    \begin{itemize}
        \item For every $\bpt \in \partial \Omega$, there exists a local parametrization in \eqref{bdcurve} with $(\gx(t^*),\gy(t^*))=\bpt$ such that $\gx$ and $\gy$ have continuous derivatives of order up to six, and $(\gx' (t^*), \gy' (t^*)) \neq (0, 0)$.

        \item The exact solution $u \in C^8 (\overline{\Omega})$, the data functions $a, f \in C^6 (\overline{\Omega})$, and boundary $g \in C^6 (\partial \Omega)$.

        \item The diffusion coefficient $a$ satisfies $\inf_{x \in \overline{\Omega}} a(x) > 0$.
    \end{itemize}
    Such regularity is needed for performing Taylor expansion is various places. In addition, at least $C^1$ boundary is needed for the construction of the boundary stencils and the admissibility condition. According to \Cref{prop: boundary c}, we assume that $0<h < h_0$ throughout this section. Our main result on convergence is as follows.
    \begin{thm}
        \label{thm: convergence}
        Let $u$ be the exact solution to the model problem \eqref{eq:PDE}, and let $u_h$ be the numerically approximated solution by solving the linear system in \eqref{eq: overall scheme}. Then there exist $0 < h_1 \leq h_0$ and a positive constant $C$ such that
        \begin{equation}
            \label{eq: convergence}
            \|u - u_h\|_{L^\infty (\Omega_h)} \le C h^6,\qquad \forall\; 0 < h < h_1,
        \end{equation}
        where the positive constant $C = \bo_{a, \gx, \gy, u}(1)$, i.e., the constant $C$ only depends on the diffusion coefficient $a$, the exact solution $u$ and the boundary curve $\partial \Omega$.
    \end{thm}

    The proof of \Cref{thm: convergence} will be presented at the end of this subsection. To prove \Cref{thm: convergence}, we shall follow a slightly modified traditional method by using the discrete maximum principle to prove the sixth-order convergence of our proposed FDM.

    Recall that $\Omega_h, \Omega^\circ_h$ and $\partial \Omega_h$ are defined in \eqref{eq: Omegah} and $\Omega^\circ_h\cup \partial \Omega_h=\Omega_h=\Omega\cap (h\Z^2)$. We define the difference operator $\mathcal{L}_h$ acting on any grid function $v_h: \Omega_h\rightarrow \R$ by
    \begin{equation}
        \label{eq: Lh}
        \mathcal{L}_h v_h(\spt) := h^{-\sigma} \sum_{p \in \SS_{\spt}} C_p(h) v_h(\spt + ph) \quad \mbox{with}\quad \sigma:=
        \begin{cases}
            2, & \mbox{if} \ \spt \in \Omega_h^\circ, \\
            0, & \mbox{if} \ \spt \in \partial \Omega_h.
        \end{cases}
    \end{equation}
    Here $\SS_{\spt} = \SS:=[-1,1]^2\cap \Z^2$ for $\spt \in \Omega_h^\circ$, and $C_p(h)$ are the real-valued stencil coefficients in Propositions~\ref{prop: interior c} or \ref{prop: boundary c} depending on $\spt\in \Omega_h^\circ$ for interior stencils or $\spt\in \partial \Omega_h$ for boundary stencils. The FDM scheme in \Cref{sec: interior,sec: boundary} to the model problem \eqref{eq:PDE} can be expressed as
    \begin{equation}
        \label{eq: overall scheme}
        \mathcal{L}_h u_h = f_h
        \quad \mbox{with} \quad f_h :=
        \begin{cases}
            h^{-2} \sum_{p \in \SS} C_p(h) F(p) & \mbox{on} \ \Omega_h^\circ, \\
            \sum_{p \in \SS_{\spt}} C_p(h) G(\psv) & \mbox{on} \ \partial \Omega_h,
        \end{cases}
    \end{equation}
    where the real-valued quantities $F(p)$ and $G(\psv)$ are defined in \eqref{Fp} and \eqref{Akm:Gm}. Despite the use of complex partial derivatives in deriving this FDM, we eventually obtain a real-valued linear system for the numerically approximated solution $u_h$, which also guarantees that $u_h$ is real-valued.

    According to \Cref{sec: choice of boundary stencil}, when $0<h < h_0$, the stencil at a boundary grid point does not include points on the true boundary $\partial \Omega$. Hence, we can treat $\mathcal{L}_h$ as a linear mapping on the space $L^\infty (\Omega_h)$. Moreover, \Cref{prop: interior c,prop: boundary c} guarantee
    \begin{equation}
        \label{eq: consistency}
        \| \mathcal{L}_h u - f_h \|_{L^\infty (\Omega_h)} = \bo_{\tilde{a}, \gx, \gy, u} (h^6).
    \end{equation}

    \begin{thm}
        \label{thm: maximum principle}
        Assume that $0<h < h_0$. Let $v_h$ be a grid function defined on $\Omega_h$ such that $\mathcal{L}_h v_h \geq 0$. Then $v_h$ takes its minimum in $\partial \Omega_h$, and its minimum must be nonnegative.
    \end{thm}

    \begin{proof}
        Suppose $v_h$ takes its minimum at $\spt \in \Omega_h^\circ$. By \Cref{prop: interior c}(ii), the interior stencil coefficients satisfy $C_p(h) < 0$ for all $p \in \mathring{\SS} = \SS_{\spt} \bs \{ (0, 0) \}$, and $\sum_{p \in \SS_{\spt}} C_p(h) = 0$. Thus,
        \begin{equation}\label{discretemp}
            0 \leq \mathcal{L}_h v_h (\spt) = \sum_{p \in \SS_{\spt}} C_p(h) v_h(\spt + ph)
            \leq \sum_{p \in \SS_{\spt}} C_p(h) v_h(\spt).
        \end{equation}
        Because $\sum_{p \in \SS_{\spt}} C_p(h) = 0$, the above inequalities imply that all inequalities in \eqref{discretemp} must be equalities.
        Hence, we conclude from \eqref{discretemp} and $C_p(h) < 0$ for all $p \in \mathring{\SS}$ that $v_h(\spt + ph) = v_h(\spt)$ for all $p \in \SS_{\spt}$. Consequently, $v_h$ must take its minimum on $\partial \Omega_h$.

        Now let $\spt \in \partial \Omega_h$ be the minimum point of $v_h$. By \Cref{prop: boundary c}, we have $C_p(h) \le 0$ for all $p\in \SS_{\spt}$ but $p \neq (0, 0)$, and $\sum_{p \in \SS_{\spt}} C_p(h)\ge \mu_c>0$. Note that \eqref{discretemp} is still true in this case. It follows from \eqref{discretemp} and $\sum_{p \in \SS_{\spt}} C_p(h)>0$ that $v_h(\spt) \geq 0$. So, the minimum of $v_h$ must be nonnegative.
    \end{proof}

    \begin{lemma}
        \label{lem: aux function}
        There exists a real-valued function $\phi$ in $\overline{\Omega}$ such that $\| \phi \|_{L^\infty (\overline{\Omega})} = \bo_{\tilde{a}, \gx, \gy} (1)$, $\|\mathcal{L}_h \phi-1\|_{L^\infty(\Omega_h^\circ)}=
        \bo_{\tilde{a}, \gx, \gy} (h)$, and $\mathcal{L}_h \phi \ge 1$ on $\partial\Omega_h$
         for all $0<h<h_0$.
    \end{lemma}

    \begin{proof}
        Fix a function $\tilde{\phi}$ on $\overline{\Omega}$ such that $-\nabla \cdot (a \nabla \tilde{\phi}) = a$, or equivalently, $\Delta \tilde{\phi} = \nabla \tilde{a} \cdot \nabla \tilde{\phi} - 1$ with $\tilde{a}:=-\ln a$. By elliptic regularity theory (e.g., \cite[Chapter~6]{gilbarg1977elliptic}), the derivatives of $\tilde{\phi}$ are bounded by the derivatives of $\tilde{a}$ and $\|\tilde{\phi}\|_{L^\infty(\overline{\Omega})}=\bo_{\tilde{a}, \gx, \gy} (1)$. As an analog of \cref{eq: consistency}, we have
        \begin{equation*}
            \Big\| \mathcal{L}_h \tilde{\phi} - h^{-2} \sum_{p \in \SS_{\spt}} C_p(h) F^\phi(p) \Big\|_{L^\infty (\Omega_h^\circ)}
            = \bo_{\tilde{a}, \gx, \gy, \tilde{\phi}} (h^6) = \bo_{\tilde{a}, \gx, \gy} (h^6),
        \end{equation*}
        where $F^\phi(p)$ is obtained by replacing $\tilde{f}$ with $-1$ in $F(p)$. By symbolic calculation, we can obtain
        \begin{equation*}
            h^{-2} \sum_{p \in \SS_{\spt}} \tilde{C}_p(h) F^\phi(p) = 6 + \bo_{\tilde{a}} (h^2),
        \end{equation*}
        where $\tilde{C}_p(h):=\sum_{j=0}^7 \vec{c}_p(j) h^j$ with the column vectors $\vec{c}_{p}$ given in Appendix \ref{app: Interior stencil coefficients}. According to the proof of \Cref{prop: interior c}, there exists a polynomial $Q(h) = 1 + \bo_{\tilde{a}}(h)$ such that $C_p(h) = \tilde{C}_p(h) Q(h) + \bo_{\tilde{a}} (h^8)$ holds for each $p\in \SS_{\spt}$. Therefore, we have
        \begin{equation*}
            h^{-2} \sum_{p \in \SS_{\spt}} C_p(h) F^\phi(p) = 6 + \bo_{\tilde{a}} (h),
        \end{equation*}
        which implies $\| \mathcal{L}_h \tilde{\phi} - 6 \|_{L^\infty (\Omega_h^\circ)} = \bo_{\tilde{a}, \gx, \gy} (h)$.

        For the boundary case, as an analog of \cref{eq: consistency}, we have
        \begin{equation}
            \label{eq: L_h2 tilde_phi}
            \Big\| \mathcal{L}_h \tilde{\phi} - \sum_{p \in \SS_{\spt}} C_p(h) G^\phi(\psv) \Big\|_{L^\infty (\partial \Omega_h)}
            = \bo_{\tilde{a}, \gx, \gy, \tilde{\phi}} (h^6) = \bo_{\tilde{a}, \gx, \gy} (h^6),
        \end{equation}
        where $G^\phi(\psv)$ is obtained by replacing $\tilde{f}$ and $g$ with $-1$ and $\tilde{\phi} \big|_{\partial \Omega}$ in $G(\psv)$. Clearly $G(\psv) = \bo_{\tilde{f}, g, \tilde{a}, \gx, \gy} (1)$, which implies $G^\phi(\psv) = \bo_{\tilde{a}, \gx, \gy} (1)$. In view of \eqref{eq: L_h2 tilde_phi}, we get $\mathcal{L}_h \tilde{\phi} = \bo_{\tilde{a}, \gx, \gy} (1)$. Hence, we proved
        \begin{equation}\label{tildephi}
        \| \mathcal{L}_h \tilde{\phi} - 6 \|_{L^\infty (\Omega_h^\circ)} = \bo_{\tilde{a}, \gx, \gy} (h)
        \quad \mbox{and}\quad
        \| \mathcal{L}_h \tilde{\phi} \|_{L^\infty (\partial \Omega_h)} = \bo_{\tilde{a}, \gx, \gy} (1).
        \end{equation}

        Define $M_0 := \frac{1}{6\mu_c} \| \mathcal{L}_h \tilde{\phi} \|_{L^\infty (\partial \Omega_h)}$, where the positive constant $\mu_c$ is as in \Cref{def:admissible} and \Cref{prop: boundary c}(iii). Then
        $M_0=\bo_{\tilde{a},\gx,\gy}(1)$ by \eqref{tildephi}.
        Consider $\phi := \frac{1}{6} \tilde{\phi} + M_0$ in $\overline{\Omega}$. Noting that $\sum_{p \in \SS_{\spt}} C_p(h) = 0$ in item (iii) of \Cref{prop: interior c} for all $\spt\in \Omega_h^\circ$, we must have
        \begin{equation*}
            \mathcal{L}_h \phi(\spt)=\frac{1}{6}\mathcal{L}_h \tilde{\phi}(\spt)+M_0 \mathcal{L}_h 1=\frac{1}{6}\mathcal{L}_h \tilde{\phi}(\spt)
        \end{equation*}
        for all $\spt\in \Omega_h^{\circ}$. Now it follows directly from the first identity in \eqref{tildephi} that
        \begin{equation*}
            \|\mathcal{L}_h \phi-1\|_{L^\infty(\Omega_h^\circ)}
            = \left\| \frac{1}{6}\mathcal{L}_h \tilde{\phi}-1 \right\|_{L^\infty(\Omega_h^\circ)}
            =\bo_{\tilde{a}, \gx, \gy} (h).
        \end{equation*}
        On the other hand, for every $\spt \in \partial \Omega_h$, noting that $\mu_c \leq \sum_{p \in \SS_{\spt}} C_p(h)=\bo_{\tilde{a}, \gx, \gy}(1)$ by \Cref{prop: boundary c}, we have
        \[
            \mathcal{L}_h \phi(\spt)
            = \frac{1}{6} \mathcal{L}_h \tilde{\phi} + \sum_{p \in \SS_{\spt}} C_p(h) M_0
            \ge \mu_c M_0-\frac{1}{6} \|\mathcal{L}_h \tilde{\phi}\|_{L^\infty(\overline{\Omega})} \ge 1.
        \]
        This proves $\mathcal{L}_h \phi(\spt)\ge 1$ for all $\spt\in \partial \Omega_h$.
    \end{proof}

    We are now ready to prove \Cref{thm: convergence}.

    \begin{proof}[Proof of \Cref{thm: convergence}]
     Recall that $\tilde{a}:=-\ln a$ is defined in \eqref{pde:mod}.
        Let $\phi$ be the auxiliary function in \Cref{lem: aux function}. By $\|\mathcal{L}_h \phi-1\|_{L^\infty(\Omega_h^\circ)} = \bo_{\tilde{a},\gx,\gy}(h)$ and $\mathcal{L}_h \phi\ge 1$ on $\partial \Omega_h$ in \Cref{lem: aux function}, there exists $h_1\in (0,h_0)$ such that $\mathcal{L}_h \phi\ge 1/2$ on $\Omega_h$ for all $0 < h < h_1$.

        We first prove that the linear operator $\mathcal{L}_h: L^\infty(\Omega_h)\rightarrow L^\infty(\Omega_h)$, defined in \eqref{eq: Lh}, must satisfy
        \begin{equation}\label{Lh:bound}
            \left\| \mathcal{L}_h^{-1} \right\|_{L^\infty (\Omega_h)}\le 2 \|\phi\|_{L^\infty(\overline{\Omega})},\qquad \forall\; 0 < h < h_1.
        \end{equation}
        Let $w_h$ be any grid function on $\Omega_h$ and define another grid function
        \[
            v_h := 2 M_w \phi + w_h
            \quad \mbox{with}\quad
            M_w:=\|\mathcal{L}_h w_h \|_{L^\infty (\Omega_h)}.
        \]
        Then $\mathcal{L}_h v_h = 2 M_w \mathcal{L}_h \phi + \mathcal{L}_h w_h\ge M_w + \mathcal{L}_h w_h \geq 0$ for all $0 < h < h_1$, due to $\mathcal{L}_h \phi\ge 1/2$. By \Cref{thm: maximum principle}, we must have $\min_{\Omega_h} v_h \geq 0$. Similarly, consider $v_h = 2 M_w \phi - w_h$ instead. Then the same argument shows that $\mathcal{L}_h v_h\ge 0$ and $\min_{\Omega_h} v_h \geq 0$ holds. Consequently, we proved $2M_w\phi \pm w_h\ge 0$ on $\Omega_h$ and hence
        \begin{equation*}
            \|w_h\|_{L^\infty (\Omega_h)}
            \leq 2M_w \| \phi \|_{L^\infty (\overline{\Omega})}
            =
            2 \| \phi \|_{L^\infty (\overline{\Omega})} \cdot \|\mathcal{L}_h w_h\|_{L^\infty (\Omega_h)}, \quad \forall \, 0 < h < h_1
        \end{equation*}
        for all grid functions $w_h$ on $\Omega_h$. This proves that $\mathcal{L}_h^{-1}$ is bounded and satisfies \eqref{Lh:bound}.

        Note that $\mathcal{L}_h u_h=f_h$. By the consistency in \eqref{eq: consistency} and the boundedness of $\mathcal{L}_h^{-1}$ in \eqref{Lh:bound}, we have
        \[
            \|u - u_h\|_{L^\infty (\Omega_h)}
            =\|\mathcal{L}_h^{-1} (\mathcal{L}_h u - f_h)\|_{L^\infty (\Omega_h)}
            \le \|\mathcal{L}_h^{-1}\|_{L^\infty(\Omega_h)}
            \|\mathcal{L}_h u - f_h\|_{L^\infty (\Omega_h)}=\bo_{\tilde{a},\gx,\gy,u}(h^6),
        \]
        where we also used $\| \phi \|_{L^\infty (\overline{\Omega})} = \bo_{\tilde{a}, \gx, \gy} (1)$ in \Cref{lem: aux function}.
        This proves \eqref{eq: convergence} for all $0 < h < h_1$.
    \end{proof}

    \subsection{A high-order approximation of $\nabla u$}
    \label{sec: nabla u_h}

    In this section, we derive a fifth-order accurate approximation in the $\infty$-norm of the gradient $\nabla u$ from the numerically approximated solution $u_h$ without solving additional equations. For any stencil centered at $\spt \in \Omega_h$ with its associated base point $\bpt$, we perform a local approximation for $\grad u(\bpt)$ using the already computed numerical solution $u_h$ from a set of points $\spt + ph$, $p \in \widehat{\SS} \supseteq \SS_{\spt}$. In this process we do not need to solve any linear system to obtain the approximated gradient. In the next subsection, we prove that this gradient approximation exhibits a suboptimal sixth-order superconvergence in the $1$-norm.

    We first discuss the case when $\spt \in \Omega_h^\circ$. Note that $\bpt=\spt$, i.e., the base point $\bpt$ agrees with the stencil center $\spt$. To approximate $\nabla u(\bpt)$, it is sufficient to look at how $\partial_x u(\bpt) = \frac{1}{2} \partial_{\C}^{(1, 0)} u(\bpt) + \frac{1}{2} \partial_{\C}^{(0, 1)} u(\bpt)$ is approximated. As an analog of \cref{eq: interior scheme 0}, we look for a set of real-valued coefficients $C_p(h)$, $p \in \widehat{\SS}$ such that
    \begin{equation}
        \label{eq: interior nabla u}
        \sum_{p \in \widehat{\SS}} C_p(h) u(\spt + ph) = h \left( \partial_{\C}^{(1, 0)} u(\bpt) + \partial_{\C}^{(0, 1)} u(\bpt) \right) + \sum_{p \in \widehat{\SS}} C_p(h) F(p) + \bo_{\tilde{a}, u} (h^{M + 2})
    \end{equation}
    holds, where $F(p)$ is defined in \eqref{Fp}. This is equivalent to computing
    \begin{equation}
        \label{eq: interior nabla u rearranged}
        \partial_x u(\bpt)
        = \frac{1}{h} \sum_{p \in \widehat{\SS}} C_p(h) u(\spt + ph) - \frac{1}{h} \sum_{p \in \widehat{\SS}} C_p(h) F(p) + \bo_{\tilde{a}, u} (h^{M + 1}).
    \end{equation}
    Since the numerical solution $u_h$ satisfies $\|u-u_h\|_{L^\infty(\Omega)}=\bo_{\tilde{a},\gx,\gy,u}(h^6)$ by \Cref{thm: convergence}, we obtain
    \begin{equation*}
        \partial_x u(\bpt)
        = \frac{1}{h} \sum_{p \in \widehat{\SS}} C_p(h) u_h(\spt + ph) - \frac{1}{h} \sum_{p \in \widehat{\SS}} C_p(h) F(p) + \bo_{\tilde{a}, \gx, \gy, u} (h^{\min\{ M + 1, 5 \}}).
    \end{equation*}
    Thus, as long as the stencil $\widehat{\SS}$ and the coefficients $C_p(h)$ are known, we can use the right-hand side of the above identity to approximate $\partial_x u(\bpt)$ with the accuracy order $\min(M+1,5)$.

    In the same way as \Cref{lem: interior constraint of c}, we can prove that the coefficients $C_p(h) := \sum_{j = 0}^{M + 1} c_{p, j} h^j \in \R$ satisfy \eqref{eq: interior nabla u} if and only if
    \begin{equation}
    \label{eq: interior nabla constraints}
        \begin{aligned}
            \sum_{p \in \widehat{\SS}} A^{m+n}_{m,n}(p) c_{p, j}
            &= \td(j) \big( \td(m - 1) \td (n) + \td(m) \td (n - 1) \big) - \sum_{k = 0}^{j - 1} \sum_{p \in \widehat{\SS}} A^{m + n+j- k}_{m, n}(p) c_{p, k}, \\
            & \qquad \qquad \forall \, j=0,\ldots, M + 1, \ (m, n) \in \LN^{M+1-j}_{M + 1 - j},
        \end{aligned}
    \end{equation}
    and we can take the real and imaginary parts to get a real linear system. If we choose $\widehat{\SS} = \SS_{\spt} = [-1,1]^2 \cap \Z^2$, then the maximum possible $M$ is 3, that is, the original stencil $\spt + ph$, $p \in \mathcal{S}$ only yields at most fourth-order accurate numerical $\nabla u$. To reach the maximum potential of fifth order, we can choose $\widehat{\SS} = \SS \cup \{ (\pm 2, 0), (0, \pm 2) \}$ and consider only the grid points $\spt$ such that $\spt + h \widehat{\SS}\subset \Omega_h$. In each of these two cases, we present one particular set of coefficients $C_p(h)$ satisfying \cref{eq: interior nabla constraints} in Appendix \ref{app: coefficients partial_x}.

    Now we consider $\spt \in \partial \Omega_h$. Note that $\bpt\in \partial \Omega$ as in \eqref{vec:shift} and we have an exact formula for $\re (e^{\ii \theta} \partial_{\C}^{(1, 0)} u(\bpt))$ in \cref{Re:Im:u}. Then we can approximate $\nabla u(\bpt)$ by using $\im (e^{\ii \theta} \partial_{\C}^{(1, 0)} u(\bpt))$ according to the identities
    \begin{align*}
        & \partial_x u(\bpt) = \frac{\cos \theta}{|z_0(t^*)|} \dd{t} \tilde{g}(t^*) + 2 \sin \theta \im (e^{\ii \theta} \partial_{\C}^{(1, 0)} u(\bpt)), \\
        & \partial_y u(\bpt) = \frac{\sin \theta}{|z_0(t^*)|} \dd{t} \tilde{g}(t^*) - 2 \cos \theta \im (e^{\ii \theta} \partial_{\C}^{(1, 0)} u(\bpt)).
    \end{align*}

    The way to approximate $\im (e^{\ii \theta} \partial_{\C}^{(1, 0)} u(\bpt))$ is the same as the interior case. In summary,
    \begin{equation*}
        \im (e^{\ii \theta} \partial_{\C}^{(1, 0)} u(\bpt))
        = \frac{1}{h} \sum_{p \in \widehat{\SS}} C_p(h) u_h(\spt + ph) - \frac{1}{h} \sum_{p \in \widehat{\SS}} C_p(h) G(\psv) + \bo_{\tilde{a}, \gx, \gy, u} (h^{\min\{ M + 1, 5 \}}),
    \end{equation*}
    where $G(p)$ is defined in \cref{Akm:Gm} and $C_p(h) := \sum_{j = 0}^M c_{p, j} h^j \in \R$ satisfies the linear system
    \begin{align*}
        -\sum_{p \in \SS_{\spt}} A^m_m(\psv) c_{p, j}
        = -\td(j) \td(m - 1) + \sum_{k = 0}^{j - 1} \sum_{p \in \SS_{\spt}} A^{m + j - k}_m(\psv) c_{p, k}, \\
        \qquad \qquad \forall \, j=0, \ldots, M, \ m = 1, \ldots, M + 1 - j.
    \end{align*}
    We take $M = 4$ and $\widehat{\SS} = \SS_{\spt}$. Note that the above equation only differs from \cref{eq: boundary constraint of c} on the right-hand side. According to \Cref{sec: choice of boundary stencil}, this linear system is away from being singular, so the coefficients $C_p(h)$ always exist and are bounded. Therefore, we can achieve fifth-order accurate approximation in the $\infty$-norm of the gradient $\nabla u$ from the numerical solution $u_h$.

    \subsection{Superconvergence of numerical gradient $\nabla u$}
    \label{sec: superconvergence}

    Denote the numerical gradient by $\nabla u_h = (\partial_x u_h, \partial_y u_h)$ as we discussed in \Cref{sec: nabla u_h}. Besides, we define
    \begin{equation}
        \label{eq: norm}
        \|v_h\|_{L^q_h (U_h)} := \bigg( \frac{1}{\# U_h} \sum_{\spt \in U_h} |v_h(\spt)|_\infty^q \bigg)^{1/q}, \ 1 \leq q < \infty, \quad
        \|v_h\|_{L^\infty_h (U_h)} := \sup_{\spt \in U_h} |v_h(\spt)|_\infty.
    \end{equation}
    Here $U_h$ is a finite subset of $\overline{\Omega}$ with $\# U_h$ elements, $v_h: U_h \to \R^n$ is any grid vector function and $|\cdot|_q$ stands for the $\ell^q$ norm of a vector. We shall use the second set of stencil coefficients (denoted by $C(p)$) in Appendix \ref{app: coefficients partial_x} to approximate $\partial_x u$ at interior grid points, which satisfies \cref{eq: interior nabla u rearranged} with $M = 5$. Define $\widehat{\Omega}_h$ to be the set of all associated base points $\bpt$ such that $\spt + \widehat{\SS} h \subset \Omega_h$, where $\widehat{\SS} = \SS \cup \{ (\pm 2, 0), (0, \pm 2) \}$ is the extended stencil in \Cref{sec: nabla u_h}.
    Note that we can only evaluate $\nabla u_h$ on the set $\widehat{\Omega}_h$.

    The main theorem is stated as follows. We shall first establish some necessary auxiliary results and then we prove \Cref{thm: superconvergence} in detail at the end of this subsection.

    \begin{thm}
        \label{thm: superconvergence}
        Let $u$ be the exact solution to the model problem \eqref{eq:PDE}, let $u_h$ be the numerically approximated solution by solving the linear system in \eqref{eq: overall scheme}, and denote $\nabla u_h$, $\widehat{\Omega}_h$ as above. Then
        \begin{equation}
            \label{eq: superconvergence}
            \|\nabla u - \nabla u_h\|_{L^q_h (\widehat{\Omega}_h)} = \bo_{\tilde{a}, \gx, \gy, u} \big( h^{5 + 1/q} (\log h)^{ \max\{2/q - 1, 0\} } \big), \quad \forall \, 1 \leq q \leq \infty.
        \end{equation}
    \end{thm}

    If $U_h \subseteq \Omega_h$, we can define
    \begin{equation*}
        U_h^\circ := \{ \spt \in U_h \setsp \spt + h\SS \subset U_h \}, \quad \partial U_h = U_h \bs U_h^\circ.
    \end{equation*}
    This aligns with the definition in \eqref{eq: Omegah}.
    We further define the discrete derivatives as
    \begin{equation*}
        \partial_p v_h(\spt) := \frac{1}{|p|_2 h} \left( v_h(\spt + ph) - v_h(\spt) \right), \ p \in \mathring{\SS}, \quad
        \nabla_h v_h(\spt) := (\partial_p v_h(\spt))_{p \in \mathring{\SS}},
    \end{equation*}
    where $v_h: U_h \to \R$ and $\spt \in U_h^\circ$.

    \begin{lemma}
        \label{lem: summation by parts}
        For any subset $U_h$ of $\Omega_h$, any $p \in \mathring{\SS}$ and any grid functions $v_h$, $w_h$ on $\Omega_h$, we have
        \begin{equation}
            \label{eq: summation by parts}
            \left| \langle \partial_p v_h, w_h \rangle_{L^2_h (U_h^\circ)} - \langle v_h, \partial_{-p} w_h \rangle_{L^2_h (U_h^\circ)} \right|
            \leq M_0 \|v_h\|_{L^\infty_h (\partial U_h \cup \partial U_h^\circ)} \|w_h\|_{L^\infty_h (\partial U_h \cup \partial U_h^\circ)},
        \end{equation}
        where $M_0 = \frac{\# \partial U_h + \# \partial U_h^\circ}{\# U_h^\circ \cdot |p|_2 h}$. In the case of $U_h = \Omega_h$, we have $M_0 = \bo_{\Omega}(1)$.
    \end{lemma}

    \begin{proof}
         Let $U_h^\circ \Delta (U_h^\circ - ph)$ be the symmetric difference of the sets $U_h^\circ$ and $U_h^\circ - ph$. Then
        \begin{align*}
            \langle \partial_p v_h, w_h \rangle_{L^2_h (U_h^\circ)}
            &= \frac{1}{\# U_h^\circ \cdot |p|_2 h} \sum_{\spt \in U_h^\circ} \left( v_h(\spt + ph) - v_h(\spt) \right) w_h(\spt) \\
            &= \frac{1}{\# U_h^\circ \cdot |p|_2 h} \left( \sum_{\spt \in U_h^\circ + ph} v_h(\spt) w_h(\spt - ph) - \sum_{\spt \in U_h^\circ} v_h(\spt) w_h(\spt) \right) \\
            &= \langle v_h, \partial_{-p} w_h \rangle_{L^2_h (U_h^\circ)} + \frac{1}{\# U_h^\circ \cdot |p|_2 h} \sum_{\spt \in U_h^\circ \Delta (U_h^\circ - ph)} \sigma(\spt) v_h(\spt + ph) w_h(\spt),
        \end{align*}
        where $\sigma(\spt) = 1$ if $\spt \in U_h^\circ \bs (U_h^\circ - ph)$ and $\sigma(\spt) = -1$ if $\spt \in (U_h^\circ - ph) \bs U_h^\circ$. Note that $U_h^\circ \Delta (U_h^\circ - ph) \subseteq \partial U_h \cup \partial U_h^\circ$, so \eqref{eq: summation by parts} holds with $M_0 = \frac{\# \partial U_h + \# \partial U_h^\circ}{\# U_h^\circ \cdot |p|_2 h}$. When $U_h = \Omega_h$, we have $\# \partial \Omega_h$, $\# \partial \Omega_h^\circ = \bo_{\Omega} (h)$ and $\# \Omega_h^\circ = \bo_{\Omega} (h^2)$, which imply $M_0 = \bo_{\Omega}(1)$.
    \end{proof}

    \begin{lemma}
        \label{lem: nabla_h L^2}
        Let $u$ be the exact solution to the model problem \eqref{eq:PDE}, and let $u_h$ be the numerically approximated solution by solving the linear system in \eqref{eq: overall scheme}. Then
        \begin{equation}
            \label{eq: nabla_h L^2}
            \|\nabla_h (u - u_h)\|_{L^2_h (\Omega_h^\circ)} = \bo_{\tilde{a}, \gx, \gy, u} (h^{11/2})
            \quad \mbox{and} \quad
            \|\phi_h \nabla_h (u - u_h)\|_{L^2_h (\Omega_h^\circ)} = \bo_{\tilde{a}, \gx, \gy, u} (h^6 (\log h)^{1/2}),
        \end{equation}
        where $\phi_h(\spt) = (\dist(\spt, \partial \Omega) + h)^{1/2}$.
    \end{lemma}

    \begin{proof}
        \textbf{Step 1:} In this proof we use a generic constant $C$ to bound any quantity of order $\bo_{\tilde{a}, \gx, \gy, u} (1)$. Take $v_h = h^{-6} (u - u_h)$. Then \Cref{thm: convergence} implies $\|v_h\|_{L^\infty_h (\Omega_h)} \leq C$. Moreover, according to \cref{eq: consistency}, we have $\| \mathcal{L}_h v_h \|_{L^\infty_h (\Omega_h)} \leq C$. By definition \eqref{eq: overall scheme} of $\mathcal{L}_h$ and \Cref{prop: interior c}, we have
        \begin{equation*}
            \mathcal{L}_h v_h (\spt)
            = h^{-2} \sum_{p \in \SS} c_{p, 0} v_h(\spt + ph)
            + h^{-1} \sum_{p \in \SS} (c_{p, 1} + q c_{p, 0}) v_h(\spt + ph) + \bo_{\tilde{a}, \gx, \gy, u} (1)
        \end{equation*}
        for $\spt \in \Omega_h^\circ$ and some $q = q(\spt) = \bo_{\tilde{a}} (1)$, where the coefficients $c_{p, 0}$ and $c_{p, 1}$ are those given in Appendix \ref{app: Interior stencil coefficients}. Define the operators
        \begin{equation*}
            \mathcal{L}_{h, 0} v_h (\spt)
            = h^{-2} \sum_{p \in \SS} c_{p, 0} v_h(\spt + ph)
            \quad \mbox{and} \quad
            \mathcal{L}_{h, 1} v_h (\spt)
            = h^{-1} \sum_{p \in \SS} (c_{p, 1} + q c_{p, 0}) v_h(\spt + ph),
        \end{equation*}
        then
        \begin{equation}
            \label{eq: residual}
            \| \mathcal{L}_{h, 0} v_h + \mathcal{L}_{h, 1} v_h \|_{L^\infty_h (\Omega_h^\circ)} \leq C.
        \end{equation}

        \textbf{Step 2:} (Estimate on $\mathcal{L}_{h, 0} v_h$) For $p \in \mathring{\SS}$, denote $\omega_p = 2$ if $p$ has a zero component, and $\omega_p = 1$ otherwise. One can directly verify that $\mathcal{L}_{h, 0} = \sum_{p \in \mathring{\SS}} \omega_p \partial_{-p} \partial_p$. Therefore, for any grid function $\psi_h$, we obtain from \Cref{lem: summation by parts} and the boundedness of $v_h$ that
        \begin{equation*}
            \langle \mathcal{L}_{h, 0} v_h, \psi_h v_h \rangle_{L^2_h (\Omega_h^\circ)}
            = \sum_{p \in \mathring{\SS}} \omega_p \langle \partial_p v_h, \partial_p(\psi_h v_h) \rangle_{L^2_h (\Omega_h^\circ)} + \bo_{\tilde{a}, \gx, \gy, u} (h^{-1}) \|\psi_h v_h\|_{L^\infty_h (\partial \Omega_h \cup \partial \Omega_h^\circ)}.
        \end{equation*}
        Define the translation operator $T_p: v_h(\spt) \mapsto v_h(\spt + ph)$, then $\partial_p (\psi_h v_h) = \psi_h \partial_p v_h + T_p v_h \cdot \partial_p \psi_h$. It follows that
        \begin{align*}
            \langle \partial_p v_h, \partial_p(\psi_h v_h) \rangle_{L^2_h (\Omega_h^\circ)}
            &= \| \psi_h^{1/2} \partial_p v_h \|_{L^2_h (\Omega_h^\circ)}^2 + \langle \partial_p v_h, T_p v_h \cdot \partial_p \psi_h \rangle_{L^2_h (\Omega_h^\circ)} \\
            &\geq \frac{1}{2} \| \psi_h^{1/2} \partial_p v_h \|_{L^2_h (\Omega_h^\circ)}^2 - C \| \psi_h^{-1/2} \partial_p \psi_h \|_{L^2_h (\Omega_h^\circ)}^2.
        \end{align*}
        Combining last two equations, we obtain
        \begin{equation}
            \label{eq: L_h0 estimate general}
            \langle \mathcal{L}_{h, 0} v_h, \psi_h v_h \rangle_{L^2_h (\Omega_h^\circ)}
            \geq \frac{1}{2} \| \psi_h^{1/2} \nabla_h v_h \|_{L^2_h (\Omega_h^\circ)}^2 - C \| \psi_h^{-1/2} \nabla_h \psi_h \|_{L^2_h (\Omega_h^\circ)}^2 - Ch^{-1} \|\psi_h v_h\|_{L^\infty_h (\partial \Omega_h \cup \partial \Omega_h^\circ)}.
        \end{equation}
        Taking $\psi_h \equiv 1$, we immediately obtain
        \begin{equation}
            \label{eq: L_h0 estimate 1}
            \langle \mathcal{L}_{h, 0} v_h, v_h \rangle_{L^2_h (\Omega_h^\circ)}
            \geq \frac{1}{2} \| \nabla_h v_h \|_{L^2_h (\Omega_h^\circ)}^2 - Ch^{-1}.
        \end{equation}
        Now we take $\psi_h = \phi_h^2$, and it is clear that $\|\psi_h v_h\|_{L^\infty_h (\partial \Omega_h \cup \partial \Omega_h^\circ)} \leq Ch$. Note that $\dist(\cdot, \partial \Omega)$ is $1$-Lipschitz continuous. Together with the mean value theorem, we can obtain
        \begin{equation*}
            \psi_h^{-1} (\spt) \cdot |\nabla_h \psi_h (\spt)|_\infty^2
            \leq C \dist(\spt, \partial \Omega)^{-1}.
        \end{equation*}
        For $n \in \Z$, the number of points in $\Omega_h^\circ$ with $2^n h \leq \dist(\spt, \partial \Omega) < 2^{n + 1} h$ is bounded by $C 2^{-n} h^{-1}$, and the number is $0$ if $n < 0$ or $n > C \log h$. Therefore,
        \begin{equation*}
            \| \psi_h^{-1/2} \nabla_h \psi_h \|_{L^2_h (\Omega_h^\circ)}^2
            \leq Ch^2 \cdot \sum_{n = 0}^{C \log h} C h^{-2} \leq C \log h.
        \end{equation*}
        Substituting into \cref{eq: L_h0 estimate general}, we finally get
        \begin{equation}
            \label{eq: L_h0 estimate 2}
            \langle \mathcal{L}_{h, 0} v_h, \phi_h v_h \rangle_{L^2_h (\Omega_h^\circ)}
            \geq \frac{1}{2} \| \phi_h \nabla_h v_h \|_{L^2_h (\Omega_h^\circ)}^2 - C\log h.
        \end{equation}

        \textbf{Step 3:} (Estimate on $\mathcal{L}_{h, 1} v_h$) For $p \in \mathring{\SS}$, denote $\omega'_p = 4$ if $p$ has a zero component, and $\omega'_p = \sqrt{2}$ otherwise. Moreover, we take the quantity $\omega_p$ from Step 2 and denote $\partial_p^*$ to be the directional derivative of a smooth function in the direction $p / |p|_2$. A direct calculation yields $\mathcal{L}_{h, 1} = \sum_{p \in \mathring{\SS}} (\omega_p \partial_p^* \tilde{a} - q \omega'_p) \partial_p$. Now, for any grid function $\psi_h$, we use the boundedness of $\nabla \tilde{a}$ and Young's inequality to obtain
        \begin{equation*}
            \langle \mathcal{L}_{h, 1} v_h, \psi_h v_h \rangle_{L^2_h (\Omega_h^\circ)}
            \leq C \langle \nabla_h v_h, \psi_h v_h \rangle_{L^2_h (\Omega_h^\circ)}
            \leq \frac{1}{4} \| \psi_h^{1/2} \nabla_h v_h \|_{L^2_h (\Omega_h^\circ)}^2 + C \| \psi_h^{1/2} v_h \|_{L^2_h (\Omega_h^\circ)}^2.
        \end{equation*}
        Either $\psi_h \equiv 1$ or $\psi_h = \phi_h^2$ yields
        \begin{equation}
            \label{eq: L_h1 estimate}
            \langle \mathcal{L}_{h, 1} v_h, \psi_h v_h \rangle_{L^2_h (\Omega_h^\circ)}
            \leq \frac{1}{4} \| \psi_h^{1/2} \nabla_h v_h \|_{L^2_h (\Omega_h^\circ)}^2 + C.
        \end{equation}
        Combining \cref{eq: residual,eq: L_h0 estimate 1,eq: L_h0 estimate 2,eq: L_h1 estimate}, we obtain
        \begin{equation*}
            \|\nabla_h v_h\|_{L^2_h (\Omega_h^\circ)} \leq C h^{-1/2}
            \quad \mbox{and} \quad
            \|\phi_h \nabla_h v_h\|_{L^2_h (\Omega_h^\circ)} \leq C (\log h)^{1/2}.
        \end{equation*}
        This implies \eqref{eq: nabla_h L^2}.
    \end{proof}

    \begin{cor}
        \label{cor: nabla_h L^q}
        Let $u$ be the exact solution to the model problem \eqref{eq:PDE}, and let $u_h$ be the numerically approximated solution by solving the linear system in \eqref{eq: overall scheme}. Then
        \begin{equation}
            \label{eq: nabla_h L^q}
            \|\nabla_h (u - u_h)\|_{L^q_h (\Omega_h^\circ)} = \bo_{\tilde{a}, \gx, \gy, u} \big( h^{5 + 1/q} (\log h)^{ \max\{2/q - 1, 0\} } \big), \quad \forall \, 1 \leq q \leq \infty.
        \end{equation}
    \end{cor}

    \begin{proof}
        Take the function $\phi_h$ in \Cref{lem: nabla_h L^2}. Using the same proof as \Cref{lem: nabla_h L^2}, we can show $\| \phi_h^{-1} \|_{L^2_h (\Omega_h^\circ)} = \bo_{\tilde{a}, \gx, \gy, u} ((\log h)^{1/2})$. Thus,
        \begin{equation}
            \label{eq: nabla_h L^1}
            \|\nabla_h (u - u_h)\|_{L^1_h (\Omega_h^\circ)}
            \leq \|\phi_h \nabla_h (u - u_h)\|_{L^2_h (\Omega_h^\circ)} \| \phi_h^{-1} \|_{L^2_h (\Omega_h^\circ)}
            = \bo_{\tilde{a}, \gx, \gy, u} (h^6 \log h).
        \end{equation}
        Combining estimates \eqref{eq: nabla_h L^2}, \eqref{eq: nabla_h L^1} and the direct consequence $\|\nabla_h (u - u_h)\|_{L^\infty_h (\Omega_h^\circ)} = \bo_{\tilde{a}, \gx, \gy, u} (h^5)$ of \eqref{eq: convergence}, we can obtain \eqref{eq: nabla_h L^q} by an interpolation of $L^q$ spaces.
    \end{proof}

    We are now ready to prove the superconvergence stated in \Cref{thm: superconvergence}.

    \begin{proof}[Proof of \Cref{thm: superconvergence}]
        Considering $\|\nabla u - \nabla u_h\|_{L^\infty_h (\widehat{\Omega}_h)} = \bo_{\tilde{a}, \gx, \gy, u} (h^5)$ and $\# \partial \Omega_h + \# \partial \Omega_h^\circ = \bo_{\Omega} (h^{-1})$, it is sufficient to prove \eqref{eq: superconvergence} with $\widehat{\Omega}_h$ replaced by $\widehat{\Omega}_h \cap (\Omega_h^\circ)^\circ$. Moreover, due to symmetry, we only need to prove the convergence for $\partial_x u - \partial_x u_h$.

        At interior grid points, $\partial_x u_h$ is defined in \Cref{sec: nabla u_h} by
        \begin{equation*}
            \partial_x u_h(\bpt)
            = \frac{1}{h} \sum_{p \in \widehat{\SS}} C_p(h) u_h(\spt + ph) - \frac{1}{h} \sum_{p \in \widehat{\SS}} C_p(h) F(p).
        \end{equation*}
        Since \cref{eq: interior nabla u rearranged} holds with $M = 5$, we obtain
        \begin{equation*}
            \partial_x u(\bpt) - \partial_x u_h(\bpt)
            = \frac{1}{h} \sum_{p \in \widehat{\SS}} C_p(h) (u(\spt + ph) - u_h(\spt + ph)) + \bo_{\tilde{a}, u} (h^6).
        \end{equation*}
        Writing $v_h = u - u_h$ and using the explicit value of $C(p)$ in Appendix \ref{app: coefficients partial_x} we can see that
        \begin{align*}
            \frac{1}{h} \sum_{p \in \widehat{\SS}} C_p(h) v_h(\spt + ph)
            &= \sum_{p \in \mathring{\SS}} \omega_p \partial_p v_h(\spt) + \bo_{\tilde{a}}(1) \|v_h\|_{L^\infty (\Omega_h)} \\
            &\quad \, + \frac{1}{60} \partial_{(1, 0)} v_h(\spt + (1, 0)h) - \frac{1}{60} \partial_{(-1, 0)} v_h(\spt + (-1, 0)h),
        \end{align*}
        where $(\omega_p)_{p \in \mathring{\SS}} = \left( -\frac{1}{5\sqrt{2}}, -\frac{17}{60}, -\frac{1}{5\sqrt{2}}, 0, 0, \frac{1}{5\sqrt{2}}, \frac{17}{60}, \frac{1}{5\sqrt{2}} \right)$. It follows from \Cref{thm: convergence} that
        \begin{equation*}
            |\partial_x u(\bpt) - \partial_x u_h(\bpt)|
            \leq \bo(1) \sum_{p \in \{ (0, 0), (-1, 0), (1, 0) \}} |\nabla_h v_h (\spt + ph))|_\infty
            + \bo_{\tilde{a}, \gx, \gy, u} (h^6).
        \end{equation*}
        For $\bpt = \spt \in \widehat{\Omega}_h \cap (\Omega_h^\circ)^\circ$ and $p \in \{ (0, 0), (-1, 0), (1, 0) \}$, we must have $\spt + ph \in \Omega_h^\circ$. Hence,
        \begin{equation*}
            \|\partial_x u - \partial_x u_h\|_{L^q_h (\widehat{\Omega}_h \cap (\Omega_h^\circ)^\circ)} = \bo(1) \|\nabla_h v_h\|_{L^q_h (\Omega_h^\circ)} + \bo_{\tilde{a}, \gx, \gy, u} (h^6), \quad 1 \leq q \leq \infty.
        \end{equation*}
        Now the estimate \eqref{eq: superconvergence} is a consequence of \Cref{cor: nabla_h L^q}.
    \end{proof}

    \section{Numerical Experiments}
    \label{sec: numerical}

    In this section we present several numerical experiments to illustrate the effectiveness of our proposed scheme and discuss some implementation details of our proposed FDM scheme.

    \subsection{Evaluation of derivatives using function values}
    \label{sec: derivatives}

    The scheme proposed in this article requires frequent evaluation of high-order derivatives. In many applications, it is impossible to obtain an expression of a function. Instead, we can only measure them at certain places. Therefore, it is essential to have an accurate estimate of the derivatives only using function values. Note that all required derivatives can be calculated prior to setting up the linear system \eqref{eq: overall scheme} in the FDM.

    One way to evaluate the derivatives is the moving least squares method proposed in \cite{levin1998approximation}. Suppose $p^* \in \R$ or $\R^2$ is the point at which we would like to evaluate the derivatives of a function $F$. Let $p_1, \ldots, p_K$ be a set of points in $\R$ or $\R^2$. We will approximate the derivatives of $F$ using the function values $F(p_k)$, $1 \leq k \leq K$. Define a diagonal matrix
    \[
        D := 2\diag \big( \eta(|p_1 - p^*|), \ldots, \eta(|p_K - p^*|) \big) \in \R^{K \times K} \quad \mbox{with}\quad \eta(r) = e^{r^2 / h^2}.
    \]
    For $M' \in \mathbb{N}$, we denote by $\Pi_{M'}$ the space of polynomials of total order no more than $M'$. Take
    \[
        q_j(p) := (p - p^*)^j, \quad 0 \leq j \leq M'
    \]
    to be a basis of $\Pi_{M'}$ for the 1D case. For the 2D case, we set $p^* = (x^*, y^*)$, $J := \frac{1}{2} (M' + 1)(M' + 2)$ and take polynomials $q_j, j=1,\ldots, J$ to form a basis of $\Pi_{M'}$ as follows:
    \[
        q_j(x, y) := (x - x^*)^m (y - y^*)^n, \quad 0 \leq m + n \leq M'
        \quad\mbox{with}\quad  j = \frac{1}{2} (m + n) (m + n + 1) + n+1.
    \]
    Let $E := (q_j(p_k))_{1 \leq k \leq K, \, 1 \leq j \leq J} \in \R^{K \times J}$ be a $K\times J$ matrix. Then, according to \cite{levin1998approximation}, the $\omega$-th derivative ($\omega \in \mathbb{N}$ or $\mathbb{N}^2$) is approximated via the formula
    \begin{equation}
        \label{eq: moving least squares}
        F^{(\omega)} (p^*) \approx (F(p_1), \ldots, F(p_K)) D^{-1} E (E^T D^{-1} E)^{-1} (q_1^{(\omega)} (p^*), \ldots, q_J^{(\omega)} (p^*))^T.
    \end{equation}

    Numerical differentiation is prone to round-off errors, and this is worsened by taking the inverse of the matrix $E^T D^{-1} E$ in \eqref{eq: moving least squares}. To mitigate this problem, we try to combine symbolic and numerical calculation in this process. We make a few simplifications as follows. First, we fix some integer $L \in \mathbb{N}$. Then we set the points $\{p_k: 1 \leq k \leq K\}$ by $\{p^* + \ell h/L: -L \leq \ell \leq L\}$ for the 1D case and $\{p^* + (\ell_1 h/L, \ell_2 h/L): -L \leq \ell_1, \ell_2 \leq L\}$ for the 2D case. Now, each component of the matrix $E$ can be expressed as a monomial of a single variable $h$, and the vector $(q_1^{(\omega)} (p^*), \ldots, q_J^{(\omega)} (p^*))$ is a constant vector with only one nonzero element. Furthermore, we take the function $\eta \equiv \frac{1}{2}$. In this case $D$ becomes the identity matrix. Finally, the term
    \begin{equation*}
        D^{-1} E (E^T D^{-1} E)^{-1} (q_1^{(\omega)} (p^*), \ldots, q_J^{(\omega)} (p^*))^T
    \end{equation*}
    can be symbolically calculated in advance. The evaluation of derivatives \eqref{eq: moving least squares} simply becomes a direct linear combination of $F(p_k)$, $1 \leq k \leq K$.

    For $\omega \in \mathbb{N}$ or $\mathbb{N}^2$, let $|\omega|$ be the sum of all components of $\omega$. Since we use polynomials of degree up to $M'$ in the moving least squares algorithm, we expect that the approximation of $F^{(\omega)} (p^*)$ has an accuracy order of $\bo(h^{M' + 1 - |\omega|})$. To correspond with \cref{u:taylor,u:taylor:bdr}, we set $M' = 7$ if the derivative is evaluated at an interior grid point, and $M' = 5$ if the derivative is evaluated at a point on $\partial \Omega$. In addition, in all numerical examples, we take $L = 8$ for differentiating 1D functions and $L = 4$ for differentiating 2D functions.

    \subsection{Examples}
    In this section, we present several numerical examples in different perspectives. \Cref{ex: ex1,ex: ex2} deals with prescribed exact solution, while in \Cref{ex: ex3,ex: ex4} the exact solutions are unknown. Domains with complicated geometries are involved In \Cref{ex: ex2,ex: ex4}. In addition, \Cref{ex: ex1} validates our FDM by a comparison with existing methods in the literature. For the rest of the examples, we try to represent diverse scenarios by casually selecting the functions in the PDE with certain oscillation. In the examples we also compare the results of the proposed sixth-order FDM with a second and a fourth-order method. We use the same strategy for constructing the stencil coefficients in lower-order methods, which are summarized in Appendices \ref{app: Interior stencil coefficients} and \ref{app: order 2 4}. 

    Let $u$ be the exact solution to the model problem \eqref{eq:PDE}. For the accuracy orders $M = 2, 4$ or $6$, we let $u_h^{[M]}$ be the numerical solution computed from our proposed $M$-th order schemes. We measure the relative numerical errors in the $q$-norm (i.e., $L^q_h$ norm in \eqref{eq: norm}) by
    \begin{equation}
        \label{eq: errors}
        e_{q, h}^{[M]}:= \|u - u_h^{[M]}\|_{L^q_h (\Omega_h)} / \|u\|_{L^q (\Omega)}, \quad
        e_{\nabla, q, h}^{[M]}:= \| \nabla u - \nabla u_h^{[M]} \|_{L^q_h (\widehat{\Omega}_h)} / \| \nabla u \|_{L^q (\widehat{\Omega}_h)}.
    \end{equation}
    If the exact solution $u$ is unknown, we take a sufficiently small mesh size $h_{\text{ref}}$ and take the reference solution $u^{[M]}_{h_{\text{ref}}}$ in place of the exact solution $u$. If $h$ is an integer multiple of $h_{\text{ref}}$, then $\Omega_h \subseteq \Omega_{h_{\text{ref}}}$ and $\widehat{\Omega}_h \cap \Omega_h^\circ \subseteq \widehat{\Omega}_{h_{\text{ref}}} \cap \Omega_{h_{\text{ref}}}^\circ$, so we can use \eqref{eq: errors} with a slight change of the domain for calculating errors.
    We use the following two methods to estimate the convergence order:
    \begin{itemize}
        \item[(a)] We estimate the local convergence order at the grid size $h$ by dividing the errors with grid sizes $2h$ and $h$, and then we take the average with multiple $h$ values.

        \item[(b)] We perform linear regression on the data $(-\log_{10} h, -\log_{10} e_h)$ with multiple $h$ ($e_h$ is one of the errors in \eqref{eq: errors}). The coefficient of the linear part is taken as the convergence order.
    \end{itemize}

    \begin{figure}[b!]
        \centering
        \begin{minipage}[t]{0.36 \linewidth}
            \vspace{0pt} % used to align the top of minipages
            \includegraphics[width = \linewidth]{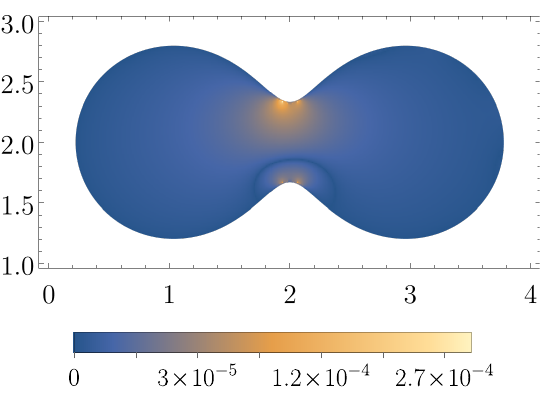}
        \end{minipage}
        \hspace{0.1 \linewidth}
        \begin{minipage}[t]{0.36 \linewidth}
            \vspace{0pt}
            \includegraphics[width = \linewidth]{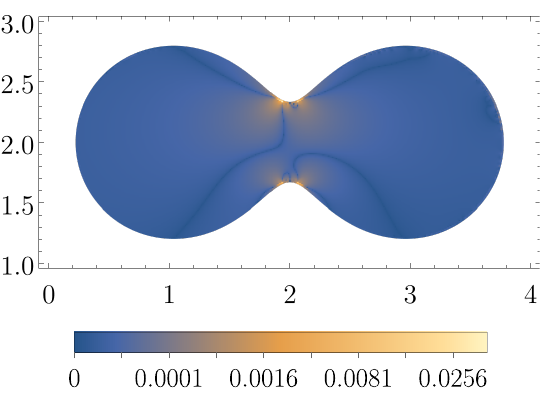}
        \end{minipage}
        \newline
        \begin{minipage}[t]{0.46 \linewidth}
            \vspace{0pt}
            \includegraphics[width = \linewidth]{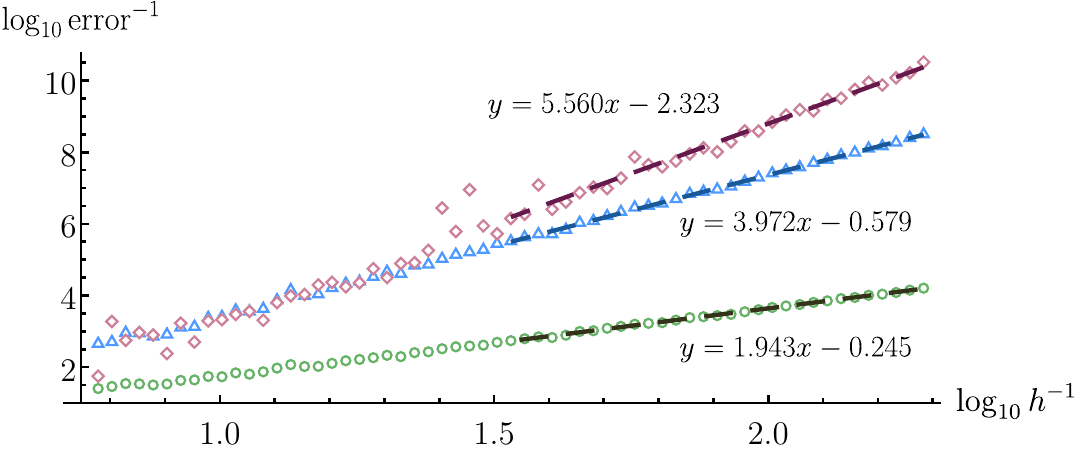}
        \end{minipage}
        \begin{minipage}[t]{0.46 \linewidth}
            \vspace{0pt}
            \includegraphics[width = \linewidth]{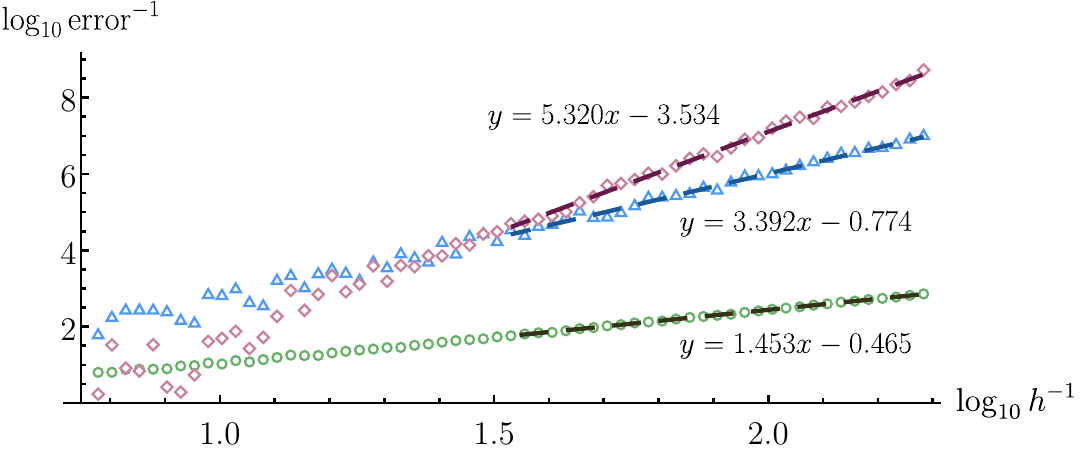}
        \end{minipage}
        \caption{The numerical errors $|u - u_h^{[M]}|$ (first) and $|\partial_x u - \partial_x u_h^{[M]}|$ (second) using $M = 6$ and $h = \frac{1}{96}$ in \Cref{ex: ex1}, and the relative errors $e_{2, h_j}^{[M]}$ (third) and $e_{\nabla, 2, h_j}^{[M]}$ (last) with green, blue and magenta data points representing $M = 2, 4, 6$, respectively. Here $h_j = \frac{1}{6} \times 2^{-j / 12}$, $0 \leq j \leq 60$ and the linear fits are taken from $36 \leq j \leq 60$. Note the nonlinear scaling on first two plots.}
        \label{fig: ex1 errors}
    \end{figure}

    Finally, we discuss some aspects on the implementation of our FDM. First we talk about the choice of the base point $\bpt$ on the boundary. For a boundary grid point $\spt\in \partial \Omega_h$, there must exist $\bpt\in (\spt+[-h,h]^2)\cap \partial \Omega$ such that the line segment from $\spt$ to $\bpt$ is horizontal, vertical or $\pm 45^\circ$. The base point $\bpt$ is taken so that $\|\bpt-\spt\|$ is the smallest among them. Next, as we have mentioned in \Cref{sec: choice of boundary stencil} and \Cref{app: verification of stencil coefs}, when the grid size $h$ is not sufficiently small, we may not be able to construct the stencil according to the 6 stencil types in \Cref{sec: choice of boundary stencil}, or the zeroth order stencil coefficients $\vec{c}_0$ may not be admissible. In the first case, except for $\bpt$, we randomly pick 5 points from a few grid points closest to $\bpt$ to form a stencil. In the second case we normalize $\vec{c}_0$ by $\|\vec{c}_0\|_{\ell_1} = 1$. Moreover, if the augmented matrix $\mathbb{A}^*_0$ in \eqref{extra:eqs} is ill-conditioned or its determinant is below a certain threshold, then we should re-choose the stencil points. These considerations for not sufficiently small $h$ aim to decrease the errors induced by the Taylor expansion and stabilize the numerical results.

    \subsection{Two numerical examples with prescribed solution}

    \begin{example}\label{ex: ex1}
        \rm This example is taken from Section 6.1.2 in \cite{clp21}. Let $\Omega = \{ (x, y) \in \R^2 \setsp ((x - 1)^2 + (y - 2)^2 - 0.75^2) ((x - 3)^2 + (y - 2)^2 - 0.75^2) = 0.3 \}$ and $a(x, y) = 1$. The exact solution $u$ to the model problem \eqref{eq:PDE} is prescribed as $u(x, y) = e^{x + 2y}$. The functions $f$, $g$ in \eqref{eq:PDE} are induced by $u$ through \eqref{eq:PDE}. The numerical results are presented in \Cref{fig: ex1 errors} and \Cref{table: ex1 comparison,table: ex1 convergence u_h,table: ex1 convergence nabla u_h}.
    \end{example}
    \begin{table}[h!]
        \vspace{-6pt}
        \begin{center}
        \begin{NiceTabular}{|c||c|c|c|}[cell-space-limits=2pt]
            \hline
            $\frac{1}{h}$ & $e_{\infty, h}^{[4]}$ & $e_{\infty, h}^{[6]}$ & $e_{\infty, h}^{[4]}$ in \cite{clp21}\\ \hline\hline
            $30$ & 1.512E$-$5 & 9.216E$-$6 & 2.46E$-$5 \\ \hline
            $60$ & 1.083E$-$6 & 5.032E$-$7 & 1.63E$-$6 \\ \hline
            $80$ & 5.436E$-$7 & 6.978E$-$8 & 5.31E$-$7 \\ \hline
        \end{NiceTabular}
        \end{center}
        \vspace{-9pt}
        \caption{\footnotesize A comparison of \Cref{ex: ex1} between our FDM scheme and the numerical results under a fourth-order scheme in \cite{clp21}. Errors $e^{[M]}_{q, h}$ are defined in \eqref{eq: errors}.}
        \label{table: ex1 comparison}
    \end{table}
    \begin{table}[h!]
        \vspace{-6pt}
        \begin{center}
        \begin{NiceTabular}{|c||c|c|c||c|c|c||c|c|c|}[cell-space-limits=2pt]
            \hline
            & \Block{1-3}{$M = 6$} & & & \Block{1-3}{$M = 4$} & & & \Block{1-3}{$M = 2$} & & \\ \hline
            $\frac{1}{h}$ & $e_{\infty, h}^{[6]}$ & ord$_{\infty}$ & ord$_2$ & $e_{\infty, h}^{[4]}$ & ord$_{\infty}$ & ord$_2$ & $e_{\infty, h}^{[2]}$ & ord$_{\infty}$ & ord$_2$ \\ \hline \hline
            $12$       & 3.139E$-$3  &      &      & 6.201E$-$4 &      &      & 3.232E$-$2 &      &      \\ \hline
            $24$       & 5.312E$-$5  & 5.89 & 6.47 & 3.992E$-$5 & 3.96 & 4.12 & 9.584E$-$3 & 1.75 & 1.85 \\ \hline
            $48$       & 9.523E$-$7  & 5.80 & 5.87 & 2.931E$-$6 & 3.77 & 4.03 & 2.755E$-$3 & 1.80 & 1.96 \\ \hline
            $96$       & 4.839E$-$8  & 4.30 & 5.19 & 1.860E$-$7 & 3.98 & 4.02 & 7.180E$-$4 & 1.94 & 1.95 \\ \hline
            $192$      & 6.327E$-$10 & 6.26 & 6.40 & 1.155E$-$8 & 4.01 & 4.01 & 1.879E$-$4 & 1.93 & 2.00 \\ \hline \hline
            Average    &             & 5.56 & 5.98 &            & 3.93 & 4.05 &            & 1.86 & 1.94 \\ \hline
            Linear fit &             & 5.18 & 5.56 &            & 3.96 & 3.97 &            & 1.89 & 1.94 \\ \hline
        \end{NiceTabular}
        \end{center}
        \vspace{-9pt}
        \caption{Convergence order estimates for $u_h$ in \Cref{ex: ex1}. ord$_q$ indicates the estimate of convergence order using $q$-norm. Errors $e^{[M]}_{q, h}$ are defined in \eqref{eq: errors}. The linear fits are performed in the same way as \Cref{fig: ex1 errors}.}
        \label{table: ex1 convergence u_h}
    \end{table}
    \begin{table}[h!]
        \vspace{-6pt}
        \begin{center}
        \begin{NiceTabular}{|c||c|c||c|c||c|c|}[cell-space-limits=2pt]
            \hline
            $\frac{1}{h}$ & $e_{\nabla, \infty, h}^{[6]}$ & order & $e_{\nabla, 2, h}^{[6]}$ & order & $e_{\nabla, 1, h}^{[6]}$ & order \\ \hline \hline
            $12$       & 1.493E$-$1 &      & 1.905E$-$2 &      & 1.792E$-$3  &      \\ \hline
            $24$       & 1.978E$-$3 & 6.24 & 1.390E$-$4 & 7.10 & 1.101E$-$5  & 7.35 \\ \hline
            $48$       & 1.067E$-$4 & 4.21 & 3.949E$-$6 & 5.13 & 1.979E$-$7  & 5.80 \\ \hline
            $96$       & 7.108E$-$6 & 3.91 & 1.109E$-$7 & 5.15 & 4.449E$-$9  & 5.47 \\ \hline
            $192$      & 1.911E$-$7 & 5.22 & 1.878E$-$9 & 5.88 & 5.538E$-$11 & 6.33 \\ \hline \hline
            Average    &            & 4.89 &            & 5.82 &             & 6.24 \\ \hline
            Linear fit &            & 4.54 &            & 5.32 &             & 5.69 \\ \hline
        \end{NiceTabular}
        \end{center}
        \vspace{-9pt}
        \caption{Convergence order estimates for approximated gradient $\nabla u_h$ with $M = 6$ in \Cref{ex: ex1}. Errors $e^{[M]}_{\nabla, q, h}$ are defined in \eqref{eq: errors}. The linear fits are performed in the same way as \Cref{fig: ex1 errors}.}
        \label{table: ex1 convergence nabla u_h}
    \end{table}
    \begin{table}[h!]
        \normalsize
        \begin{center}
        \begin{NiceTabular}{|c||c|c|c||c|c|c|}[cell-space-limits=2pt]
            \hline
            & \Block{1-3}{\Cref{ex: ex1} with $h = \frac{1}{48}$} & & & \Block{1-3}{\Cref{ex: ex3} with $h = \frac{1}{60}$} & & \\ \cline{2-7}
            & $M = 6$ & $M = 4$ & $M = 2$ & $M = 6$ & $M = 4$ & $M = 2$ \\ \hline \hline
            Mesh configuration & 0.457 s & 0.517 s & 0.440 s & 0.018 s & 0.015 s & 0.016 s \\ \hline
            Function evaluation & 0.542 s & 0.117 s & 0.027 s & 0.817 s & 0.192 s & 0.045 s \\ \hline
            Solving linear system & 0.846 s & 0.843 s & 0.811 s & 1.902 s & 1.841 s & 1.825 s \\ \hline
            Estimating numerical gradient & 0.055 s & 0.029 s & 0.009 s & 0.087 s & 0.036 s & 0.014 s \\ \hline\hline
            & \Block{1-3}{\Cref{ex: ex1} with $h = \frac{1}{192}$} & & & \Block{1-3}{\Cref{ex: ex3} with $h = \frac{1}{240}$} & & \\ \hline
            Mesh configuration & 0.983 s & 1.010 s & 1.029 s & 0.195 s & 0.192 s & 0.183 s \\ \hline
            Function evaluation & 9.974 s & 2.285 s & 0.606 s & 16.78 s & 3.419 s & 0.753 s \\ \hline
            Solving linear system & 45.95 s & 46.83 s & 45.42 s & 113.70 s & 111.02 s & 113.61 s \\ \hline
            Estimating numerical gradient & 0.691 s & 0.268 s & 0.092 s & 1.134 s & 0.424 s & 0.150 s \\ \hline
        \end{NiceTabular}
        \end{center}
        \vspace{-9pt}
        \caption{Computation time in \Cref{ex: ex1,ex: ex3} under different $M$. ``Function evaluation'' refers to the evaluation of any symbolic quantities occurred in the FDM, including the estimate of the derivatives. The linear system \eqref{eq: overall scheme} is solved using sparse QR method in cuSPARSE library.}
        \label{table: time}
    \end{table}
    
\begin{figure}[h!]
        \centering
        \begin{minipage}[t]{0.3 \linewidth}
            \vspace{0pt} % used to align the top of minipages
            \includegraphics[width = \linewidth]{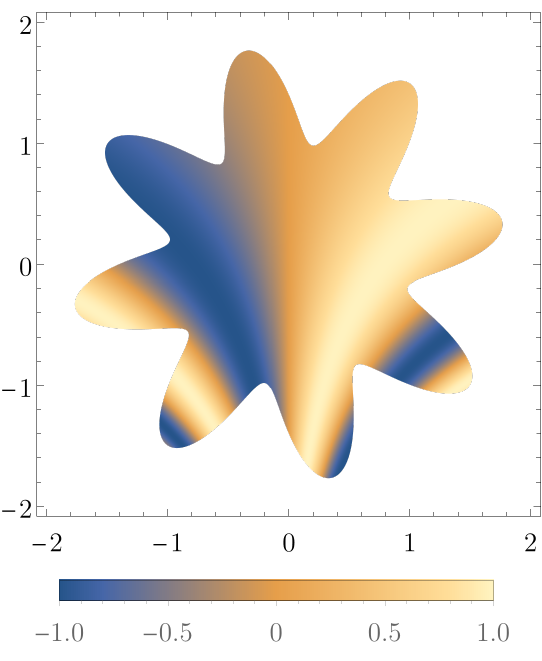}
        \end{minipage}
        \hspace{0.16 \linewidth}
        \begin{minipage}[t]{0.3 \linewidth}
            \vspace{0pt}
            \includegraphics[width = \linewidth]{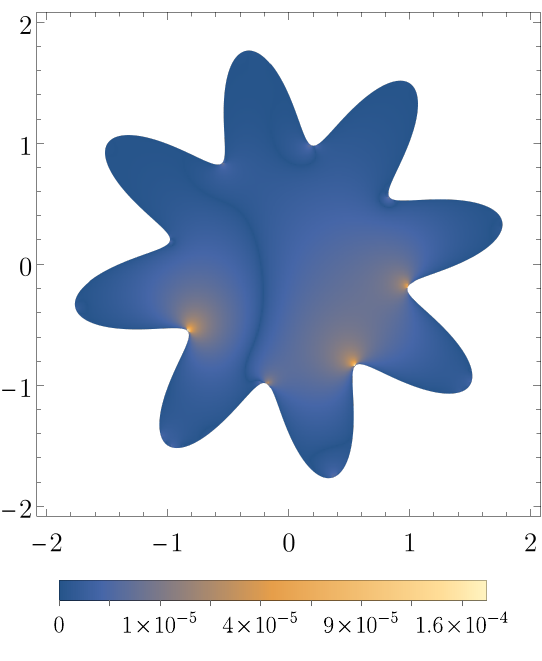}
        \end{minipage}
        \newline
        \begin{minipage}[t]{0.46 \linewidth}
            \vspace{0pt}
            \includegraphics[width = \linewidth]{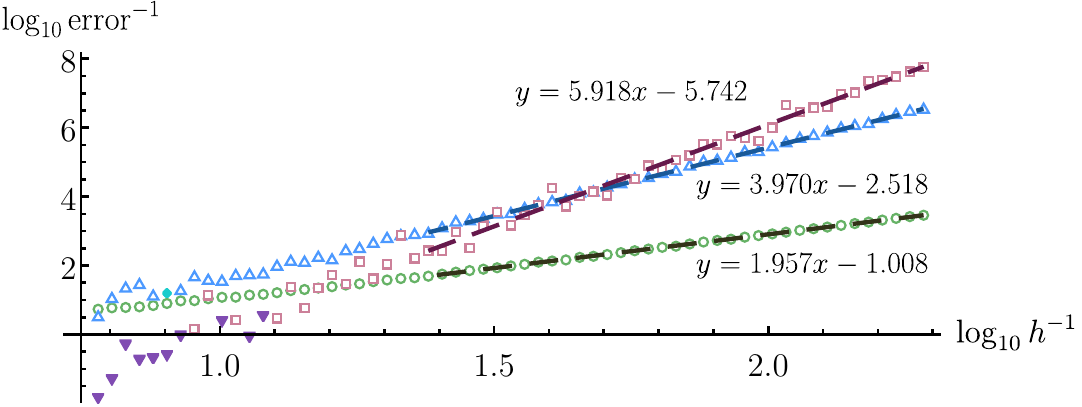}
        \end{minipage}
        \begin{minipage}[t]{0.46 \linewidth}
            \vspace{0pt}
            \includegraphics[width = \linewidth]{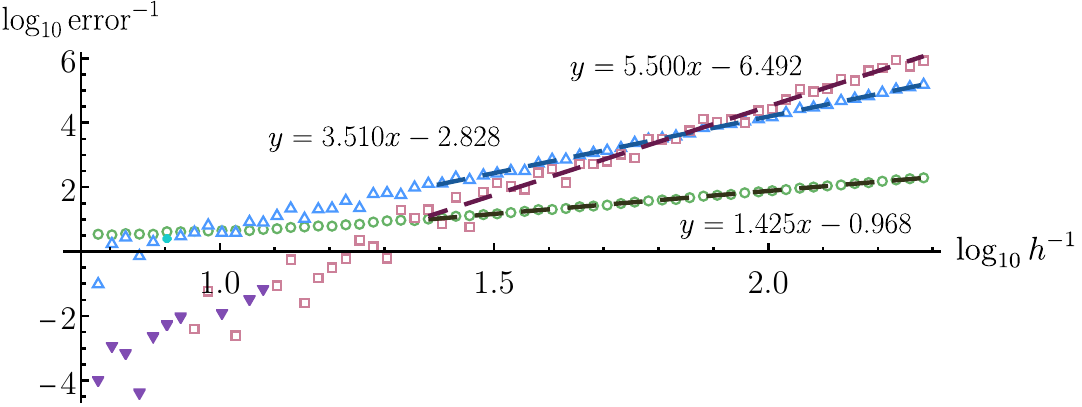}
        \end{minipage}
        \caption{The numerical solution $u_h^{[M]}$ (first) and the error $|u - u_h^{[M]}|$ (second) with $M = 6$ and $h = \frac{1}{96}$ in \Cref{ex: ex2}, and the relative errors $e_{2, h_j}^{[M]}$ (third) and $e_{\nabla, 2, h_j}^{[M]}$ (last). Green, blue and magenta data points represent the errors from $M = 2$, $4$ and $6$ respectively, and solid data points indicate the use of random boundary stencils. Here $h_j = \frac{1}{6} \times 2^{-j / 12}$, $0 \leq j \leq 60$ and the linear fit is taken from $30 \leq j \leq 60$.}
        \label{fig: ex2 errors}
    \end{figure}

    \newpage
    \begin{example}\label{ex: ex2}
    \rm
        Let $\Omega \subset \R^2$ be the region enclosed by the curve $(\gx(t), \gy(t))$ with $\gx(t) = (1.4 + 0.4 \sin(8t)) \cos t$ and $\gy(t) = (1.4 + 0.4 \sin(8t)) \sin t$ for $t\in [0, 2\pi]$. Let $a(x, y) = \arctan \left( \frac{x + 3}{y + 2} \right)$, $u(x, y) = \sin(2x e^{-y})$, and the functions $f$ and $g$  are induced by $u$ through \eqref{eq:PDE}. The results are presented in \Cref{fig: ex2 errors} and \Cref{table: ex2 convergence}.
    \end{example}

    \begin{table}[h!]
        \begin{center}
        \begin{NiceTabular}{|c||c|c||c|c||c|c||c|c|}[cell-space-limits=2pt]
            \hline
            $\frac{1}{h}$ & $e_{\infty, h}^{[6]}$ & order & $e_{2, h}^{[6]}$ & order & $e_{\nabla, \infty, h}^{[6]}$ & order & $e_{\nabla, 1, h}^{[6]}$ & order \\ \hline \hline
            $12$       & 6.428E+0 &      & 3.067E$-$1 &      & 1.582E$+$2 &      & 9.493E$-$1 &      \\ \hline
            $24$       & 1.054E$-$1 & 5.93 & 3.713E$-$3 & 6.37 & 6.918E$-$1 & 7.84 & 3.271E$-$3 & 8.28 \\ \hline
            $48$       & 2.745E$-$3 & 5.26 & 7.121E$-$5 & 5.70 & 5.921E$-$2 & 3.55 & 6.648E$-$5 & 5.62 \\ \hline
            $96$       & 1.782E$-$4 & 3.95 & 2.394E$-$6 & 4.89 & 2.598E$-$3 & 4.51 & 1.604E$-$6 & 5.37 \\ \hline
            $192$      & 1.578E$-$6 & 6.82 & 1.747E$-$8 & 7.10 & 1.704E$-$4 & 3.93 & 1.652E$-$8 & 6.60 \\ \hline \hline
            Average    &            & 5.49 &            & 6.02 &            & 4.96 &            & 6.44 \\ \hline
            Linear fit &            & 5.25 &            & 5.92 &            & 4.57 &            & 5.98 \\ \hline
        \end{NiceTabular}
        \end{center}
        \vspace{-9pt}
        \caption{Convergence order estimates with $M = 6$ in \Cref{ex: ex2}. Errors $e^{[M]}_{q, h}$ and $e^{[M]}_{\nabla, q, h}$ are defined in \eqref{eq: errors}. The linear fits are performed in the same way as \Cref{fig: ex2 errors}.}
        \label{table: ex2 convergence}
    \end{table}

    \newpage
    \subsection{Two numerical examples without explicit solution}

\begin{figure}[h!]
        \centering
        \begin{minipage}[t]{0.34 \linewidth}
            \vspace{0pt} % used to align the top of minipages
            \includegraphics[width = \linewidth]{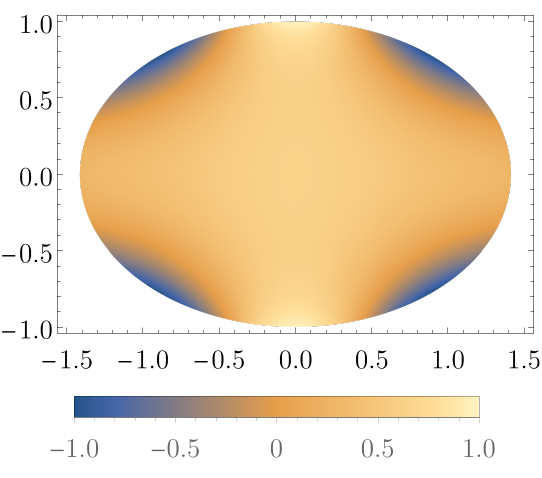}
        \end{minipage}
        \hspace{0.12 \linewidth}
        \begin{minipage}[t]{0.34 \linewidth}
            \vspace{0pt}
            \includegraphics[width = \linewidth]{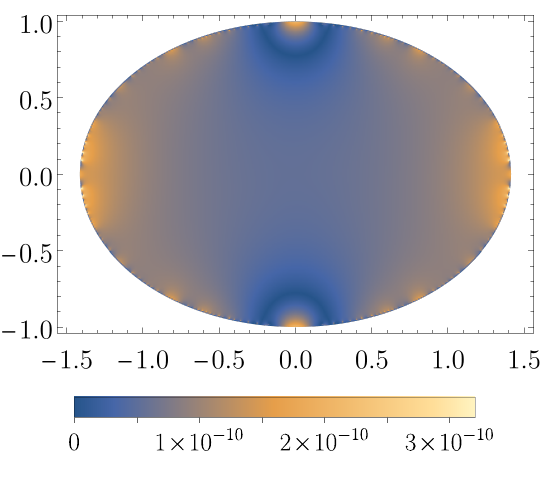}
        \end{minipage}
        \newline
        \begin{minipage}[t]{0.46 \linewidth}
            \vspace{0pt}
            \includegraphics[width = \linewidth]{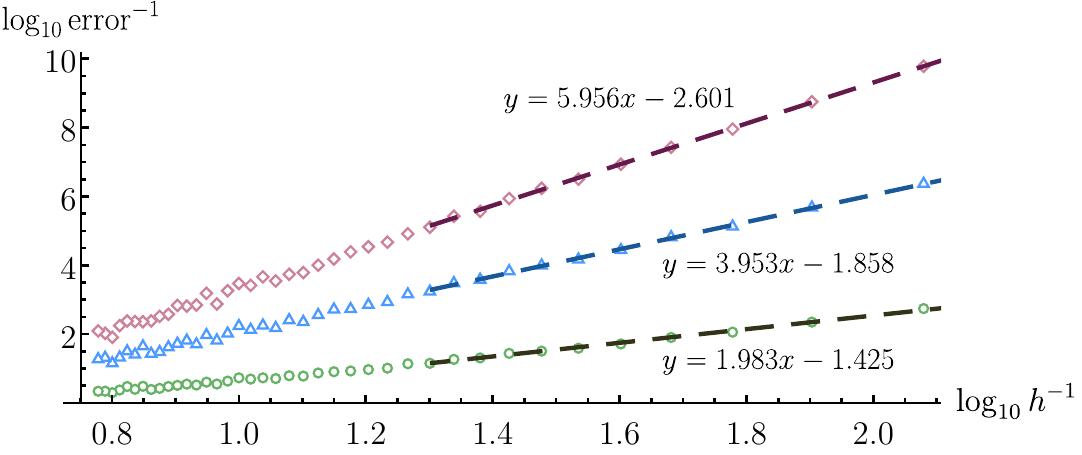}
        \end{minipage}
        \begin{minipage}[t]{0.46 \linewidth}
            \vspace{0pt}
            \includegraphics[width = \linewidth]{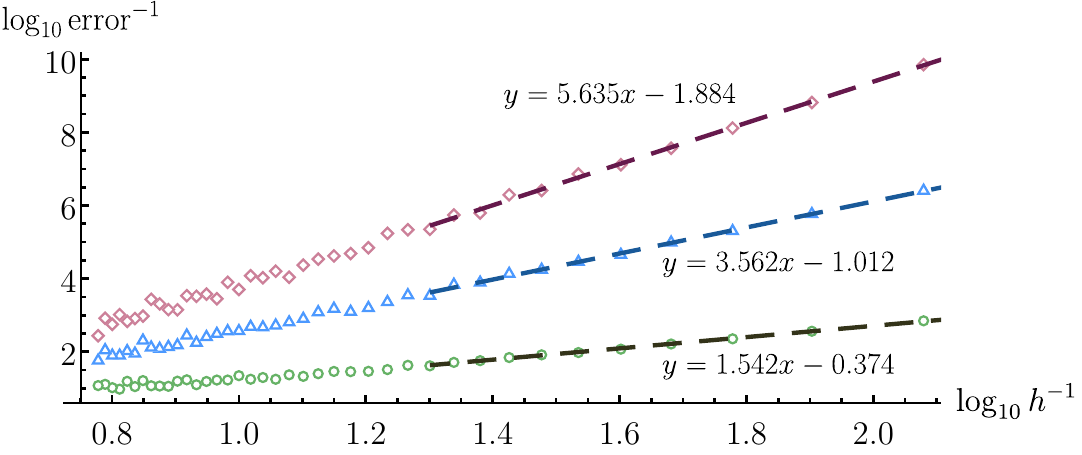}
        \end{minipage}
        \caption{The reference solution $u_{h_{\text{ref}}}^{[M]}$ with $M = 6$ and $h_\text{ref} = \frac{1}{240}$ (first), the error $|u_{h_{\text{ref}}}^{[M]} - u_h^{[M]}|$ with $M = 6$ and $h = \frac{1}{120}$ (second), and the relative errors $e_{2, h_j}^{[M]}$ (third) and $e_{\nabla, 2, h_j}^{[M]}$ (last) in \Cref{ex: ex3}. Here $h_j = j \cdot h_\text{ref}$, $2 \leq j \leq 40$ and the linear fits are taken from $2 \leq j \leq 12$.}
        \label{fig: ex3 errors}
    \end{figure}

    \begin{example}\label{ex: ex3}
        \rm Let $\Omega = \{ (x, y) \in \R^2 \, :\, \frac{1}{2} x^2 + y^2 < 1 \}$. We set
        \begin{equation*}
            a(x, y) = e^{-x^2 - y^2}, \quad
            f(x, y) = 1, \quad \mbox{and} \quad
            g(t) = \cos(5\cos(t)).
        \end{equation*}
        The exact solution $u$ is unknown. We take $h_{\text{ref}} = \frac{1}{240}$ and plot the reference solution $u^{[6]}_{h_{\text{ref}}}$ in \Cref{fig: ex3 errors}. The numerical results are presented in  \Cref{fig: ex3 errors} and \Cref{table: ex3 convergence}.
    \end{example}

    \begin{table}[h!]
        \begin{center}
        \begin{NiceTabular}{|c||c|c|c||c|c|c||c|c|c|}[cell-space-limits=2pt]
            \hline
            & \Block{1-3}{$M = 6$} & & & \Block{1-3}{$M = 4$} & & & \Block{1-3}{$M = 2$} & & \\ \hline
            $\frac{1}{h}$ & $e_{\nabla, \infty, h}^{[6]}$ & ord$_{\infty}$ & ord$_1$ & $e_{\nabla, \infty, h}^{[4]}$ & ord$_{\infty}$ & ord$_1$ & $e_{\nabla, \infty, h}^{[2]}$ & ord$_{\infty}$ & ord$_1$ \\ \hline \hline
            $15$       & 2.246E$-$5  &      &      & 1.104E$-$3 &      &      & 4.717E$-$2 &      &      \\ \hline
            $30$       & 6.388E$-$7  & 5.14 & 5.95 & 1.060E$-$4 & 3.38 & 3.96 & 2.733E$-$2 & 0.79 & 1.72 \\ \hline
            $60$       & 1.758E$-$8  & 5.18 & 5.82 & 1.543E$-$5 & 2.78 & 3.80 & 1.616E$-$2 & 0.76 & 1.57 \\ \hline
            $120$      & 9.826E$-$10 & 4.16 & 6.01 & 1.841E$-$6 & 3.07 & 4.05 & 8.951E$-$3 & 0.85 & 2.13 \\ \hline \hline
            Average    &             & 4.83 & 5.93 &            & 3.08 & 3.94 &            & 0.80 & 1.81 \\ \hline
            Linear fit &             & 4.76 & 5.82 &            & 2.92 & 3.81 &            & 0.78 & 1.79 \\ \hline
        \end{NiceTabular}
        \end{center}
        \vspace{-9pt}
        \caption{Convergence order estimates for $u_h$ in \Cref{ex: ex3}. ord$_q$ indicates the estimate of convergence order using $q$-norm. Errors $e^{[M]}_{\nabla, q, h}$ are defined in \eqref{eq: errors}. The linear fits are performed in the same way as \Cref{fig: ex3 errors}.}
        \label{table: ex3 convergence}
    \end{table}

\begin{figure}[h!]
        \centering
        \begin{minipage}[t]{0.34 \linewidth}
            \vspace{0pt} % used to align the top of minipages
            \includegraphics[width = \linewidth]{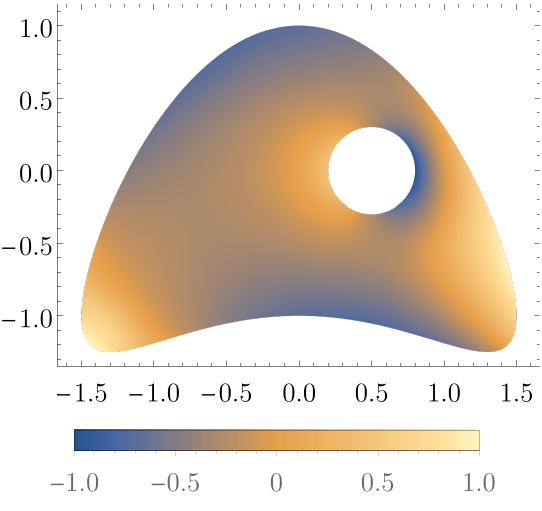}
        \end{minipage}
        \hspace{0.12 \linewidth}
        \begin{minipage}[t]{0.34 \linewidth}
            \vspace{0pt}
            \includegraphics[width = \linewidth]{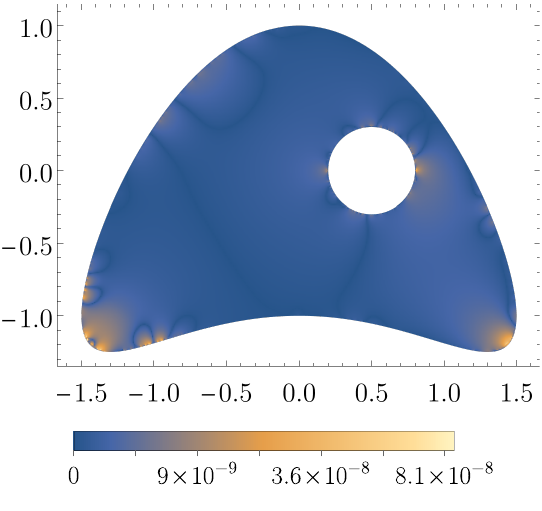}
        \end{minipage}
        \newline
        \begin{minipage}[t]{0.46 \linewidth}
            \vspace{0pt}
            \includegraphics[width = \linewidth]{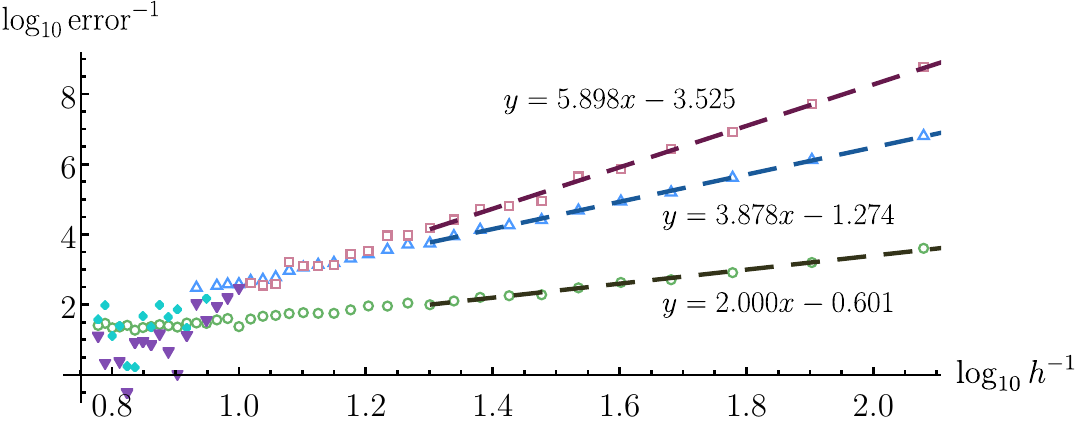}
        \end{minipage}
        \begin{minipage}[t]{0.46 \linewidth}
            \vspace{0pt}
            \includegraphics[width = \linewidth]{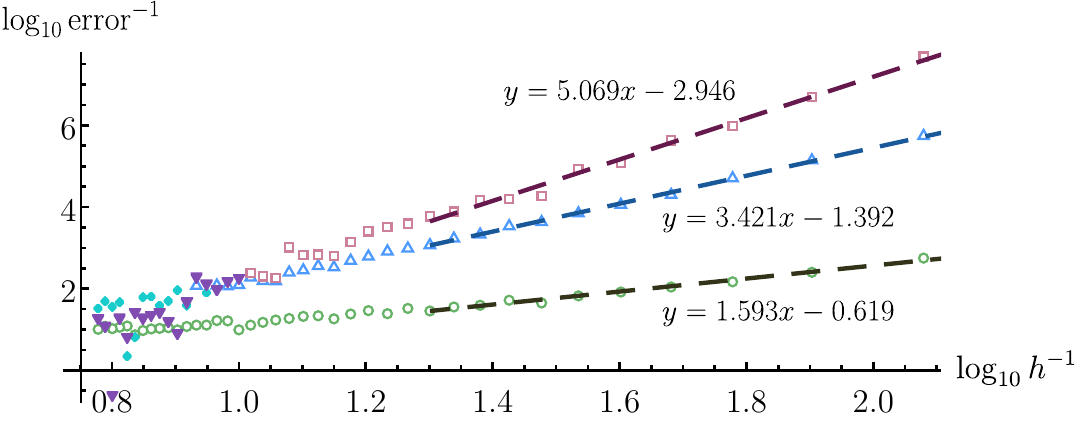}
        \end{minipage}
        \vspace{-9pt}
        \caption{The reference solution $u_{h_{\text{ref}}}^{[M]}$ with $M = 6$ and $h_\text{ref} = \frac{1}{240}$ (first), the error $|u_{h_{\text{ref}}}^{[M]} - u_h^{[M]}|$ with $M = 6$ and $h = \frac{1}{120}$ (second), and the relative errors $e_{2, h_j}^{[M]}$ (third) and $e_{\nabla, 2, h_j}^{[M]}$ (last) in \Cref{ex: ex4}. The grid size $h_j$ and the linear fitting procedure are the same as \Cref{fig: ex3 errors}. Solid data points indicate the use of random boundary stencils.}
        \label{fig: ex4 errors}
    \end{figure}

    \newpage
    \begin{example}\label{ex: ex4}
        \rm
        Let $\Omega\subset \R^2$ be the region between the curves $(\gx^{in}, \gy^{in})$ and
        $(\gx^{out}, \gy^{out})$, where
        \begin{alignat*}{2}
            & \gx^{out} (t) = 1.5 \cos t, \quad
            && \gy^{out} (t) = \sin t - \cos^2 t, \\
            & \gx^{in} (t) = 0.3 \cos t + 0.5, \quad
            && \gy^{in} (t) = 0.3 \sin t.
        \end{alignat*}
        and $g \big|_{\partial \Omega^{out}} (t) = e^{\sin(2t + 1)}$, $g \big|_{\partial \Omega^{in}} (t) = 1 - \cos t$, and
        \[
            a(x, y) = \sin (4xy) + 1.5, \quad
            f(x, y) = \sin (0.5 + x + x^2 - 2y^2).
        \]
        The exact solution $u$ is unknown. We take $h_{\text{ref}} = \frac{1}{240}$ and plot the reference solution $u^{[6]}_{h_{\text{ref}}}$ in \Cref{fig: ex4 errors}. The results are presented in \Cref{fig: ex4 errors} and \Cref{table: ex4 convergence}.
    \end{example}

    \begin{table}[h!]
        \vspace{-3pt}
        \begin{center}
        \begin{NiceTabular}{|c||c|c||c|c||c|c||c|c|}[cell-space-limits=2pt]
            \hline
            $\frac{1}{h}$ & $e_{\infty, h}^{[6]}$ & order & $e_{2, h}^{[6]}$ & order & $e_{\nabla, \infty, h}^{[6]}$ & order & $e_{\nabla, 1, h}^{[6]}$ & order \\ \hline \hline
            $15$       & 2.714E$-$3 &      & 3.558E$-$4 &      & 1.435E$-$3 &      & 5.035E$-$4 &      \\ \hline
            $30$       & 1.572E$-$4 & 4.11 & 1.107E$-$5 & 5.01 & 3.736E$-$4 & 1.94 & 1.630E$-$5 & 4.95 \\ \hline
            $60$       & 3.009E$-$6 & 5.71 & 1.201E$-$7 & 6.53 & 1.312E$-$5 & 4.83 & 2.792E$-$7 & 5.87 \\ \hline
            $120$      & 3.177E$-$8 & 6.57 & 1.691E$-$9 & 6.15 & 2.666E$-$7 & 5.62 & 5.284E$-$9 & 5.72 \\ \hline \hline
            Average    &            & 5.46 &            & 5.89 &            & 4.13 &            & 5.51 \\ \hline
            Linear fit &            & 5.31 &            & 5.90 &            & 4.16 &            & 5.47 \\ \hline
        \end{NiceTabular}
        \end{center}
        \vspace{-9pt}
        \caption{Convergence order estimates with $M = 6$ in \Cref{ex: ex4}. Errors $e^{[M]}_{q, h}$ and $e^{[M]}_{\nabla, q, h}$ are defined in \eqref{eq: errors}. The linear fits are performed in the same way as \Cref{fig: ex3 errors}.}
        \label{table: ex4 convergence}
    \end{table}

    \begin{remark}\rm
        As we can see from \Cref{sec: boundary}, the position of the tangent line plays a fundamental role in the construction of the scheme. We can expect that the numerical solution will deviate from the exact solution if the tangent line does not align well with the boundary $\partial \Omega$. This happens when the grid size is not small enough, or the boundary has a large curvature at some point, which can be seen from the examples above.
    \end{remark}

    \section{Conclusion and Discussion}
    \label{sec: conclusion}

    In this article, we proposed a compact $9$-point finite difference method and proved its sixth-order convergence using the discrete maximum principle. Additionally, we derive a gradient approximation $\nabla u$ directly from  $u_h$ without solving auxiliary equations such that it achieves a superconvergence of $\bo(h^{5 + 1/q} (\log h)^{\max \{2/q - 1, 0\} })$ under the $q$-norm. The proposed scheme is also efficient in that each stencil near the boundary utilizes no more than $8$ points and generally has only $6$ stencil configurations. The stencil coefficients of the scheme can be efficiently obtained either by the analytic expression given in the Appendices or by solving some small linear systems. Moreover, all the derivatives involved can be suitably approximated using function values only. The effectiveness of the method is confirmed by various numerical examples.

    Our method can be easily generalized to the convection-diffusion equation, that is,
    \[
        \begin{cases}
            -\nabla \cdot (a \nabla u) + b \cdot \nabla u = f & \mbox{in} \ \Omega, \\
            u = g & \mbox{on} \ \partial \Omega.
        \end{cases}
    \]
    As an analog of \cref{pde:mod}, the above equation is equivalent to
    \[
        \Delta u = \left( \nabla \tilde{a} + \frac{b}{a} \right) \cdot \nabla u + \tilde{f}\qquad \mbox{with}\quad \tilde{a}:=-\ln a,\quad \tilde{f}:=-\frac{f}{a},
    \]
    which has no essential difference from the pure diffusion case.
    We believe that the same strategy can also be applied to the equation $-\nabla \cdot (a \nabla u) + b \cdot \nabla u + cu = f$ with $c \geq 0$.
    Moreover, instead of the Dirichlet boundary condition $u=g$ on $\partial \Omega$, the techniques developed in \Cref{sec: boundary} can be extended to the Robin (or Neumann) boundary condition $\frac{\partial u}{\partial \nv}+\alpha u=g$ on $\partial \Omega$, where $\nv$ is the outward unit normal vector and $a, g$ are smooth functions on $\partial \Omega$.
    In this case, the left-hand side of \cref{eq: boundary constraint of c} becomes $\sum_{p \in \SS_{\spt}} \frac{2}{m!} \re \left( (p_r + \ii  p_i) e^{-\ii  \theta} \right)^m c_{p, j}$ for $j=0,\ldots, M - 1$ and $m=1, \ldots, M - j$. Moreover, the proof of convergence for the case of Robin or Neumann boundary condition needs to be modified accordingly. We shall address these issues elsewhere.

    \appendix
    \section{Examples of Explicitly Presented Interior Stencil Coefficients}
    \label{app: Interior stencil coefficients}

    Recall that the reference stencil $\SS=[-1,1]^2\cap \mathbb{Z}^2$ is ordered in \eqref{S:order}. For $M=6$, we now present one possible particular real-valued solution to $\mathbb{A}_j \vec{c}_j = \vec{b}_j$ for $j=0,\ldots, 7$ satisfying items (i) and (iii) of \Cref{prop: interior c}, whose general nontrivial solutions have 24 free parameters. For simplicity of presentation, we shall use the notation $\tilde{a}^{(m,n)} := \partial^{(m,n)} \tilde{a}(\bpt)$. Moreover, we introduce an operator $\star: \ \tilde{a}^{(m,n)} \mapsto \tilde{a}^{(n,m)}$ which preserves addition, multiplication and scalar multiplication.
       {
        \allowdisplaybreaks
        \footnotesize
        \begin{flalign*}
            & \vec{c}_0 =[-1, -4, -1, -4, 20, -4, -1, -4, -1]; & \\
            & \vec{c}_1 = \left[-\tfrac{1}{2} ( \tilde{a}^{(0, 1)} + \tilde{a}^{(1, 0)}), \, -2  \tilde{a}^{(1, 0)}, \, \tfrac{1}{2} (\tilde{a}^{(0, 1)} - \tilde{a}^{(1, 0)}), \, -2  \tilde{a}^{(0, 1)}, \, 0,
            2 \tilde{a}^{(0, 1)}, \, \tfrac{1}{2} (- \tilde{a}^{(0, 1)} +  \tilde{a}^{(1, 0)}), \, 2 \tilde{a}^{(1, 0)}, \, \tfrac{1}{2} (\tilde{a}^{(0, 1)} + \tilde{a}^{(1, 0)}) \right]; \\
            & \vec{c}_2 =\left[d_2 + d_3, \, d_1 - d_3, \, -d_2 + d_3, \, -d_1 - d_3, \, 0, \, -d_1 - d_3, \, -d_2 + d_3, \, d_1 - d_3, \, d_2 + d_3\right], \ \mbox{where} \\
            & \qquad d_1 := \frac{1}{2} ( \tilde{a}^{(2, 0)} - \tilde{a}^{(0, 2)}) + \frac{3}{4} ( [\tilde{a}^{(1, 0)}]^2 - [\tilde{a}^{(0, 1)}]^2), \
            d_2 := \frac{1}{2}  \tilde{a}^{(1, 1)} + \frac{3}{4}  \tilde{a}^{(1, 0)}  \tilde{a}^{(0, 1)}, \\
            & \qquad d_3 := \frac{2}{25} ( \tilde{a}^{(2, 0)} + \tilde{a}^{(0, 2)}) + \frac{3}{25} ( [\tilde{a}^{(1, 0)}]^2 + [\tilde{a}^{(0, 1)} ]^2); \\
            & \vec{c}_3 =\left[ d_4 + d_4^\star, \, d_5, \, d_4 - d_4^\star, \, d_5^\star, \, 0, \, -d_5^\star, \, -d_4 + d_4^\star, \, -d_5, \, -d_4 - d_4^\star \right], \ \mbox{where} \\
            & \qquad d_4 := \frac{1}{600} \big( 13 [\tilde{a}^{(1, 0)}]^3 - 12 \tilde{a}^{(1, 0)} [\tilde{a}^{(0, 1)}]^2 + 24 \tilde{a}^{(1, 0)} \tilde{a}^{(0, 2)} - 11 \tilde{a}^{(2, 0)} \tilde{a}^{(1, 0)} + 115 \tilde{a}^{(1, 1)} \tilde{a}^{(0, 1)} - 115 \tilde{a}^{(1, 2)} - 15 \tilde{a}^{(3, 0)} \big), \\
            & \qquad d_5 := \frac{1}{600} \big( 37 [\tilde{a}^{(1, 0)}]^3 + 87 \tilde{a}^{(1, 0)} [\tilde{a}^{(0, 1)}]^2 - 174 \tilde{a}^{(1, 0)} \tilde{a}^{(0, 2)} + 46 \tilde{a}^{(2, 0)} \tilde{a}^{(1, 0)} - 80 \tilde{a}^{(1, 1)} \tilde{a}^{(0, 1)} + 80 \tilde{a}^{(1, 2)} - 120 \tilde{a}^{(3, 0)} \big); \\
            & \vec{c}_4 =\left[ d_6 + d_6^\star, \, d_7 - d_7^\star, \, -d_6 - d_6^\star, \, -d_7 + d_7^\star, \, 0, \, -d_7 + d_7^\star, \, -d_6 - d_6^\star, \, d_7 - d_7^\star, \, d_6 + d_6^\star \right], \ \mbox{where} \\
            & \qquad d_6 := \frac{1}{1200} \big( 13 [\tilde{a}^{(1, 0)}]^3 \tilde{a}^{(0, 1)} - 15 \tilde{a}^{(3, 0)} \tilde{a}^{(0, 1)} - 11 \tilde{a}^{(2, 0)} \tilde{a}^{(1, 0)} \tilde{a}^{(0, 1)} - 11 [\tilde{a}^{(1, 0)}]^2 \tilde{a}^{(1, 1)} - 8 \tilde{a}^{(2, 0)} \tilde{a}^{(1, 1)} \\
            & \qquad\qquad + 60 \tilde{a}^{(3, 1)} - 75 \tilde{a}^{(2, 1)} \tilde{a}^{(1, 0)} \big), \\
            & \qquad d_7 := \frac{1}{2400} \big( 31 [\tilde{a}^{(1, 0)}]^4 - 104 [\tilde{a}^{(1, 0)}]^2 \tilde{a}^{(2, 0)} - 100 \tilde{a}^{(3, 0)} \tilde{a}^{(1, 0)} + 60 \tilde{a}^{(1, 0)} \tilde{a}^{(1, 2)} + 44 [\tilde{a}^{(2, 0)}]^2 + 80 \tilde{a}^{(4, 0)} \big); \\
            & \vec{c}_5 =\left[ 0, \, d_8, \, 0, \, d_8^\star, \, 0, \, -d_8^\star, \, 0, \, -d_8, \, 0 \right], \ \mbox{where} \\
            & \qquad d_8 := \frac{1}{4800} \big( 22 [\tilde{a}^{(0, 1)}]^2 \tilde{a}^{(1, 2)}-31 ([\tilde{a}^{(0, 1)}]^4+[\tilde{a}^{(1, 0)}]^4) \tilde{a}^{(1, 0)} - 40 (\tilde{a}^{(1, 4)}+\tilde{a}^{(0, 3)} \tilde{a}^{(1, 1)}+\tilde{a}^{(5, 0)}) \\
            & \qquad\qquad + 42 ([\tilde{a}^{(0, 1)}]^2 \tilde{a}^{(1, 0)} \tilde{a}^{(2, 0)}+[\tilde{a}^{(1, 0)}]^3 \tilde{a}^{(2, 0)}-[\tilde{a}^{(0, 1)}]^2 [\tilde{a}^{(1, 0)}]^3+[\tilde{a}^{(0, 1)}]^2 \tilde{a}^{(3, 0)}) - 62 [\tilde{a}^{(0, 1)}]^3 \tilde{a}^{(1, 1)} \\
            & \qquad\qquad + 44 (\tilde{a}^{(0, 1)} \tilde{a}^{(1, 1)} \tilde{a}^{(2, 0)}-\tilde{a}^{(1, 2)} \tilde{a}^{(2, 0)} - \tilde{a}^{(0, 2)} \tilde{a}^{(1, 2)} - [\tilde{a}^{(0, 2)}]^2 \tilde{a}^{(1, 0)}) + 80 (\tilde{a}^{(0, 1)} \tilde{a}^{(1, 3)}-\tilde{a}^{(0, 4)} \tilde{a}^{(1, 0)}) \\
            & \qquad\qquad +84 (\tilde{a}^{(0, 2)} [\tilde{a}^{(1, 0)}]^3-\tilde{a}^{(0, 2)} \tilde{a}^{(1, 0)} \tilde{a}^{(2, 0)}-\tilde{a}^{(0, 2)} \tilde{a}^{(3, 0)}-\tilde{a}^{(2, 0)} \tilde{a}^{(3, 0)}) + 100 (\tilde{a}^{(0, 3)}+\tilde{a}^{(2, 1)}) \tilde{a}^{(0, 1)} \tilde{a}^{(1, 0)} \\
            & \qquad\qquad +104 [\tilde{a}^{(0, 1)}]^2 \tilde{a}^{(0, 2)} \tilde{a}^{(1, 0)} +120 \tilde{a}^{(1, 0)} [\tilde{a}^{(2, 0)}]^2 +122 [\tilde{a}^{(1, 0)}]^2 \tilde{a}^{(3, 0)}+160 (-\tilde{a}^{(1, 0)} \tilde{a}^{(2, 2)}+\tilde{a}^{(0, 1)} \tilde{a}^{(3, 1)}-\tilde{a}^{(3, 2)}) \\
            & \qquad\qquad +164 \tilde{a}^{(0, 1)} \tilde{a}^{(0, 2)} \tilde{a}^{(1, 1)} -242 \tilde{a}^{(0, 1)} [\tilde{a}^{(1, 0)}]^2 \tilde{a}^{(1, 1)}+280 \tilde{a}^{(1, 0)} [\tilde{a}^{(1, 1)}]^2+302 [\tilde{a}^{(1, 0)}]^2 \tilde{a}^{(1, 2)}+360 \tilde{a}^{(1, 1)} \tilde{a}^{(2, 1)} \big); \\
            & \vec{c}_6 = \vec{c}_7 = 0.
        \end{flalign*}
    }

    For a fourth-order scheme with $M=4$, a particular solution satisfying items (i) and (iii) of \cref{prop: interior c} with $M = 4$ is given by: $\vec{c}_0, \vec{c}_1$ are the same as the case $M = 6$, $\vec{c}_4 = \vec{c}_5 = 0$ and
    \begin{align*}
        & \vec{c}_2 =\left[ d_1, \, d_2, \, -d_1, \, -d_2, \, 0, \, -d_2, \, -d_1, \, d_2, \, d_1 \right], \ \mbox{where} \\
        & \quad d_1 := \frac{1}{2} \tilde{a}^{(1, 1)} - \frac{1}{4} \tilde{a}^{(1, 0)} \tilde{a}^{(0, 1)},
        d_2 := \frac{1}{2} \tilde{a}^{(2, 0)} - \frac{1}{2} \tilde{a}^{(0, 2)} - \frac{1}{4} [\tilde{a}^{(1, 0)}]^2 + \frac{1}{4} [\tilde{a}^{(0, 1)}]^2; \\
        & \vec{c}_3 =\left[ 0, \, d_3, \, 0, \, d_3^\star, \, 0, \, -d_3^\star, \, 0, \, -d_3, \, 0 \right], \ \mbox{where} \\
        & \quad d_3 := \frac{1}{8} \Big( ([\tilde{a}^{(1, 0)}]^2 + [\tilde{a}^{(0, 1)}]^2 - 2\tilde{a}^{(0, 2)} - 2\tilde{a}^{(1, 1)}) \tilde{a}^{(1, 0)} + ([\tilde{a}^{(1, 0)}]^2 - 2\tilde{a}^{(1, 1)}) \tilde{a}^{(0, 1)} - 2\tilde{a}^{(1, 2)} - 2\tilde{a}^{(3, 0)} \Big).
    \end{align*}

    For a second-order scheme with $M=2$, a particular solution satisfying items (i) and (iii) of \cref{prop: interior c} with $M = 2$ is given by $\vec{c}_2 = \vec{c}_3 = 0$ and
    \begin{align*}
        & \vec{c}_0 = [0, -1, 0, -1, 4, -1, 0, -1, 0]; \\
        & \vec{c}_1 = [0, \, -\tfrac{1}{2}  \tilde{a}^{(1, 0)}, \, 0, \, -\tfrac{1}{2} \tilde{a}^{(0, 1)}, \, 0, \tfrac{1}{2} \tilde{a}^{(0, 1)}, \, 0, \, \tfrac{1}{2} \tilde{a}^{(1, 0)}, \, 0].
    \end{align*}

    \section{Stencil Coefficients for Approximating $\partial_x u$ at Interior Grid Points}
    \label{app: coefficients partial_x}

    Here we present one possible particular real-valued solution to the linear system \eqref{eq: interior nabla constraints} for approximating $\partial_x  u$. We discuss two cases: $\widehat{\SS} = \SS = [-1,1]^2 \cap \Z^2$ with $M = 3$ (fourth-order), and $\widehat{\SS} = \SS \cup \{ (\pm 2, 0), (0, \pm 2) \}$ with $M = 4$ (fifth-order). We use the same convention and notation as in Appendix \ref{app: Interior stencil coefficients}. For $M = 4$, the ordering of the set $\widehat{\SS}$ is given by the ordering of $\SS$ in \eqref{S:order} followed by $(-2, 0)$, $(0, -2)$, $(0, 2)$, $(2, 0)$.

    A fourth-order approximation of $u_x$ from numerical $u_h$ using the original reference stencil $\SS$ is
    \begin{align*}
        & \vec{c}_0 = \left[ -\tfrac{1}{12}, \, -\tfrac{1}{3}, \, -\tfrac{1}{12}, \, 0, \, 0, \, 0, \, \tfrac{1}{12}, \, \tfrac{1}{3}, \, \tfrac{1}{12} \right]; \\
        & \vec{c}_1 = \frac{1}{24} [ 2\tilde{a}^{(1, 0)} - \tilde{a}^{(0, 1)}, \, -4 \tilde{a}^{(1, 0)}, \, 2\tilde{a}^{(1, 0)} + \tilde{a}^{(0, 1)}, \, 0, \, 0, \, 0, \, 2\tilde{a}^{(1, 0)} + \tilde{a}^{(0, 1)}, \, -4 \tilde{a}^{(1, 0)}, \, 2\tilde{a}^{(1, 0)} - \tilde{a}^{(0, 1)}]; \\
        & \vec{c}_2 = [d_1 + d_2, \, 0, \, d_1 - d_2, \, 0, \, 0, \, 0, \, -d_1 + d_2, \, 0, \, -d_1 - d_2], \ \mbox{where} \\
        & \qquad d_1 = \frac{1}{24} ([\tilde{a}^{(1, 0)}]^2 + \tilde{a}^{(2, 0)}), \
        d_2 = \frac{1}{24} (\tilde{a}^{(1, 1)} + \tilde{a}^{(1, 0)} \tilde{a}^{(0, 1)}); \\
        & \vec{c}_3 = \vec{c}_4 = 0.
    \end{align*}

    A fifth-order approximation of $u_x$ using the extended reference stencil $\SS \cup \{ (\pm 2, 0), (0, \pm 2) \}$ is
    {\footnotesize
    \begin{align*}
        & \vec{c}_0 = \left[ -\tfrac{1}{10}, \, -\tfrac{4}{15}, \, -\tfrac{1}{10}, \, 0, \, 0, \, 0, \, \tfrac{1}{10}, \, \tfrac{4}{15}, \, \tfrac{1}{10}, -\tfrac{1}{60}, \, 0, \, 0, \, \tfrac{1}{60} \right]; \\
        & \vec{c}_1 = \frac{1}{40} [ -\tilde{a}^{(1, 0)} - 2\tilde{a}^{(0, 1)}, \, -8 \tilde{a}^{(1, 0)}, \, -\tilde{a}^{(1, 0)} + 2\tilde{a}^{(0, 1)}, \, 0, \, 20\tilde{a}^{(1, 0)}, \, 0, \, \\
        & \qquad -\tilde{a}^{(1, 0)} + 2\tilde{a}^{(0, 1)}, \, -8 \tilde{a}^{(1, 0)}, \, -\tilde{a}^{(1, 0)} - 2\tilde{a}^{(0, 1)}, \, 0, \, 0, \, 0, \, 0]; \\
        & \vec{c}_2 = [d_1 + d_2, \, d_3, \, d_1 - d_2, \, 0, \, 0, \, 0, \, -d_1 + d_2, \, -d_3, \, -d_1 - d_2, \, 0, \, 0, \, 0, \, 0], \ \mbox{where} \\
        & \qquad d_1 = \frac{1}{240} \big( -[\tilde{a}^{(1, 0)}]^2 - 2[\tilde{a}^{(0, 1)}]^2 + 4 \tilde{a}^{(0, 2)} \big), \
        d_2 = \frac{1}{80} \big( 4 \tilde{a}^{(1, 1)} - \tilde{a}^{(1, 0)} \tilde{a}^{(0, 1)} \big), \\
        & \qquad d_3 = \frac{1}{60} \big( [\tilde{a}^{(0, 1)}]^2 - [\tilde{a}^{(1, 0)}]^2 + 6 \tilde{a}^{(2, 0)} - 2 \tilde{a}^{(0, 2)} \big); \\
        & \vec{c}_3 = [d_4 + d_5, \, d_6, \, d_4 - d_5, \, 0, \, 0, \, 0, \, d_4 - d_5, \, d_6, \, d_4 + d_5, \, 0, \, 0, \, 0, \, 0], \ \mbox{where} \\
        & \qquad d_4 = \frac{1}{480} \big( -[\tilde{a}^{(0, 1)}]^2 \tilde{a}^{(1, 0)} + 2 \tilde{a}^{(0, 2)} \tilde{a}^{(1, 0)} + [\tilde{a}^{(1, 0)}]^3 + 4 \tilde{a}^{(0, 1)} \tilde{a}^{(1, 1)} - 4 \tilde{a}^{(1, 2)} + 2 \tilde{a}^{(1, 0)} \tilde{a}^{(2, 0)} + 12 \tilde{a}^{(3, 0)} \big), \\
        & \qquad d_5 = \frac{1}{480} \big( \tilde{a}^{(0, 1)} (2 \tilde{a}^{(0, 2)} - [\tilde{a}^{(1, 0)}]^2) - 2 (\tilde{a}^{(0, 3)} - 2 \tilde{a}^{(1, 0)} \tilde{a}^{(1, 1)} + 5 \tilde{a}^{(2, 1)}) \big), \\
        & \qquad d_6 = \frac{1}{240} \big( [\tilde{a}^{(0, 1)}]^2 \tilde{a}^{(1, 0)} - 2 \tilde{a}^{(0, 2)} \tilde{a}^{(1, 0)} - [\tilde{a}^{(1, 0)}]^3 - 4 \tilde{a}^{(0, 1)} \tilde{a}^{(1, 1)} + 4 \tilde{a}^{(1, 2)} - 2 \tilde{a}^{(1, 0)} \tilde{a}^{(2, 0)} - 12 \tilde{a}^{(3, 0)} \big); \\
        & \vec{c}_4 = [d_7 + d_8, \, 0, \, d_7 - d_8, \, 0, \, 0, \, 0, \, -d_7 + d_8, \, 0, \, -d_7 - d_8, \, 0, \, 0, \, 0, \, 0], \ \mbox{where} \\
        & \qquad d_7 = \frac{1}{960} \big( [\tilde{a}^{(1, 0)}]^4 + 3 [\tilde{a}^{(1, 0)}]^2 \tilde{a}^{(2, 0)} + \tilde{a}^{(1, 0)} (-3 \tilde{a}^{(0, 1)} \tilde{a}^{(1, 1)} + 3 \tilde{a}^{(1, 2)} + 11 \tilde{a}^{(3, 0)}) \\
        & \qquad \qquad + 4 (-[\tilde{a}^{(1, 1)}]^2 + [\tilde{a}^{(2, 0)}]^2 - \tilde{a}^{(0, 1)} \tilde{a}^{(2, 1)} + \tilde{a}^{(2, 2)} + \tilde{a}^{(4, 0)}) \big), \\
        & \qquad d_8 = \frac{1}{960} \big( -\tilde{a}^{(0, 3)} \tilde{a}^{(1, 0)} - 4 [\tilde{a}^{(0, 1)}]^2 \tilde{a}^{(1, 1)} - 4 \tilde{a}^{(0, 2)} \tilde{a}^{(1, 1)} + [\tilde{a}^{(1, 0)}]^2 \tilde{a}^{(1, 1)} + 4 \tilde{a}^{(1, 3)} + 4 \tilde{a}^{(1, 1)} \tilde{a}^{(2, 0)} \\
        & \qquad \qquad - \tilde{a}^{(1, 0)} \tilde{a}^{(2, 1)} + \tilde{a}^{(0, 1)} (\tilde{a}^{(0, 2)} \tilde{a}^{(1, 0)} + [\tilde{a}^{(1, 0)}]^3 + 2 \tilde{a}^{(1, 0)} \tilde{a}^{(2, 0)} + 12 \tilde{a}^{(3, 0)}) + 4 \tilde{a}^{(3, 1)} \big); \\
        & \vec{c}_5 = 0.
    \end{align*}
    }
    The above stencil coefficients together with $\vec{c}_6 = 0$ satisfies the linear system \eqref{eq: interior nabla constraints} with $M = 5$.

    \section{Existence of Admissible Solutions $\vec{c}_0$ to $\mathbb{A}^*_0 \vec{c}_0=\vec{b}^*_0$ Given in \Cref{sec: choice of boundary stencil}}
    \label{app: verification of stencil coefs}

    In this section, we verify our claim in \Cref{sec: choice of boundary stencil} that when $h$ is small enough, we can obtain a unique stable admissible zeroth-order solution $\vec{c}_0$ from $\mathbb{A}^*_0 \vec{c}_0 = \vec{b}_0^*$ satisfying all the conditions in \Cref{def:admissible}. Note that $\mathbb{A}^*_0$ depends on the stencil $\SS_{\spt}$, the base point $\bpt$ and the tangent angle $\theta$. However, when we discussed the construction of the stencil $\SS_{\spt}$, we only considered the position of the directed tangent line $L_{\bpt}$ and did not care about the exact location of the base point $\bpt$ on the line. This is due to the following result, which states that as long as the augmented data $\mathbb{A}^*_{0; k}$ and $\vec{b}_0^* (k)$ for $k=1,\ldots, \# \SS_{\spt} - 5$ in \eqref{extra:eqs} only depends on the position of $L_{\bpt}$, then so does the solution $\vec{c}_0$ to the augmented linear system $\mathbb{A}_0^* \vec{c}_0 = \vec{b}_0^*$. In other words, if there are two identical stencils $\spt_1 + h \SS_{\spt_1}$ and $\spt_2 + h \SS_{\spt_2}$ with (possibly different) base points $\bpt_1$, $\bpt_2$ such that $L_{\bpt_1} = L_{\bpt_2}$, then the corresponding solutions $\vec{c}_0$ must be the same.

    \begin{prop}
        \label{prop: independent of base point}
        Let $M \in \mathbb{N}$, $L$ be a straight line in $\R^2$ with direction angle $\theta \in (-\pi, \pi]$, $\spt \in \R^2$, and $\SS_{\spt}$ be a finite set of $\mathbb{R}^2$ with $\# \SS_{\spt} \geq M + 1$. For any point $\bpt$ on $L$, define $\psv = (\psv_r, \psv_i) = p + (\spt - \bpt) / h$ for $p \in \SS_{\spt}$ and an associated matrix $\mathbb{A}_0$ by
        \begin{equation*}
            \mathbb{A}_0 = \left( 2\im \left( (\psv_r + \ii \psv_i) e^{-\ii  \theta} \right)^m \right)_{1 \leq m \leq M + 1, \ p \in \SS_{\spt}}.
        \end{equation*}
        Let $\vec{b}_0 = (0, \ldots, 0) \in \R^{M + 1}$. Now we augment the linear system $\mathbb{A}_0 \vec{c}_0 = \vec{b}_0$ into a square linear system $\mathbb{A}_0^* \vec{c}_0 = \vec{b}_0^*$ as in \cref{extra:eqs}. If the augmented data $\mathbb{A}^*_{0; k}$ and $\vec{b}_0^* (k)$ in \eqref{extra:eqs}
         for $k=1, \ldots, \# \SS_{\spt} - M - 1$ do not depend on the choice of $\bpt \in L$, then the same is true for the solution $\vec{c}_0$ to  $\mathbb{A}_0^* \vec{c}_0 = \vec{b}_0^*$.
    \end{prop}

    \begin{proof}
        Consider an arbitrary base point $\tbpt\in L$ and define its associated matrix
        \[
        \tilde{\mathbb{A}}_0:=\left( 2\im \left( (\tilde{p}_r + \ii \tilde{p}_i) e^{-\ii  \theta} \right)^m \right)_{1 \leq m \leq M + 1, \ p \in \SS_{\spt}} \quad \mbox{with}\quad \tilde{p}=(\tilde{p}_r, \tilde{p}_i):=p+(\spt-\tbpt)/h.
        \]
        Because both $\bpt$ and $\tbpt$ lie on the line $L$ with the tangent angle $\theta$,
        we must have $\bpt=\tbpt+ e^{i\theta}r$ with $r=|\tbpt-\bpt|$ or $r=-|\tbpt-\bpt|$ depending on whether the vector from $\bpt$ to $\tbpt$ agrees with the selected direction of $L$.
        Consequently, for any $p\in \R^2$, $\tilde{p}_r+\ii \tilde{p}_i = (\psv_r+\ii \psv_i)+e^{\ii \theta}r h^{-1}$. Therefore, noting that $\im((\tilde{p}_r+\ii \tilde{p}_i)e^{-\ii\theta})^m = \im((\psv+\ii \psv)+rh^{-1})^m$, we conclude that $\tilde{\mathbb{A}}_0=\mathbb{B} \mathbb{A}_0$ with $\mathbb{B}:=( \binom{m}{n} (r/h)^{m-n})_{1\le m, n\le M+1}$, where $\binom{m}{n}:=0$ for $m<n$ and $\binom{m}{n}:=\frac{m!}{n!(m-n)!}$ for $m\ge n$. Because $\mathbb{B}$ is a lower triangular square matrix with unit diagonal, $\mathbb{B}$ is invertible. Due to $\vec{b}_0=0$ and $\tilde{\mathbb{A}}_0=\mathbb{B} \mathbb{A}_0$, we conclude that $\tilde{\mathbb{A}}_0\vec{c}_0=0$ is equivalent to $\mathbb{A}_0\vec{c}_0=0$, sharing the same solution space of $\vec{c}_0$. Because the augmented linear equations are independent of the choice of $\tbpt\in L$, we conclude that the solution $\vec{c}_0$ to $\tilde{\mathbb{A}}_0^* \vec{c}_0 = \vec{b}_0^*$ is independent of the choice of $\tbpt\in L$.
    \end{proof}

    The above result shows that the choice of the stencil $\SS_{\spt}$ and the property of the matrix $\mathbb{A}_0^*$ are only related to the local geometry of the grid $h\Z^2$, the region $\Omega$ and the tangent line $L_{\bpt}$ near $\spt$. To study the admissibility of the zeroth-order coefficients $\vec{c}_0=\{c_{p, 0}\}_{p\in \SS_{\spt}}$, we will not perform analysis for the specific stencil $\spt + h \SS_{\spt}$ and tangent line $L_{\bpt}$ at a boundary grid point $\spt \in \partial \Omega_h$; instead, we consider a point $\spt$ and a generic line $L$ tangent to $\partial \Omega$ satisfying
    \begin{equation}
        \label{eq: general spt L}
        \spt + [-h, h]^2 \cap \overline{\Omega} \neq \emptyset, \qquad
        \spt + [-h, h]^2 \cap L \neq \emptyset,
    \end{equation}
    and we construct the stencil $\spt + h \SS_{\spt}$ according to \Cref{sec: choice of boundary stencil}. The conditions in \eqref{eq: general spt L} are naturally satisfied under the specific construction of $\spt \in \partial \Omega_h$ and $L = L_{\bpt}$.

    The position of a directed line $L$ with direction angle $\theta_L \in (-\pi, \pi]$, relative to the point $\spt$, can be described with two parameters $\tau$ and $d$ as follows:
    \begin{equation}
        \label{eq: tau d}
        \tau =
        \begin{cases}
            \tan \theta_L, & k = 1, 2, 5, \\
            \tan (\theta_L - \frac{\pi}{4}), & k = 3, 4, 6,
        \end{cases}
        \quad \mbox{and} \quad d =
        \begin{cases}
            \frac{1}{h} |\overrightarrow{AC}|, & k = 1, 2, 5, \\
            \frac{1}{\sqrt{2}h} |\overrightarrow{AC}|, & k = 3, 4, 6.
        \end{cases}
    \end{equation}
    Here $1 \leq k \leq 6$ is the type of the stencil $\spt + h \SS_{\spt}$, $A = \spt$, and the point $C$ is shown in \Cref{fig: stencil types,fig: other inside point cases}. Under the assumption $(\spt + h\SS_{\spt}) \cap \hp{L} \subset \Omega$, where $\hp{L}$ is the open half plane to the left of $L$, we denote the parameter space of the pair $(\tau, d)$ for stencil type $k$ as $\mathcal{P}_k (0)$. We list the parameter space in the second column of \Cref{table: parameter space 1}. Note that $\mathcal{P}_k (0)$ does not depend on $h$, and in the set $\mathcal{P}_2(0)$, we purposefully included the case where $\spt + ph \notin \hp{L}$ for $p = (-1, -1)$ and $(0, -1)$. Lifting the assumption $(\spt + h\SS_{\spt}) \cap \hp{L} \subset \Omega$, we denote the parameter space as $\mathcal{P}_k (h)$. We aim to show that when $h$ is sufficiently small, then $\mathcal{P}_k (h)$ is ``close enough'' to $\mathcal{P}_k (0)$. If this is true, by verifying that the solution $\vec{c}_0$ to $\mathbb{A}_0^* \vec{c}_0=\vec{b}_0^*$ is admissible for all parameters $(\tau, d)$ in a set slightly larger than $\mathcal{P}_k(0)$, then $\vec{c}_0$ is admissible for all parameters $(\tau, d) \in \mathcal{P}_k(h)$ for sufficiently small $h$. In particular, for the specific grid $\Omega_h$ and the boundary stencils on it, the solution $\vec{c}_0$ to $\mathbb{A}_0^* \vec{c}_0=\vec{b}_0^*$ is admissible as in \Cref{def:admissible} if we set the grid size sufficiently small. The same argument holds for the existence, uniqueness and numerical stability of the solution $\vec{c}_0$.

    \begin{table}[h!]
        \begin{center}
        \begin{NiceTabular}{|c|c|c|}[cell-space-limits=4pt]
            \hline
            Stencil type & $\mathcal{P}_k (0)$ & $\mathcal{P}_k (\infty)$ \\ \hline\hline
             1& \Block{1-1}{$\tau \in [0, 1)$ \\ $d \in \left( 0, \max \{\tau, 1 - \tau\} \right]$} & \Block{1-1}{$\tau \in \R$ \\ $d \in (0, \infty)$} \\ \hline
             2& \Block{1-1}{$\tau \in (-1, 1)$ \\ $d \in (|\tau|, 1]$} & \Block{1-1}{$\tau \in \R$ \\ $d \in (|\tau|, \infty)$} \\ \hline
              3 & \Block{1-1}{$\tau \in (-\frac{1}{3}, \frac{1}{3})$ \\ $d \in \left( |\tau|, \frac{1}{2} (1 - |\tau|) \right]$} & \Block{1-1}{$\tau \in \R$ \\ $d \in (|\tau|, \infty)$} \\  \hline
             4& \Block{1-1}{$\tau \in (-1, 1)$ \\ $d \in \left( \frac{1}{2} (1 + |\tau|), 1 \right]$} & \Block{1-1}{$\tau \in (-1, 1)$ \\ $d \in \left( \frac{1}{2} (1 + |\tau|), 1 \right]$}  \\ \hline
             5 & $\varnothing$ & \Block{1-1}{$\tau \in \R$ \\ $d \in (0, \infty)$} \\ \hline
             6 & $\varnothing$ & \Block{1-1}{$\tau \in \R$ \\ $d \in (0, \infty)$} \\ \hline
        \end{NiceTabular}
        \end{center}
        \vspace{-9pt}
        \caption{\footnotesize Parameter space $\mathcal{P}_k(0)$ under the assumption $(\spt + h\SS_{\spt}) \cap \hp{L} \subset \Omega$, and the largest parameter space $\mathcal{P}_k(\infty) \supseteq \bigcup_{h > 0} \mathcal{P}_k (h)$. The definitions of $\tau$ and $d$ are given in \eqref{eq: tau d}.}
        \label{table: parameter space 1}
    \end{table}

    To begin with, we denote $\overline{\mathcal{P}_k(0)}$ to be the usual closure of $\mathcal{P}_k(0)$ for $1 \leq k \leq 4$, and to be the set $\{\tau = d = 0\}$ for $k = 5, 6$. We also take a set $\mathcal{P}_k (\infty) \supseteq \bigcup_{h > 0} \mathcal{P}_k (h)$ and present it in the third column of \Cref{table: parameter space 1}. According to the last condition in \eqref{eq: general spt L}, $\tau$ does not take $\infty$ and the infimum of $d$ is the same as in $\mathcal{P}_k (0)$. Otherwise, the stencil will not follow the designated stencil type as certain grid points fall outside of $\hp{L}$. Due to the same condition, for stencil type 4 we have $\spt + (1, -1)h \notin \hp{L}$. We can therefore set $\mathcal{P}_4 (\infty) = \mathcal{P}_4 (0)$.

    We adopt a topological approach. Observe that the set of directed lines forms a topological manifold $\mathcal{M}$ homeomorphic to $\mathbb{S} \times \R$, where $\mathbb{S}$ is the unit circle. This manifold has an atlas $\{ (\phi_{kj}, \mathcal{M}_{kj}, \R^2) \}_{1 \leq k \leq 6, \, 1 \leq j \leq J_k}$ given by the definitions below.
    \begin{itemize}
        \item $\phi_k$ is the map $L \mapsto (\tau, hd)$ given by \eqref{eq: tau d} for stencil type $k$, except that point $A = (0, 0)$.

        \item Through rotation and flipping, one can transform one type of stencil into another type not listed in \Cref{fig: inside grid points cases,fig: other inside point cases}. For stencil type $k$, $J_k$ is defined to be the number of different types of stencils through such transformation.

        \item A stencil transformed from type $k$ can be obtained from applying a linear transform $T_{kj}$ to the original type-$k$ stencil. Now we set the map $\phi_{kj}$ to be $\phi_{k} \circ T_{kj}^{-1}$. We also set $T_{k1} = I_2$, the identity map on $\R^2$.

        \item $\mathcal{M}_{kj} := \phi_{kj}^{-1} (\R^2)$ is an open subset of $\mathcal{M}$.
    \end{itemize}

    In addition, for any $L \in \mathcal{M}$, define
    \begin{equation}
        \label{eq: psi_L}
        \psi_L: \R^2 \to \R, \ (x, y) \mapsto ((x, y) - \bpt_L) \cdot (\sin \theta_L, -\cos \theta_L),
    \end{equation}
    where the point $\bpt_L \in L$. The function $\psi_L$ represents the coordinate of a point $(x, y)$ normal to the direction of $L$, and its definition does not depend on the choice of $\bpt_L$. Moreover, $\psi_L(p) < 0$ if and only if $p \in \hp{L}$. Now, $\mathcal{M}$ can be embedded into $\R^9$, given by the mapping
    \begin{align*}
        \Psi: L \in \mathcal{M} \mapsto (\psi_L(p))_{p \in \SS} \in \R^9,
    \end{align*}
    where $\SS  = [-1, 1]^2 \cap \Z^2$ with its usual ordering given in \eqref{S:order}. This embedding $\Psi$ is used exactly as the criteria to classify the cases of the grid points within $\hp{L}$. For $k \neq 2$, $\Psi \circ \phi_{kj} (\overline{\mathcal{P}_k(0)})$ is merely the intersection of $\Phi(\mathcal{M})$ and the product of several intervals of $\R_\leq := (-\infty, 0]$ or $\R_\geq := [0, \infty)$ in a certain order. For example,
    \begin{align}
        \label{eq: Psi circ phi_11 P_1(0)}
        & \Psi \circ \phi_{11} (\overline{\mathcal{P}_1 (0)}) = \Phi(\mathcal{M}) \cap \prod_{\ell \in \{1, 4, 7, 8\}} \R_\geq^{(\ell)} \times \prod_{\ell \in \{2, 3, 5, 6, 9\}} \R_\leq^{(\ell)}, \\
        & \Psi \circ \phi_{61} (\overline{\mathcal{P}_6 (0)}) = \Phi(\mathcal{M}) \cap \prod_{\ell \in \{1, 4, 7, 8, 9\}} \R_\geq^{(\ell)} \times \prod_{\ell \in \{2, 3, 5, 6\}} \R_\leq^{(\ell)}, \nonumber
    \end{align}
    where the superscript $(\ell)$ indicates that the $\ell$-th component of $\Psi(L)$ belong to that interval. When $k = 2$, the set $\Psi \circ \phi_{21} (\overline{\mathcal{P}_2 (0)})$ is given by
    \begin{align*}
        \Phi(\mathcal{M}) \cap \prod_{\ell \in \{2, 3, 5, 6, 8, 9\}} \R_\leq^{(\ell)} \times \left( \R_\geq^{(1)} \times \R_\geq^{(4)} \times \R_\geq^{(7)} \cup \R_\leq^{(1)} \times \R_\geq^{(4)} \times \R_\geq^{(7)} \cup \R_\geq^{(1)} \times \R_\geq^{(4)} \times \R_\leq^{(7)} \right).
    \end{align*}

    Now we formulate and prove the result that $\mathcal{P}_k(h)$ approaches $\mathcal{P}_k(0)$.

    \begin{lemma}
        \label{lem: P_k(h) approaches P_k(0)}
        Suppose the domain $\Omega \subseteq \R^2$ has $C^1$ boundary. Define the sets $\mathcal{P}_k(0)$, $\mathcal{P}_k (h)$ and $\mathcal{P}_k (\infty)$ as above. Then for any $1 \leq k \leq 6$ and any open set $\mathcal{P} \supseteq \overline{\mathcal{P}_k(0)}$, there exists $h_* = \bo_{\gx, \gy} (1)$ such that $\mathcal{P}_k (h) \subseteq \mathcal{P} \cap \mathcal{P}_k (\infty)$ for all $0 < h < h_*$.
    \end{lemma}

    \begin{proof}
        Fix $1 \leq k \leq 6$. It is enough to prove that $\mathcal{P}_k (h) \subseteq \mathcal{P}$ when $h$ is small enough. This is equivalent to
        \begin{equation}
            \label{eq: equivalence P_k(h) approaches P_k(0)}
            \Psi(\mathcal{M}) \backslash \Psi \circ \phi_{k1} (\mathcal{P}) \subseteq \Psi(\mathcal{M}) \backslash \Psi \circ \phi_{k1} (\mathcal{P}_k (h)),
        \end{equation}
        when $h$ is small enough.

        Since $\mathcal{P}$ is an open set containing $\overline{\mathcal{P}_k(0)}$, the boundaries of these sets have a positive $L^\infty$ distance. It follows that the boundaries of $\Psi \circ \phi_{k1} (\mathcal{P})$ and $\Psi \circ \phi_{k1} (\overline{\mathcal{P}_k(0)})$ have a positive $L^\infty$ distance as well, which means that
        \begin{equation*}
            \Psi \circ \phi_{k1} (\mathcal{P}) \supseteq \left\{ \xi \in \Psi(\mathcal{M}): \|\xi - \eta\|_{\infty} < \epsilon_0 \ \mbox{for some} \ \eta \in \Psi \circ \phi_{k1} (\overline{\mathcal{P}_k(0)}) \right\}.
        \end{equation*}
        for some $\epsilon_0 > 0$. This implies
        \begin{equation*}
            \Psi(\mathcal{M}_0) \backslash \Psi \circ \phi_{k1} (\mathcal{P}) \subseteq \left\{ \xi \in \Psi(\mathcal{M}_0): \|\xi - \eta\|_{\infty} \geq \epsilon_0 \ \mbox{for all} \ \eta \in \Psi \circ \phi_{k1} (\overline{\mathcal{P}_k(0)}) \right\},
        \end{equation*}
        where $\mathcal{M}_0$ is the set of all directed lines $L$ so that $L \cap [-1, 1]^2 \neq \emptyset$.
        Hence, to prove \eqref{eq: equivalence P_k(h) approaches P_k(0)}, we only need to prove the following statement: given $\xi \in \Psi(\mathcal{M}_0)$ such that $\|\xi - \eta\|_{\infty} \geq \epsilon_0 \ \mbox{for all} \ \eta \in \Psi \circ \phi_{k1} (\overline{\mathcal{P}_k(0)})$, we have $\xi \notin \Psi \circ \phi_{k1} (\mathcal{P}_k (h))$.

        For any tangent line $L^*$ on $\partial \Omega$, let $\bpt_L$ be the tangent point of $L^*$. Let $\sv \in [-1, 1]^2$ and set $\spt = \bpt_L + h\sv$, then $\spt$ and $L^*$ satisfy the conditions in \eqref{eq: general spt L}. If we fix $L^*$ and $\sv$, then $L := (L^* - \spt) / h \in \mathcal{M}$ is independent of $h > 0$. The set of lines $L$ obtained from all possible tangent lines $L^*$ and $\sv$ is identical to $\mathcal{M}_0$.

        Let $\xi \in \Psi(\mathcal{M}_0)$ such that $\|\xi - \eta\|_{\infty} \geq \epsilon_0$ for all $\eta \in \Psi \circ \phi_{k1} (\overline{\mathcal{P}_k(0)})$. From the above discussion, there exists a tangent line $L^*$ and $\sv \in [-1, 1]^2$ such that $\xi = \Psi(L)$. We can find $h_1 = \bo_{\gx, \gy} (1) > 0$, so that $(\spt + [-h, h]^2) \cap \partial \Omega \neq \emptyset$ consists of a single segment of curve when $0 < h < h_1$. In other words,
        \begin{equation*}
            (\gx, \gy)^{-1} \left( (\spt + [-h, h]^2) \cap \partial \Omega \right) = [t_1(h), t_2(h)]
        \end{equation*}
        for some $t_1(h) \leq t_2(h)$. Since $(\gx, \gy)$ is a $C^1$ curve, we have $t_2(h) - t_1(h) = \bo_{\gx, \gy}(h)$. By Taylor expansion at the tangent point $\bpt_L$, we can obtain $\psi_{L^*} (x', y') = \so_{\gx, \gy}(h)$ for any point $(x', y') = (\gx(t'), \gy(t'))$ on this segment of curve (here $\so$ follows the same convention as $\bo$). Hence, there exists $h_2 = \bo_{\gx, \gy} (1) \in (0, h_1)$, so that
        \begin{equation}
            \label{eq: close to tangent line}
            \max_{t' \in [t_1(h), t_2(h)]} \left| \psi_{L^*} (\gx(t'), \gy(t')) \right| \leq \frac{1}{2} \epsilon_0 h, \quad \forall 0 < h < h_2.
        \end{equation}
        This implies that, when $0 < h < h_2$, any point $(x', y') \in \spt + [-h, h]^2$ such that $\psi_{L^*} (x', y') > \frac{1}{2} \epsilon_0 h$ is outside the region $\Omega$.

        In the remaining proof, we suppose $k = 1$. The same method applies for all $1 \leq k \leq 6$. Since $\xi = (\xi_\ell)_{1 \leq \ell \leq 9} = \Psi(L)$ satisfies $\|\xi - \eta\|_{\infty} \geq \epsilon_0$ for all $\eta \in \Psi \circ \phi_{11} (\overline{\mathcal{P}_1(0)})$, from \cref{eq: Psi circ phi_11 P_1(0)} we know that there exists an index $1 \leq \ell_0 \leq 9$ so that
        \begin{equation*}
            \xi_{\ell_0}
            \begin{cases}
                \geq  \epsilon_0, & \mbox{if} \ \ell_0 \in \{1, 4, 7, 8\}, \\
                \leq -\epsilon_0, & \mbox{if} \ \ell_0 \in \{2, 3, 5, 6, 9\}.
            \end{cases}
        \end{equation*}
        From the definition of $\Psi$, we know that $\xi_{\ell_0} = \psi_L (p)$ for the $\ell_0$-th element $p \in \SS$. Since $L = (L^* - \spt) / h$, we obtain $\psi_{L^*} (\spt + ph) = h \psi_L (p)$. It follows that $\mathrm{sgn} (\xi_{\ell_0}) \psi_{L^*} (\spt + ph) \geq \epsilon_0 h$.

        If $\ell_0 \in \{2, 3, 5, 6, 9\}$, then the stencil point $\spt + ph$ is in $\hp{L^*}$, which shows that this stencil will not be of type 1. If $\ell_0 \in \{1, 4, 7, 8\}$, then $\psi_{L^*} (\spt + ph) \geq \epsilon_0 h$. Together with the fact that $\spt + ph \in \spt + [-h, h]^2$, this implies the grid point $\spt + ph$ is outside the region $\Omega$ when $0 < h < h_2$. In this case, the stencil is not of type 1 either. Therefore, $\xi \notin \Psi \circ \phi_{11} (\mathcal{P}_1 (h))$. This completes the proof of all claims.
    \end{proof}

    Finally, we have found a set $\mathcal{P}$ which is the intersection of $\mathcal{P}_k(\infty)$ and an open set containing $\overline{\mathcal{P}_k(0)}$, and verified that there exists a unique admissible solution $\vec{c}_0$ to $\mathbb{A}_0^*\vec{c}_0=\vec{b}_0^*$ for all parameters $(\tau, d) \in \mathcal{P}$. The set $\mathcal{P}$ is listed in the second column of \Cref{table: parameter space 2}. The quantity $\mu_c$ in item (iii) of \Cref{prop: boundary c} is numerically calculated from $\inf_{(\tau, d) \in \mathcal{P}} \sum_{p \in \SS} c_{p, 0}$, and is shown in the third column of \Cref{table: parameter space 2}. For stability, we verified that the matrix $\mathbb{A}_0^*$ is well-conditioned for each stencil type. In the last column of \Cref{table: parameter space 2}, we present the supremum $M_\kappa$ of the $L^\infty$ condition number $\kappa (\mathbb{A}_0^*)$ over all parameters in $\mathcal{P}$. According to \Cref{lem: P_k(h) approaches P_k(0)} and the discussion before it, the unique solution $\vec{c}_0$ to $\mathbb{A}_0^*\vec{c}_0=\vec{b}_0^*$ is always stable and admissible if we set the grid size $h$ sufficiently small.

    \begin{table}[h!]
        \begin{center}
        \begin{NiceTabular}{|c|c|c|c|}[cell-space-limits=4pt]
            \hline
            Stencil type $k$ & $\mathcal{P}$ & $\mu_c$ & $M_\kappa$ \\ \hline\hline
            1 & \Block{1-1}{$\tau \in (-\frac{1}{6}, \frac{7}{6})$ \\ $d \in \left( 0, \max \{ \frac{1}{6} + \tau, \frac{7}{6} - \tau \} \right)$} & 0.591 & 89.3786 \\ \hline
            2 & \Block{1-1}{$\tau \in (-\frac{11}{10}, \frac{11}{10})$ \\ $d \in (|\tau|, \frac{11}{10})$} & 0.247 & 308.050 \\ \hline
            3 & \Block{1-1}{$\tau \in (-\frac{1}{2}, \frac{1}{2})$ \\ $d \in \left( |\tau|, \frac{1}{2} (\frac{3}{2} - |\tau|) \right)$} & 0.355 & 491.000 \\  \hline
            4 & \Block{1-1}{$\tau \in (-1, 1)$ \\ $d \in \left( \frac{1}{2} (1 + |\tau|), 1 \right]$} & 0.050 & 324.498 \\ \hline
            5 & \Block{1-1}{$\tau \in (-\frac{1}{6}, \frac{1}{6})$ \\ $d \in (0, \frac{1}{6} - \tau)$} & 0.875 & 39.4634 \\ \hline
            6 & \Block{1-1}{$\tau \in (-\frac{1}{10}, \frac{1}{10})$ \\ $d \in (0, \frac{1}{10} - \tau)$} & 0.852 & 600.507 \\ \hline
        \end{NiceTabular}
        \end{center}
        \vspace{-9pt}
        \caption{\footnotesize Second column: parameter space $\mathcal{P}$ in which the zeroth-order solutions $\vec{c}_0=\{c_{p, 0}\}_{p\in \SS_{\spt}}$ are admissible as described in \Cref{def:admissible}; Third column: $\mu_c := \inf_{(\tau, d) \in \mathcal{P}} \sum_{p \in \SS} c_{p, 0}$; Fourth column: $M_\kappa := \sup_{(\tau, d) \in \mathcal{P}} \kappa (\mathbb{A}_0^*)$. The matrix $\mathbb{A}_0^*$ is defined in \eqref{extra:eqs} and the vectors $\mathbb{A}^*_{0;k}$ and $\vec{b}^*_0$ in \eqref{extra:eqs} for $k=1,\ldots, \#\SS_{\spt}-5$  are determined by the extra constraints in \Cref{table: stencil information,table: other stencil information}.}
        \label{table: parameter space 2}
    \end{table}

    \section{Second and Fourth-order Schemes at Boundary Grid Points}
    \label{app: order 2 4}

    In this section, we briefly talk about the essential changes to the proposed sixth-order FDM scheme at boundary grid points in order to get a second or fourth-order scheme.

    \subsection{Second-order FDM scheme}
    We set $\SS = \{ \spt, \spt + ph \}$ for some $p \in \Z^2$ with $p$ near $(0,0)$ and $\spt, \spt+ph\in \Omega$. In this case, all solutions to \cref{eq: boundary constraint of c} are given by
    \begin{equation*}
        c_{p, 0} = -\frac{\sv_r \sin \theta - \sv_i \cos \theta}{\psv_r \sin \theta - \psv_i \cos \theta} c_{(0, 0), 0},
    \end{equation*}
    where $\sv$ and $\psv$ are defined in \cref{vec:shift} with $\sv = (\sv_r, \sv_i)$ and $\psv = (\psv_r, \psv_i)$. Using the same normalization $c_{(0, 0), 0} = 1$ as before, we get
    \begin{equation}
        \label{eq: 2nd order boundary c}
        C_{(0, 0)}(h) = 1,
        \quad C_p(h) = -\frac{\sv_r \sin \theta - \sv_i \cos \theta}{\psv_r \sin \theta - \psv_i \cos \theta}.
    \end{equation}
    In order to let the coefficients satisfy the properties in \Cref{prop: boundary c}, we only need to set a threshold $\mu_c>0$, and then for each boundary stencil point $\spt$, we look for a desired point $\spt + ph$ satisfying
    \begin{equation}
        \label{eq: 2nd order condition}
        0 \leq -C_p(h) = \frac{\sv_r \sin \theta - \sv_i \cos \theta}{\psv_r \sin \theta - \psv_i \cos \theta} \leq 1 - \mu_c.
    \end{equation}

    The condition \eqref{eq: 2nd order condition} can be very easily satisfied. To see this, we adopt the function $\psi_L$ defined in \cref{eq: psi_L}, where $L$ is the tangent line. Using $\bpt_L = \bpt$, we see that \cref{eq: 2nd order condition} is equivalent to $0 \leq \psi_L(\spt) / \psi_L(\spt + ph) \leq 1 - \mu_c$. Since $\psi_L$ represents the coordinate of a point perpendicular to the tangent line, and the point $\spt$ is always inside the tangent line, the above condition just means that the perpendicular coordinate of $\spt + ph$ should be at least $\frac{1}{1 - \mu_c}$ times that of $\spt$. Such a point can be found at ease.

    In practice, one only needs to iterate through several grid points adjacent to $\spt$, calculate the coefficients according to \cref{eq: 2nd order boundary c}, and verify directly whether condition \eqref{eq: 2nd order condition} holds.

    \subsection{Fourth-order FDM scheme}
    Same as in the sixth-order scheme, we adopt 6 types of different stencils according to the points inside the tangent line and the boundary. We describe the choice of the stencil and the extra constraints in \Cref{table: stencil information order 4}, where the definition of the parameters $(\tau, d)$ are taken in the same way as \cref{eq: tau d}. We still use the parameter spaces $\mathcal{P}_k (0)$, $\mathcal{P}_k (\infty)$ and $\mathcal{P}$ in \Cref{table: parameter space 1,table: parameter space 2} from the sixth-order scheme. Under this parameter space, we present \Cref{table: parameter space order 4} as an analog of \Cref{table: parameter space 2} for the fourth-order scheme.

    \begin{table}[h!]
        \begin{center}
        \begin{NiceTabular}{|c|c|c|c|c|}[cell-space-limits=4pt]
            \hline
            Stencil type & Case & $\#\SS_{\spt}$ & $\SS_{\spt}$ & Extra constraints \\ \hline\hline
            1 & I    & 5 & \Block{1-1}{$\{ (0, 0), (-1, 0),$ \\ $(-1, 1), (0, 1), (1, 1) \}$} & $\vec{c}_0(2) - \vec{c}_0(5) = 0$ \\ \hline
            \Block{3-1}{2} & II  & \Block{3-1}{6} & \Block{3-1}{$\{ (0, 0), (-1, 0), (1, 0),$ \\ $(-1, 1), (0, 1), (1, 1) \}$} & \Block{3-1}{$\vec{c}_0(2) - \vec{c}_0(6) = -\frac{|\overrightarrow{AC}|}{5 h}$ \\ $\vec{c}_0(3) - \vec{c}_0(4) = -\frac{|\overrightarrow{AC}|}{5 h}$} \\ \cline{2-2}
              & IV   &   &  & \\ \cline{2-2}
              & VIII &   &  & \\ \hline
            3 & III  & 6 & \Block{1-1}{$\{ (0, 0), (-1, -1), (1, 1),$ \\ $(-1, 0), (0, 1), (-1, 1) \}$} & \Block{1-1}{$\vec{c}_0(2) = -\frac{1}{20}$ \\ $\vec{c}_0(3) = -\frac{1}{20}$} \\ \hline
            4 & V    & 6 & \Block{1-1}{$\{ (0, 0), (0, -1), (1, 0),$ \\ $(-1, 0), (0, 1), (-1, 1) \}$} & \Block{1-1}{$\vec{c}_0(2) - \vec{c}_0(5) = 0$ \\ $\vec{c}_0(3) - \vec{c}_0(4) = 0$} \\ \hline
            5 & VI   & 4 & \Block{1-1}{$\{ (0, 0), (0, 2),$ \\ $(-1, 1), (1, 1) \}$} & N/A \\ \hline
            6 & VII  & 4 & \Block{1-1}{$\{ (0, 0), (-1, 0),$ \\ $(0, 1), (-1, 1) \}$} & N/A \\ \hline
        \end{NiceTabular}
        \end{center}
        \vspace{-9pt}
        \caption{\footnotesize The stencil types and the extra equations in \eqref{extra:eqs} for the fourth-order FDM scheme. Cases I -- VIII are the same as the sixth-order scheme, and the number $|\overrightarrow{AC}|$ above is the distance between points $A$ and $C$ in \Cref{fig: stencil types}.}
        \label{table: stencil information order 4}
    \end{table}

    \begin{table}[h!]
        \begin{center}
        \begin{NiceTabular}{|c|c|c|}[cell-space-limits=4pt]
            \hline
            Stencil type $k$ & $\mu_c$ & $M_\kappa$ \\ \hline\hline
            1 & 0.564 & 24.1369 \\ \hline
            2 & 0.239 & 61.9390 \\ \hline
            3 & 0.366 & 43.5262 \\ \hline
            4 & 0     & 38.0575 \\ \hline
            5 & 0.874 & 9.33333 \\ \hline
            6 & 0.852 & 14.9249 \\ \hline
        \end{NiceTabular}
        \end{center}
        \vspace{-9pt}
        \caption{\footnotesize The constants $\mu_c := \inf_{(\tau, d) \in \mathcal{P}} \sum_{p \in \SS} c_{p, 0}$ and the condition number $M_\kappa := \sup_{(\tau, d) \in \mathcal{P}}$ $\kappa (\mathbb{A}_0^*)$ for the fourth-order FDM scheme. The parameter spaces $\mathcal{P}$ are the same as in the sixth-order scheme.}
        \label{table: parameter space order 4}
    \end{table}

    Readers should be aware that $\mu_c = 0$ for stencil type 4, which violates the admissibility condition $\mu_c > 0$. Indeed, $(1, -\frac{1}{4}, -\frac{1}{4}, -\frac{1}{4}, -\frac{1}{4}, 0)$ is the unique zeroth-order stencil coefficients under the designated stencil. These coefficients satisfy the admissibility conditions (i), (ii) and $\sum_{p \in \SS} c_{p, 0} = 0$, which is characteristic of the interior stencil coefficients (see \Cref{prop: interior c}). In this situation, we modify the higher-order stencil coefficients using \cref{q:modcoeff} instead of \eqref{eq: lambda_m boundary}. To prove the fourth-order convergence, we only need to treat type-4 boundary stencil as an interior stencil. Besides, it is impossible if all boundary stencils are of type 4. We omit the detailed discussion.

\bibliographystyle{unsrtnat}
\bibliography{ref}

\begin{thebibliography}{31}
\providecommand{\natexlab}[1]{#1}
\providecommand{\url}[1]{\texttt{#1}}
\expandafter\ifx\csname urlstyle\endcsname\relax
  \providecommand{\doi}[1]{doi: #1}\else
  \providecommand{\doi}{doi: \begingroup \urlstyle{rm}\Url}\fi

\bibitem[Jensen(1972)]{jensen1972finite}
P.~S. Jensen.
\newblock Finite difference techniques for variable grids.
\newblock \emph{Computers \& Structures}, 2:\penalty0 17--29, 1972.

\bibitem[Li and Pan(2023)]{li2023high}
Z.~Li and K.~Pan.
\newblock High order compact schemes for flux type {BC}s.
\newblock \emph{SIAM J. Sci. Comput.}, 45:\penalty0 A646--A674, 2023.

\bibitem[Settle et~al.(2013)Settle, Douglas, Kim, and Sheen]{settle2013derivation}
S.~O. Settle, C.~C. Douglas, I.~Kim, and D.~Sheen.
\newblock On the derivation of highest-order compact finite difference schemes for the one- and two-dimensional poisson equation with dirichlet boundary conditions.
\newblock \emph{SIAM J. Numer. Anal.}, 51:\penalty0 2470--2490, 2013.

\bibitem[Wang and Zhang(2009)]{wang2009sixth}
Y.~Wang and J.~Zhang.
\newblock Sixth order compact scheme combined with multigrid method and extrapolation technique for {2D} poisson equation.
\newblock \emph{J. Comput. Phys.}, 228:\penalty0 137--146, 2009.

\bibitem[Zhai et~al.(2013)Zhai, Feng, and He]{zhai2013family}
S.~Zhai, X.~Feng, and Y.~He.
\newblock A family of fourth-order and sixth-order compact difference schemes for the three-dimensional {P}oisson equation.
\newblock \emph{J. Sci. Comput.}, 54:\penalty0 97--120, 2013.

\bibitem[Zhai et~al.(2014)Zhai, Feng, and He]{zhai2014new}
S.~Zhai, X.~Feng, and Y.~He.
\newblock A new method to deduce high-order compact difference schemes for two-dimensional {P}oisson equation.
\newblock \emph{Appl. Math. Comput.}, 230:\penalty0 9--26, 2014.

\bibitem[Feng et~al.(2021)Feng, Han, and Minev]{feng2021sixth}
Q.~Feng, B.~Han, and P.~Minev.
\newblock Sixth order compact finite difference schemes for poisson interface problems with singular sources.
\newblock \emph{Comput. Math. Appl.}, 99:\penalty0 2--25, 2021.

\bibitem[Feng and Zhao(2020)]{feng2020fft}
H.~Feng and S.~Zhao.
\newblock {FFT}-based high order central difference schemes for three-dimensional {P}oisson's equation with various types of boundary conditions.
\newblock \emph{Journal of Computational Physics}, 410:\penalty0 109391, 2020.

\bibitem[Ma and Ge(2020)]{ma2020high}
T.~Ma and Y.~Ge.
\newblock High-order blended compact difference schemes for the 3{D} elliptic partial differential equation with mixed derivatives and variable coefficients.
\newblock \emph{Adv. Difference Equ.}, 2020, 2020.
\newblock Paper No. 525. 30 pp.

\bibitem[Wang et~al.(2014)Wang, Guo, and Wu]{wang2014fourth}
Y.-M. Wang, B.-Y. Guo, and W.-J. Wu.
\newblock Fourth-order compact finite difference methods and monotone iterative algorithms for semilinear elliptic boundary value problems.
\newblock \emph{Comput. Math. Appl.}, 68:\penalty0 1671--1688, 2014.

\bibitem[Shi et~al.(2021)Shi, Xie, Liang, and Fu]{shi2021high}
Y.~Shi, S.~Xie, D.~Liang, and K.~Fu.
\newblock High order compact block-centered finite difference schemes for elliptic and parabolic problems.
\newblock \emph{J. Sci. Comput.}, 87:\penalty0 1--26, 2021.

\bibitem[Feng et~al.(2022)Feng, Han, and Minev]{FHM22}
Q.~Feng, B.~Han, and P.~Minev.
\newblock A high order compact finite difference scheme for elliptic interface problems with discontinuous and high-contrast coefficients.
\newblock \emph{Appl. Math. Comput.}, 431, 2022.
\newblock Paper No. 12734. 24 pp.

\bibitem[Feng et~al.(2024)Feng, Han, and Minev]{feng2024sixth}
Q.~Feng, B.~Han, and P.~Minev.
\newblock Sixth-order hybrid finite difference methods for elliptic interface problems with mixed boundary conditions.
\newblock \emph{J. Comput. Phys.}, 497, 2024.
\newblock Paper No. 112635. 32 pp.

\bibitem[Shortley and Weller(1938)]{sw38}
G.~H. Shortley and R.~Weller.
\newblock The numerical solution of {L}aplace's equation.
\newblock \emph{J. Appl. Phys.}, 9:\penalty0 334--348, 1938.

\bibitem[Bramble and Hubbard(1964)]{bramble1964new}
J.~H. Bramble and B.~E. Hubbard.
\newblock New monotone type approximations for elliptic problems.
\newblock \emph{Math. Comp.}, 18:\penalty0 349--367, 1964.

\bibitem[Price(1968)]{price1968monotone}
H.~S. Price.
\newblock Monotone and oscillation matrices applied to finite difference approximations.
\newblock \emph{Math. Comp.}, 22:\penalty0 489--516, 1968.

\bibitem[Esmaeilzadeh and Barron(2022)]{eb22}
M.~Esmaeilzadeh and R.~M. Barron.
\newblock Numerical solution of partial differential equations in arbitrary shaped domains using cartesian cut-stencil finite difference method. {Part II}: Higher-order schemes.
\newblock \emph{Numer. Math. Theory, Methods Appl.}, 15:\penalty0 819--850, 2022.

\bibitem[Pan et~al.(2021)Pan, He, and Li]{pan2021high}
K.~Pan, D.~He, and Z.~Li.
\newblock A high order compact {FD} framework for elliptic bvps involving singular sources, interfaces, and irregular domains.
\newblock \emph{J. Sci. Comput.}, 88, 2021.
\newblock Paper No. 67. 25 pp.

\bibitem[Ren et~al.(2022)Ren, Feng, and Zhao]{RFZ22}
Y.~Ren, H.~Feng, and S.~Zhao.
\newblock A {FFT} accelerated high order finite difference method for elliptic boundary value problems over irregular domains.
\newblock \emph{J. Comput. Phys.}, 448, 2022.
\newblock Paper No. 110762. 24 pp.

\bibitem[Li et~al.(2025)Li, Zhao, Pentecost, Ren, and Guan]{Li2025spatially}
C.~Li, S.~Zhao, B.~Pentecost, Y.~Ren, and Z.~Guan.
\newblock A spatially fourth-order {C}artesian grid method for fast solutions of elliptic and parabolic problems on irregular domains with sharply curved boundaries.
\newblock \emph{Journal of Scientific Computing}, 103\penalty0 (94), 2025.

\bibitem[Samarskii and Fryazinov(1971)]{samarskii1971finite}
A.~A. Samarskii and I.~V. Fryazinov.
\newblock On finite-difference schemes for solving the dirichlet problem for an elliptic equation with variable coefficients in an arbitrary region.
\newblock \emph{USSR Comput. Math. Math. Phys.}, 11:\penalty0 109--139, 1971.

\bibitem[Ito et~al.(2005)Ito, Li, and Kyei]{ilk05}
K.~Ito, Z.~Li, and Y.~Kyei.
\newblock Higher-order, cartesian grid based finite difference schemes for elliptic equations on irregular domains.
\newblock \emph{SIAM J. Sci. Comput.}, 27:\penalty0 346--367, 2005.

\bibitem[Gibou and Fedkiw(2005)]{gf2005}
F.~Gibou and R.~Fedkiw.
\newblock A fourth order accurate discretization for the laplace and heat equations on arbitrary domains, with applications to the stefan problem.
\newblock \emph{J. Comput. Phys.}, 202\penalty0 (2):\penalty0 577--601, 2005.

\bibitem[Gibou et~al.(2002)Gibou, Fedkiw, Cheng, and Kang]{gfck2002}
F.~Gibou, R.~Fedkiw, L.~T. Cheng, and M.~Kang.
\newblock A second-order-accurate symmetric discretization of the {P}oisson equation on irregular domains.
\newblock \emph{J. Comput. Phys.}, 176\penalty0 (1):\penalty0 205--227, 2002.

\bibitem[Clain et~al.(2021)Clain, Lopes, and Pereira]{clp21}
S.~Clain, D.~Lopes, and R.~M. Pereira.
\newblock Very high-order cartesian-grid finite difference method on arbitrary geometries.
\newblock \emph{J. Comput. Phys.}, 434, 2021.
\newblock Paper No. 110217. 28 pp.

\bibitem[Varga(1966)]{varga1966discrete}
R.~S. Varga.
\newblock On a discrete maximum principle.
\newblock \emph{SIAM J. Numer. Anal.}, 3:\penalty0 355--359, 1966.

\bibitem[Shivakumar and Chew(1974)]{shivakumar1974sufficient}
P.~N. Shivakumar and K.~H. Chew.
\newblock A sufficient condition for nonvanishing of determinants.
\newblock \emph{Proc. Amer. Math. Soc.}, 43:\penalty0 63--66, 1974.

\bibitem[Plemmons(1977)]{plemmons1977m}
R.~J. Plemmons.
\newblock $m$-matrix characterizations. {I}. nonsingular $m$-matrices.
\newblock \emph{Linear Algebra Appl.}, 18:\penalty0 175--188, 1977.

\bibitem[Li and Zhang(2020)]{li2020monotonicity}
H.~Li and X.~Zhang.
\newblock On the monotonicity and discrete maximum principle of the finite difference implementation of {$C^0$ - $Q^2$} finite element method.
\newblock \emph{Numer. Math.}, 145:\penalty0 437--472, 2020.

\bibitem[Gilbarg and Trudinger(2001)]{gilbarg1977elliptic}
D.~Gilbarg and N.~S. Trudinger.
\newblock \emph{Elliptic partial differential equations of second order}.
\newblock Classics Math., Springer-Verlag, Berlin, 2001.

\bibitem[Levin(1998)]{levin1998approximation}
D.~Levin.
\newblock The approximation power of moving least-squares.
\newblock \emph{Math. Comp.}, 67:\penalty0 1517--1531, 1998.

\end{thebibliography}

\iffalse

\fi

\end{document}